\documentclass{article}

\usepackage{amsmath}
\usepackage{amsfonts}
\usepackage{amscd}
\usepackage{amssymb}

\input xypic

\newtheorem{defn}{Definition}[section]
\newtheorem{thm}[defn]{Theorem}
\newtheorem{prop}[defn]{Proposition}
\newtheorem{lemma}[defn]{Lemma}
\newtheorem{eg}[defn]{Example}
\newtheorem{claim}[defn]{Claim}
\newtheorem{cor}[defn]{Corollary}

\newtheorem{schol}[defn]{Scholium}

\newcommand{\lm}{\ensuremath{\longrightarrow}}

\DeclareMathOperator{\Rhom}{\mbox{RHom}}
\DeclareMathOperator{\Hom}{\mbox{Hom}}
\DeclareMathOperator{\shom}{\ensuremath{\mathcal{H}\mathit{om}}}
\DeclareMathOperator{\stor}{\ensuremath{\mathcal{T}\mathit{or}}}

\DeclareMathOperator{\End}{\mbox{End}}

\DeclareMathOperator{\Ext}{\mbox{Ext}}
\DeclareMathOperator{\Tor}{\mbox{Tor}}
\DeclareMathOperator{\intotimes}{\ensuremath{\underline{\otimes}_{\cala}}}
\DeclareMathOperator{\inttor}{\ensuremath{\underline{\stor}^{\cala}}}
\DeclareMathOperator{\id}{\mbox{id}\,}

\DeclareMathOperator{\im}{\mbox{im}\,}
\DeclareMathOperator{\spec}{\mbox{Spec}\,}
\DeclareMathOperator{\proj}{\mbox{Proj}\,}
\DeclareMathOperator{\supp}{\mbox{Supp}\,}
\DeclareMathOperator{\mo}{\mbox{mod}}
\DeclareMathOperator{\Mod}{\mbox{Mod}}
\DeclareMathOperator{\Bimod}{\mbox{Bimod}}

\DeclareMathOperator{\gr}{\mbox{gr}}
\DeclareMathOperator{\Gr}{\mbox{Gr}}

\DeclareMathOperator{\rk}{\mbox{rank}\,}
\DeclareMathOperator{\coker}{\mbox{coker}\,}

\DeclareMathOperator{\tr}{\mbox{tr}\,}
\DeclareMathOperator{\Hilb}{\mbox{Hilb}}

\DeclareMathOperator{\qcoh}{\mbox{QCoh}}
\DeclareMathOperator{\mmm}{\mathfrak{m}}
\DeclareMathOperator{\aaa}{\alpha}
\DeclareMathOperator{\s}{\sigma}
\DeclareMathOperator{\e}{\varepsilon}

\DeclareMathOperator{\w}{\omega}

\DeclareMathOperator{\PP}{\mathbb{P}}
\DeclareMathOperator{\Z}{\mathbb{Z}}
\DeclareMathOperator{\cale}{\mathcal{E}}
\DeclareMathOperator{\calf}{\mathcal{F}}
\DeclareMathOperator{\calm}{\mathcal{M}}
\DeclareMathOperator{\calo}{\mathcal{O}}
\DeclareMathOperator{\caln}{\mathcal{N}}
\DeclareMathOperator{\call}{\mathcal{L}}
\DeclareMathOperator{\calu}{\mathcal{U}}

\DeclareMathOperator{\cala}{\mathcal{A}}
\DeclareMathOperator{\calb}{\mathcal{B}}
\DeclareMathOperator{\ox}{\mathcal{O}_{X}}

\DeclareMathOperator{\oy}{\mathcal{O}_{Y}}

\DeclareMathOperator{\Gt}{\ensuremath{\tilde{G}}}
\DeclareMathOperator{\Ft}{\ensuremath{\tilde{F}}}

\begin{document}

\begin{center}
  \LARGE \textbf{Non-commutative Mori contractions and $\PP^1$-bundles}
\end{center}

\begin{center}
  DANIEL CHAN, ADAM NYMAN
\end{center}

\begin{center}
  {\em University of New South Wales, Western Washington University}
\end{center}

\begin{center}
e-mail address:{\em danielch@maths.unsw.edu.au, adam.nyman@wwu.edu}
\end{center}

\begin{abstract}
We give a method for constructing maps from a non-commutative scheme to a commutative projective curve. With the aid of Artin-Zhang's abstract Hilbert schemes, this is used to construct analogues of the extremal contraction of a $K$-negative curve with self-intersection zero on a smooth projective surface. This result will hopefully be useful in studying Artin's conjecture on the birational classification of non-commutative surfaces. As a non-trivial example of the theory developed, we look at non-commutative ruled surfaces and, more generally, at non-commutative $\PP^1$-bundles. We show in particular, that non-commutative $\PP^1$-bundles are smooth, have well-behaved Hilbert schemes and we compute its Serre functor. We then show that non-commutative ruled surfaces give examples of the aforementioned non-commutative Mori contractions. 
\end{abstract}

Throughout, all objects and maps are assumed to be defined over some algebraically closed base field $k$. The first author was supported by an ARC Discovery Project grant.

\vspace{2mm}
\section{Introduction}

In the last couple of decades, techniques from algebraic geometry have been succesfully applied to study non-commutative algebra giving birth to non-commutative algebraic geometry. A notable example includes Artin, Tate and Van den Bergh's study of Sklyanin algebras [ATV] using the Hilbert scheme of points for non-commutative graded algebras. A major research problem of non-commutative algebraic geometry is a ``birational classification'' of non-commutative projective surfaces [A].

The motivation for this project comes from a conjecture of Mike Artin's [A, conjecture~5.2] about this classification, which we paraphrase somewhat imprecisely as ``non-commutative surfaces are birationally ruled unless they are finite over their centre''\footnote{This lack of precision disappears if we consider non-commutative quadrics to be birationally ruled.}. Progress towards this conjecture seems to depend on i) a good understanding of birational equivalence classes, ii) a criterion for a non-commutative surface to be ruled and finally, iii) a criterion for being finite over the centre. This paper examines some ideas that may be useful in proving a criterion for ruledness for non-commutative surfaces. Again, Hilbert schemes play a central role, though in a rather different way to that in [ATV].

For the special case of non-commutative surfaces arising from orders over surfaces, the classification question has been settled using a non-commutative adaptation of Mori's minimal model program [CI]. The dichotomy in Artin's conjecture is also strongly reminiscent of the dichotomy in the Mori program.  Taking our cue from the minimal model program, our point of departure for a criterion for ruledness is the following commutative result.

\begin{thm}  \label{tcommmoricon}  
[KM, theorem~1.28(2)] 
Let $Y$ be a smooth commutative projective surface and $C \subset Y$ a curve which is extremal in the Kleiman-Mori cone and satisfies $K_Y.C < 0, C^2 = 0$. Then
\begin{enumerate}
\item there is a smooth curve $X$ and morphism $f:Y \lm X$ which contracts $C$.
\item the morphism $f$ is a $\PP^1$-fibration.
\end{enumerate}
\end{thm}
The commutative proof uses linear systems, which is the usual way one constructs maps in commutative algebraic geometry. However, this theory is not available in the non-commutative case, at least not yet. In fact, there are few methods for constructing morphisms in non-commutative algebraic geometry, and a major goal of this paper is to address this dearth.

Before explaining how we wish to generalise this result to the non-commutative setting, we need to recall some basic notions from non-commutative algebraic geometry. As is customary, we follow Grothendieck's philosophy that to study the geometry of a commutative scheme $Y$, we should study the category $\Mod Y$ of its quasi-coherent sheaves. For us, a non-commutative scheme or, to use Van den Bergh's terminology [VdB01] {\em quasi-scheme} $Y$, will be a $k$-linear Grothendieck category, that is, an abelian category with exact direct limits and a generator. We write $Y$ when we think of it geometrically and use geometric notation, and we write $\Mod Y$ when we want to think of it as a category.  A morphism $f:Y \lm X$ of quasi-schemes will then just be a pair of adjoint functors $f^*:\Mod X \lm \Mod Y, f_*:\Mod Y \lm \Mod X$, the motivating example being the usual pull-back and push-forward functors of quasi-coherent sheaves on commutative schemes $Y,X$. 

Naturally, the definition of quasi-schemes is too general to prove the types of geometric theorems we would like, so part~I is devoted to imposing various geometric conditions on the category $\Mod Y$. In particular, we define a notion of a non-commutative smooth proper $d$-fold (see section~\ref{sass}). For now, we merely note that the definition allows us to use the following geometric concepts:

\begin{itemize}
\item Intersection theory as developed by I. Mori and P. Smith [MS01].
\item Serre duality as studied by Bondal and Kapranov [BK].
\item Dimension theory.
\item Cohomology as developed by Artin and Zhang in [AZ94].
\item Hilbert schemes and base change as developed by Artin and Zhang in [AZ01].
\end{itemize}

In part~II, we introduce a method for constructing morphisms $f:Y \lm X$ from a quasi-scheme $Y$ to a commutative projective curve $X$. The basic idea can be seen as follows. First consider a morphism $f: Y \lm X$ of commutative schemes and let $\Gamma \subset Y_X := Y \times X$ be the graph of $f$. If $\pi:Y_X \lm Y$ denotes the projection, then $f^*$ can be factored as the Fourier-Mukai transform $f^* = \pi_*(\calo_{\Gamma} \otimes_X -)$. Now if $Y$ is more generally a quasi-scheme, $Y_X$ still makes sense (see section~\ref{sbbchange}) and given $\calm \in \mo Y_X$, so does $\calm \otimes_X -$. In section~\ref{sfibr} we give a definition of $\pi_*$ via relative Cech cohomology. The problem now is that $\pi_*$ is only left exact, so some conditions must be imposed to ensure the Fourier-Mukai transform $\pi_*(\calm \otimes_X -)$ is right exact, and so has a right adjoint. We identify such a sufficient condition in theorem~\ref{tpiexact}, namely, base point freedom, so named since it shares many of the features of base point freedom for linear systems. 

To see how to use this Fourier-Mukai transform, we return to our commutative guide. Given a commutative $\PP^1$-fibration $f:Y \lm X$ with fibre $C$, one can view $X$ as a component of the Hilbert scheme  corresponding to the subscheme $C\subset Y$ and $\calo_{\Gamma}$ is the corresponding universal quotient of $\oy$. Suppose now 
$Y$ is a non-commutative smooth proper surface. As in [AZ94], we will also specify a distinguished object $\oy \in \Mod Y$ which is used to define cohomology and is the analogue of the structure sheaf. We define (definition~\ref{dkneffint0}) what it means for a quotient $M \in \Mod Y$ of $\oy$ to be a $K$-non-effective rational curve with self-intersection zero. These are the analogues of the structure sheaves $\calo_C$, where $C$ is the curve in theorem~\ref{tcommmoricon}. The key point is that the definition of non-commutative smooth proper surfaces furnishes us with enough geometric concepts to make the definitions of rational, self-intersection zero etc. We can now consider the Hilbert scheme of quotients of $\oy$ corresponding to $M$. After proving a smoothness criterion for Hilbert schemes, we show (in corollary~\ref{cfamily}) that this scheme is a generically smooth projective curve $X$ and so can consider the Fourier-Mukai transform. However, the base point free condition is complicated by S. P.  Smith's discovery of ``strange points'' (see [SV]), which are in particular, zero-dimensional modules with non-zero self-intersection. The main result in part~II is the following analogue of theorem~\ref{tcommmoricon}i).

\begin{thm}  \label{tmain}  
Let $Y$ be a non-commutative smooth proper surface that has a 1-critical $K$-non-effective rational curve $M$ with self-intersection zero. Suppose furthermore that for every zero-dimensional simple quotient $P \in \Mod Y$ of $M$, we have the intersection product $M.P\geq 0$. Then there is a a morphism $f:Y \lm X$ where $X$ is a generically smooth projective curve. 
\end{thm}

The morphism $f:Y \lm X$ of non-commutative schemes that results is called a {\em non-commutative Mori contraction}. Naturally, one would like an analogue of theorem~\ref{tcommmoricon}ii) too, and this forms part of ongoing research. 

Given the desired application to Artin's conjecture, one would hope that the non-commutative ruled surfaces of Patrick and Van den Bergh (see [Pat],[VdB01p]) give non-trivial examples of the theory in part~II. Part~III is devoted to showing this. Our main result in this part is
\begin{thm}  \label{tqrsOK}  
A non-commutative ruled surface is a non-commutative smooth proper surface.
\end{thm}
In fact we show more generally that non-commutative $\PP^1$-bundles are smooth and explicitly give the Serre functor. We also show that their Hilbert schemes are well-behaved. Further, given the natural fibration of a non-commutative ruled surface $f:Y \lm X$, we show that the fibres of $f$ are $K$-non-effective rational curves with self-intersection zero, and $f$ is the associated non-commutative Mori contraction. 

A morphism $f:Y \lm X$ of quasi-schemes has very little structure in comparison to a morphism of commutative schemes. For example, even if $Y,X$ are equipped with structure sheaves $\oy,\ox$, there is no reason to suppose there is any relationship between $\oy$ and $f^*\ox$. This makes it hard to prove results like a Leray spectral sequence linking the cohomology on $Y$ with that on $X$. It seems that to be able to extract more information from a morphism of quasi-schemes, we need to restrict the morphisms under consideration. In part~IV, we study properties of non-commutative Mori contractions $f: Y \lm X$. In this case, there is a natural map $\nu:\oy \lm f^* \ox$ and the driving question here is to find criteria to guarantee $\nu$ is an isomorphism. In the process, we will also study the higher direct images $R^if_*$. We hope that the material in this part will be useful in proving an analogue of theorem~\ref{tcommmoricon}ii). 

\textbf{Notation:} Throughout this paper, $Y$ will always denote some quasi-scheme and, by default, all unadorned Ext and $\otimes$ symbols will be taken over $Y$. At the beginning of each section, we will re-state any additional hypotheses on $Y$. 

\part{Non-commutative smooth proper $d$-folds}

As is common in the non-commutative community, we shall follow Grothendieck's philosophy and do geometry via the category of quasi-coherent sheaves. To this end, we consider a {\em quasi-scheme} $Y$ which is the data of a Grothendieck category $\Mod Y$ (over $k$), that is $\Mod Y$ is a $k$-linear category with exact direct limits and a generator. Objects in $\Mod Y$ will usually be called {\em $Y$-modules}. We also let $\mo Y$ denote the full subcategory of noetherian objects.

The example to keep in mind comes from a connected graded, locally finite $k$-algebra $A = k \oplus A_1 \oplus A_2 \oplus \dots$ where {\em locally finite} means $\dim_k A_i < \infty$ for all $i$. If $\Gr A$ denotes the category of graded $A$-modules and tors, the full Serre subcategory consisting of $A_{>0}$-torsion modules, then the quotient category $Y = \proj A:=\Gr A/\mbox{tors}$ is an example of a quasi-scheme. The motivation is Serre's theorem which states that if $A$ is the (commutative) homogeneous coordinate ring of a projective scheme $Y$, then $\proj A$ is naturally equivalent to $\qcoh Y$. The quasi-schemes we will be interested in will be {\em noetherian} in the sense that $\Mod Y$ is locally noetherian, that is, they have a set of noetherian generators. 

This part is primarily concerned with addressing the question, ``Which quasi-schemes are the non-commutative analogues of smooth proper varieties?''. We start in section~\ref{sbbchange} with recounting the notion of base change and Hilbert schemes for categories developed in [AZ01]. Hilbert schemes will be a fundamental tool for us. In section~\ref{sass}, we propose a definition of a non-commutative smooth proper $d$-fold based on various geometric conditions which we impose on a quasi-scheme $Y$. The list of hypotheses is rather long so in section~\ref{sproj}, we show that in the projective case of $Y = \proj A$, these conditions do follow from hypotheses on $A$ that other authors have studied in the past. 

\section{Background on base change and Hilbert schemes}  \label{sbbchange}  

Let $R$ be a commutative $k$-algebra. Section~B of [AZ01] is devoted to the notion of base change from $k \lm R$ for arbitrary categories. In this section, we review their results and generalise their notion of $R$-objects to $X$-objects where $X$ is a separated scheme. We will always work over a field $k$ as opposed to more general commutative rings as in [AZ01]. This simplifies the treatment somewhat, since then any Grothendieck category has $k$-flat generators and the subtleties involved with defining tensor products disappear. 

Let $Y$ be a quasi-scheme defined over the ground field $k$ as always. An {\em $R$-object of $\Mod Y$} is a $Y$-module $\calm \in \Mod Y$ equipped with an $R$-action, that is, a ring morphism $R \lm \End \calm$. Morphisms of such objects are $Y$-module morphisms compatible with the $R$-action. These objects form an $R$-linear Grothendieck category [AZ01, proposition~B2.2] denoted $\Mod Y_R$. 

Given $\calm \in \Mod Y_R$ one has a tensor functor $\calm \otimes_R-: \Mod R \lm \Mod Y_R$ defined by declaring $\calm \otimes_R R = \calm$ and insisting $\calm \otimes_R-$ is right exact and commutes with direct sums [AZ01, proposition~B3.1]. We may thus define $\calm$ to be {\em $R$-flat} (or just flat) if this functor is also left exact [AZ01, section~C.1]. Also, given a morphism of commutative $k$-algebras $R \lm S$, there is a base change functor $\Mod Y_R \lm \Mod Y_S$. With this notion of base change,  Grothendieck's theory of flat descent holds [AZ01, theorem~C8.6]. 

The categories $\Mod Y_R$ thus form a stack and this allows us to do base change via arbitrary separated schemes. We copy the definition from that of quasi-coherent modules on a stack (see for example [Vistoli, Appendix]). Let $X$ be a separated scheme on which we put the small Zariski site. An $X$-object $\calm$ in $\Mod Y$ is the data of an $R$-object $\calm_R \in \Mod Y_R$ for each affine open $\spec R \subseteq X$ which is compatible with base change. If $U = \spec R$, then we will usually write $\calm(U)$ for $\calm_R$. The collection of $X$-objects in $\Mod Y$ naturally forms an abelian category with exact direct limits denoted $\Mod Y_X$. These $X$-objects may also be defined via descent data. We see immediately that tensoring with $R$-modules extends to give a bifunctor 
$$\Mod Y_X \times \Mod X \lm \Mod Y_X .$$ 
Similarly, given an affine morphism of separated schemes $g:X' \lm X$, there is a base change functor $\Mod Y_X \lm \Mod Y_{X'}$. 

We may now define a Hilbert functor [AZ01, section~E2]. Fix $F \in \mo Y$ and let $\mathcal{R}$ denote the category of all commutative $k$-algebras $R$ such that $Y_R$ is noetherian. If $Y$ is noetherian, then Hilbert's basis theorem [AZ01, theorem~B5.2] ensures this includes all algebras of finite type. The Hilbert functor $\Hilb(F):\mathcal{R} \lm \mbox{Sets}$ sends $R \in \mathcal{R}$ to the set isomorphism classes of $R$-flat quotients of $F\otimes_k R$ in $\mo Y_R$. If this functor is representable by a scheme, we shall call it the {\em Hilbert scheme of quotients of $F$}. This gives natural examples of $X$-objects, for suppose $Y$ is noetherian and $X$ is some subscheme of the Hilbert scheme. Then there exists a universal object $\calm \in \Mod Y_X$ defined as follows. Given any morphism $i:\spec R \lm X$, there exists a (tautological) $R$-flat object $i^*\calm \in \Mod Y$ corresponding to $i$, and these objects are compatible with base change. We thus obtain an $X$-object. 

Clearly it would be desirable to have $\mathcal{R}$ to be the set of all noetherian $k$-algebras so, 
following Artin-Zhang [AZ01] we make the 

\begin{defn}  \label{dsnoeth}  
If $Y_R$ is noetherian for every noetherian $R$, then we will say that $Y$ is {\em strongly noetherian}. 
\end{defn}

We will also consider the following conditions on a category.

\begin{defn}   \label{dextfinite}   
A quasi-scheme is {\em Ext-finite (respectively Hom-finite)} if for every commutative noetherian $k$-algebra $R$ and noetherian modules $M \in \mo Y, N \in \mo Y_R$ we have that $\Ext^i(M,N)$ is a finite $R$-module for all $i$ (respectively for $i=0$). 
\end{defn}

One important geometric result we have is Artin-Zhang's version of generic flatness. Let $Y$ be a strongly noetherian, Hom-finite quasi-scheme and $R$ a commutative reduced noetherian algebra. Given any $\calm \in \mo Y_R$, there is a non-zero-divisor $s \in R$ such that $\calm \otimes_R R[s^{-1}]$ is flat over $R[s^{-1}]$  (see [AZ01, theorem~C5.1, corollary~C7.4]).

We now introduce a local sections functor. Let $\iota:Z \hookrightarrow X$ be a closed embedding and $\calm \in \Mod Y_X$. We seek to define $\iota^{!}\calm \in \Mod Y_Z$. Let $\spec R \subseteq X$ be an affine open set and $I\triangleleft R$ be the ideal corresponding to $Z$. One can define the {\em annihilator of $I$} in the usual way 
$$ \caln_R := \bigcap_{t \in I} \ker(\calm_R \xrightarrow{t} \calm_R) .$$
This is compatible with flat base change so we obtain an $X$-object $\caln$ which is supported on $Z$ in the sense that $\caln_R$ is an $R/I$-object, so in particular, is zero if $\spec R \cap Z = \varnothing$. The sheaf property now ensures that this corresponds to a $Z$-object $\iota^{!}\calm \in \Mod Y_Z$. Abusing notation, we also let $\iota^{!} \calm$ denote the corresponding submodule of $\calm$. 

\begin{prop}  \label{phplanecuts}  
Let $Y$ be a strongly noetherian, Hom-finite quasi-scheme. 
Suppose $X$ is a quasi-projective scheme and $\calm \in \mo Y_X$. If $\iota: H \hookrightarrow X$ is the inclusion of a generic hyperplane, then $\iota^! \calm = 0$.
\end{prop}
\textbf{Proof.} It suffices to prove the proposition in the case where $X = \spec R$ is affine and $H \subset X$ is defined by $t=0$ for some $t\in R$. We argue by d\'{e}vissage. Let $I$ be the radical of $R$. Since $\calm$ has a finite filtration whose factors are $R/I$-objects, we may assume that $I=0$ so $R$ is reduced. By generic flatness, there is a non-zero-divisor $r \in R$ such that $\calm\otimes_R R[r^{-1}]$ is flat over $R[r^{-1}]$ so multiplication by $t$ is injective on $\calm\otimes_R R[r^{-1}]$. We need to show that multiplication by $t$ is injective on $\calm$. This follows since, by [AZ01, proposition~B6.2], the kernel of $\calm \lm \calm \otimes_R R[r^{-1}]$ is the $r$-torsion submodule of $\calm$ so we may appeal to the inductive hypothesis on $\spec R/(r^n)$ for $n \gg 0$. 

We have the following generalisation of [AZ01, corollary~B3.17, proposition~B5.1] which gives a useful sufficient criterion for $\Mod Y_X$ to be a Grothendieck category.

\begin{prop}  \label{pgenYX}  
Let $Y$ be a Hom-finite strongly noetherian quasi-scheme. 
Suppose $X$ is a quasi-projective scheme and $\ox(1)$ a very ample line bundle. 
If $\{\call_{\aaa}\}_{\aaa}$ forms a set of noetherian generators for $\Mod Y$ then $\{\call_{\aaa} \otimes_k \ox(j)\}_{j \in \Z, \aaa}$ forms a set of noetherian generators for $\Mod Y_X$. In particular, $Y_X$ is a noetherian quasi-scheme.
\end{prop}
\textbf{Proof.} First note that $\call_{\aaa} \otimes_k \ox(j)$ are noetherian since their restriction to any affine open is noetherian by [AZ01, proposition~B5.1]. Given $\calm,\caln \in \Mod Y_X$ and distinct morphisms $\phi,\phi': \calm \lm \caln$, it suffices to find a non-zero morphism of the form $\psi:\call_{\aaa} \otimes_k \ox(j) \lm \calm$ for appropriate $\aaa,j$ with the property that $\phi\psi \neq \phi'\psi$. 

As usual, we may find an affine open cover $\{\spec R_0,\ldots ,\spec R_n\}$ of $X$ where each open is the complement of a generic hyperplane. We may re-index to suppose that $\phi,\phi'$ differ on $\spec R_0$. Suppose the hyperplane $H = X - \spec R_0$ is defined in $\spec R_i$ by the zeros of $t_i \in R_i$. We let $\calm_i,\calm_{ij}$ denote the restriction of $\calm$ to $\spec R_i, \spec R_i \cap \spec R_j$ respectively.

By proposition~\ref{phplanecuts}, multiplication by $t_i$ is injective on $\calm_i$ so $\calm_i$ embeds in $\calm_{0i}$ and in fact $\calm_{0i} = \cup_j t^{-j} \calm_i$. Now the $\call_{\aaa}$ generate $\Mod Y$, so  there exists a non-zero morphism $\psi:\call_{\aaa} \lm \calm_0$ such that at least on $\spec R_0$ we have $\phi\psi \neq \phi'\psi$. Now $\call_{\alpha}$ is noetherian, so the image of $\psi$ in $\calm_{0i}$ lands in $t^{-j}\calm_i$ for $j\gg 0$.  Hence we may extend $\phi$ to a morphism $\call_{\aaa}\otimes_k \ox(-j) \lm \calm$ for $j \gg 0$ and we are done. 

\vspace{2mm}

\vspace{2mm}
\section{Some geometric hypotheses}  \label{sass}  

In this section, we motivate and propose a definition of a non-commutative smooth proper $d$-fold. The idea is to impose conditions which allow us to adapt various geometric tools to the non-commutative setting. Despite appearances, we did strive for {\em lex parsimoniae}. In the next section, we will discuss ways of simplifying the hypotheses in the case where $Y=\proj A$. Our main interest will be in surfaces, and perhaps the axioms we give are most appropriate in that context only. 

\begin{defn}  \label{dncdfold}  
For $d \geq 2$, a {\em non-commutative smooth proper $d$-fold} is a strongly noetherian quasi-scheme $Y$ which satisfies the 6 sets of hypotheses 1)-6) below. 
\end{defn}

\vspace{2mm}\noindent
{\em 1. Smooth, proper of dimension $d$}

Smooth of dimension $d$, means that the global dimension of $\Mod Y$ is $d$ so $\Ext^{d+1}_Y(-,-) = 0$. Proper means that $Y$ is Ext-finite (see definition~\ref{dextfinite}).

Given these hypotheses we can define for $N,N' \in \mo Y$ the pairing,  
\[ \xi(N,N') : = \sum_{i} (-1)^i\dim \Ext^i(N,N')  .  \]
Thus if $N,N'$ represent ``curve'' modules on a surface (so $d=2$), then the intersection product \`a la Mori-Smith [MS01] $N.N' = -\xi(N,N')$ is always finite.

\vspace{2mm}\noindent
{\em 2. Gorenstein}

We assume there exists an auto-equivalence $- \otimes \omega_Y:\mo Y \lm \mo Y$ such that $-\otimes \omega_Y[d]: D^b(Y) \lm D^b(Y)$ (where $[d]$ denotes shift in the derived category) gives Bondal-Kapranov-Serre duality [BK], that is, functorial isomorphisms in $N,N' \in \mo Y$

\[ \Ext^i(N,N') \simeq \Ext^{d-i}(N',N\otimes \omega_Y)^*   .\]

Recall that for commutative Gorenstein schemes, the canonical sheaf $\w_Y$ is invertible so $-\otimes \omega_Y$ induces an auto-equivelance. 

\vspace{2mm}\noindent
{\em 3. Compatible dimension function}

i) Firstly, we want an exact dimension function $\dim: \mo Y \lm \{-\infty,0,1, \ldots ,d\}$. Recall that exactness means,  given an exact sequence in $\mo Y$ 
\[  0 \lm M' \lm M \lm M'' \lm 0 ,\]
we have $\dim M = \max \{\dim M' , \dim M'' \}$. Traditionally, there is another condition on $\dim$ as well, but we do not seem to need it and, in fact, is not well-defined in the context of arbitrary noetherian categories.  

ii) We also want the dimension function to be compatible with the Serre functor in the sense that for every $M \in \mo Y$ we have $\dim M = \dim M\otimes \omega_Y$. 

\vspace{2mm}\noindent
{\em 4. Classical cohomology}

As part of the data, we consider a $d$-critical object $\oy \in \mo Y$ with which we compute global cohomology $H^i := \Ext^i(\oy,-)$. We call $\oy$ the {\em structure sheaf} of $Y$. By $d$-critical we mean that $\dim \oy = d$ and for any non-zero submodule $I \leq \oy$ we have $\dim \oy/I < d$ so in particular by exactness of dimension, $\dim I = d$. 

We assume further that for any zero dimensional module $P$, 
\[ i) \Ext^1(P,\oy) = 0 ,\ \ \ \mbox{and} \ \ \ ii) h^0(P) := \dim_k H^0(P) \neq 0   .\]
Note that part i) should fail for a curve $Y$ which is the reason why we have restricted our definition to $d \geq 2$. 

\vspace{2mm}\noindent
{\em 5. Halal Hilbert schemes}

We assume that for any $F \in \mo Y$, the Hilbert functor $\Hilb(F)$ is representable by a separated scheme, locally of finite type which is furthermore, a countable union of projective schemes.

\vspace{2mm}\noindent
{\em 6. No shrunken flat deformations}

If $M,M'$ are two members in a flat family parametrised by a connected scheme of finite type, then any injective map $M \hookrightarrow M'$ or surjective map $M \lm M'$ is an isomorphism. That is, we assume that a module cannot be flatly deformed into a proper submodule or quotient module. It is possible that this hypothesis can be reduced to a simpler one, but we have not studied this possibility. 

\begin{defn}  \label{dncsurface} 
A {\em non-commutative smooth proper surface} $Y$ is a non-commutative smooth proper 2-fold. 
\end{defn}

For some of the proofs, we will unfortunately need to strengthen some of the hypotheses above. 
\begin{itemize}
\item[{3iii)}] We say an exact dimension function is {\em continuous} if the following condition holds. Given any flat family $\calm\in \mo Y_X$ of $Y$-modules parametrised by an integral scheme $X$ of finite type,  we have that the function $p \mapsto \dim \calm \otimes_X k(p)$ is constant as $p$ varies over the closed points of $X$. 
\item [{3iv)}] We say an exact dimension function $\dim$ is {\em finitely partitive}, if for any $M \in \mo Y$ there exists a bound $l$  on the length of strictly decreasing chains $M=M_0 \supset M_1 \supset \ldots \supset M_j$ with factors all of dimension $\dim M$. 
\end{itemize}

We prove the following analogue of Grothendieck's vanishing theorem in the surface case. 

\begin{cor} \label{cvanish}
Let $Y$ be a non-commutative smooth proper surface and $N \in \mo Y$. 
\begin{enumerate}
\item If $\dim N < 2$ then $H^2(N) = 0$.
\item If $N$ is zero dimensional then $H^1(N) = 0$.
\end{enumerate} 
\end{cor}
\textbf{Proof.} If $\dim N <2$ then any quotient of $N$ has dimension less than 2. The assumption of compatibility of the dimension function with the Serre functor shows $\oy\otimes \w_Y$ is 2-critical too so $\Hom_Y(N,\oy \otimes \w_Y) = 0$. BK-Serre duality gives $H^2(N) = \Ext^2(\oy,N) = \Hom(N,\oy \otimes \w_Y)^* = 0$. Suppose now that $\dim N = 0$. Then $N \otimes \w_Y^{-1}$ is zero dimensional too so the assumptions on classical cohomology give $H^1(N) = \Ext^1(N,\oy \otimes \w_Y)^* = \Ext^1(N \otimes \w_Y^{-1},\oy)^* = 0$. 
\vspace{2mm} 

We will give examples of non-commutative smooth proper surfaces in the next section and in part~III.

\vspace{2mm}
\section{Projective case}  \label{sproj} 

In this section, we discuss hypotheses 1)-6) of section~\ref{sass} in the special case of $Y= \proj A$ where $A$ is a strongly noetherian, connected graded $k$-domain, finitely generated in degree one. In particular, we will define a non-commutative smooth projective $d$-fold which will hopefully be a less distasteful object to the reader than that of a non-commutative smooth proper $d$-fold! Throughout this section, we let $\mmm$ be the augmentation ideal $A_{>0}$ and $\Gamma_{\mmm}$ be the $\mmm$-torsion functor. All the $A$-modules we consider will be graded so we shall omit the adjective graded. However, $\Ext^i_A$ will denote ext groups in the ungraded category. We let $\gr A$ denote the full subcategory of $\Gr A$ whose objects are the noetherian ones. We will usually use a letter like $P$ to denote both $A$-modules and $Y$-modules. When such abuse can cause confusion, we will often let $P_{\bullet}$ denote some $A$-module which represents the $Y$-module $P$. 

We first recall Artin and Zhang's $\chi$ condition [AZ01, (C6.8)]. Let $R$ be an algebra. We say that $A \otimes_k R$ satisfies $\chi$ if for every noetherian $A \otimes_k R$-module $M$ we have that the $A \otimes_k R$-module $\Ext_A^i(A/A_{\geq n} \otimes_k R, M)_{\geq d}$ is noetherian for all $i,n,d$. We say that $A$ satisfies {\em strong $\chi$} if $A \otimes_k R$ satisfies $\chi$ for every commutative noetherian algebra $R$. This strong $\chi$ condition guarantees Ext-finiteness by [AZ01, proposition~C6.9]. 

Recall that $A$ has a balanced dualising complex if $A$ and its opposite algebra $A^o$ satisfy the $\chi$  condition and their torsion functors $\Gamma_{\mmm}, \Gamma_{\mmm^o}$ have finite cohomological dimension [VdB97, theorem~6.3]. If $A$ has a balanced dualising complex then, since we are assuming strongly noetherian, we know from [AZ01, proposition~C6.10] that $A$ satisfies strong $\chi$. If we further assume that $Y$ is smooth, then the balanced dualising complex induces the BK-Serre functor (see [NV04, theorem~A.4, corollary~A.5]).

We now consider the hypothesis of compatible dimension function in the projective case. We start with two important dimension functions in $\gr A$, namely, the Gelfand-Kirillov dimension and the canonical dimension. Let $M_{\bullet} \in \gr A$. Then the Gelfand-Kirillov dimension of $M_{\bullet}$ is the growth rate of the function $f(n) := \dim_k (\oplus_{i\leq n} M_n)$, that is
$$gk(M_{\bullet}) = \inf\{d| f(n) \leq n^d \ \mbox{for}\ n \gg 0\} .$$
More details on the Gelfand-Kirillov dimension can be found in [McR, section~8.1] or [KL]. Now under our hypotheses, the Gelfand-Kirillov dimension is exact [McR, proposition~8.3.11] though not necessarily integer valued. Thus if we assume it is integer valued, then letting $M$ denote the $Y$-module corresponding to $M_{\bullet}$, we obtain a well-defined exact dimension function on $\mo Y$ by setting $\dim M$ to be $gk(M_{\bullet}) - 1$ whenever $M_{\bullet}$ is not non-zero $\mmm$-torsion. Suppose now that the Serre functor is given by an actual noetherian $A$-bimodule $\omega_Y$ and that $-\otimes \w_Y^{-1}$ is also given by tensoring with a noetherian bimodule. Ideal invariance (see [KL, proposition~5.6]) ensures then that $\dim$ is compatible with the Serre functor. Lastly, $\dim$ is continuous since (tails of) Hilbert functions are preserved in flat families by [AZ01, lemma~E5.3]. This observation also shows there are no shrunken flat deformations.

The canonical dimension (see [YZ] for more information) can be defined for $M_{\bullet} \in \gr A$ by
$$  c.\dim(M_{\bullet}) := \max\{i| R^i\Gamma_{\mmm} M_{\bullet} \neq 0 \}  .$$
If $A$ has an Auslander dualising complex (see [YZ, definition~2.1]), then the canonical dimension is a finitely partitive exact dimension function [YZ, theorem~2.10] so we can similarly define the dimension function on $Y$ by $c.\dim -1$. We do not know if it is compatible with the Serre functor or is continuous. 

For $Y = \proj A$, the natural choice for $\oy$ is the image of $A$ in $\proj A$. It is critical for the Gelfand-Kirillov dimension since we are assuming $A$ is a domain [McR, proposition~8.3.5]. Finally, Halal Hilbert schemes holds if $A$ satisfies strong $\chi$ by [AZ01, theorem~E5.1]. 

\begin{defn}  \label{dproj}   
Let $A$ be a strongly noetherian, connected graded $k$-domain which is finitely generated in degree one. Then $Y= \proj A$ with structure sheaf $\oy = A$ is a {\em non-commutative smooth projective $d$-fold} ($d\geq 1$) if $A$ also satisfies the three hypotheses below.
\end{defn}

\begin{enumerate}
\item $Y$ is smooth of dimension $d$. 
\item There are noetherian $A$-bimodules $\omega_Y,\w_Y^{-1}$ such that $\omega_Y[d+1]$ is an Auslander balanced dualising complex of injective dimension $d+1$ and $\omega^{-1}_Y$ is inverse to $\omega_Y$ in $\proj A$, in the sense that  $- \otimes \w_Y^{-1},\ - \otimes \w_Y$ induce inverse auto-equivalences on  $\proj A$. 
\item $A$ is $gk$-Macaulay in the sense that the Gelfand-Kirillov dimension and canonical dimension coincide. We then define the dimension function on $Y$ by $\dim = gk - 1 = c.\dim -1$. 
\end{enumerate}

The preceding discussion shows that non-commutative smooth projective $d$-folds satisfy hypotheses 1)-3),5) and 6) of a non-commutative smooth proper $d$-fold. We look at the remaining hypothesis of classical cohomology. We know $\oy$ is $d$-critical. Below, we need to use local duality [VdB97, theorem~5.1] which states that, if $(-)^{\vee}$ denotes the graded $k$-vector space dual then $R\Gamma_{\mmm}(-)^{\vee} \simeq \Rhom_A(-,\w)$ where $\w$ is the dualising complex.  

\begin{prop} \label{pcandim0} 
Let $Y = \proj A$ be a non-commutative smooth projective $d$-fold and $M \in \gr A$ a non-zero module. Then $c.\dim M = 0$ if and only if $M$ is $\mmm$-torsion. Also, $c.\dim M = d+1$ if and only if $A$ is not torsion in the sense that $M \otimes_A Q(A) \neq 0$.  
\end{prop}
\textbf{Proof.} If $M$ is $\mmm$-torsion then $R\Gamma_{\mmm} M = M$ (see for example [AZ94, proposition~7.1]) so it has canonical dimension zero. Conversely suppose $c.\dim M = 0$. Let $T$ be its $\mmm$-torsion submodule and assume by way of contradiction that $T \neq M$. Now $c.\dim T \leq 0$ so in this case we can use exactness of $c.\dim$ to see that $c.\dim M/T = 0$. We may assume thus that $M$ is $\mmm$-torsion-free. Then $R\Gamma_{\mmm} M = 0 = \Rhom_A(M,\w_Y[d+1])^{\vee}$ which by duality means $M=0$, a contradiction.

We know $-\otimes \w_Y$ preserves canonical dimension so for the other result, it suffices to show that $M \otimes_A Q(A) \neq 0$ iff $\Hom_A(M,A) \neq 0$. If $M \otimes_A Q(A) \neq 0$ then it is isomorphic to $Q(A)^n$ for some $n >0$. We can thus choose a non-zero map $M \lm Q(A)$ and, by clearing denominators, we can construct a non-zero  map $M \lm A$. Conversely, if $M\lm A$ has non-zero image $I$, then $M \otimes_A Q(A)$ surjects onto $I \otimes_A Q(A) \simeq Q(A)$ so must be non-zero too. This completes the proof of the proposition. 

\vspace{2mm}

\begin{lemma}  \label{lproj} 
Suppose that $A$ satisfies $\chi$. Let $M \in \Gr A$ be an $\mmm$-torsion-free module representing a noetherian object in $\proj A$. Then $M_{\geq i} \in \gr A$ for any $i$.
\end{lemma}
\textbf{Proof.} This is well-known. We sketch the proof here. Pick $N \leq M_{\geq i}$ a noetherian $A$-module with $M_{\geq i}/N$ $\mmm$-torsion. The $\chi$ condition ensures that $\Ext^1_A(k,N)_{\geq i}$ is finite dimensional so, as a graded module is 0 in degree $\geq r$ for some $r>>0$. This shows $N_{\geq r} = M_{\geq r}$. The universal extension of $N_{\geq r}$ by $\Ext^1_A(k,N_{\geq r})_{r-1}$ is still noetherian, so we are done by downward induction on $r$. 

\vspace{2mm}

\begin{prop}  \label{pcancaseass}  
Let $Y = \proj A$ be a non-commutative smooth projective $d$-fold where $d\geq 2$. Then for any non-zero $P\in \mo Y$ with $\dim P = 0$ we have 
$$ i) \ \Ext^1_Y(P,\oy) =0 \ , \quad ii) \ h^0(P) \neq 0$$
\end{prop}
\textbf{Proof.} We let $P_{\bullet}$ be an appropriate $A$-module representing $P$. 

To prove i), consider an exact sequence of finite $A$-modules
\[  (*) \ \ \ 0 \lm A \xrightarrow{\phi} N_{\bullet} \xrightarrow{\phi'} P_{\bullet} \lm 0.\]
We need to show it splits in $\proj A$, for which it suffices to assume that $P_{\bullet}$ is  $\mmm$-torsion-free so by local duality $\Ext^{d+1}_A(P_{\bullet},\w_Y)= \Gamma_{\mmm}(P_{\bullet})^{\vee} = 0$. Now $P$ is zero dimensional so $\Ext^i_A(P_{\bullet},\w_Y)=0$ for $i<d$. Also, $A$ is maximal Cohen-Macaulay by assumption so the long exact sequence in cohomology shows that for $i \neq 0$ or $d$ we have 
\[ \Hom_A(N_{\bullet},\w_Y) = \Hom_A(A,\w_Y) = \w_Y, \ \Ext^i_A(N_{\bullet},\w_Y) = 0, \ \Ext^d_A(N_{\bullet},\w_Y) = \Ext^d_A(P_{\bullet},\w_Y)  .\]
Hence the double-Ext spectral sequence [YZ, proposition~1.7] related to $\Rhom_A(\Rhom_A(N_{\bullet},\w_Y),\w_Y)=N_{\bullet}$ has only two rows containing terms 
\[E^{p,0}_2 = \Ext^p_A(\Ext^0_A(N_{\bullet},\w_Y),\w_Y) = Ext^p_A(\w_Y,\w_Y) \]
which is $A$ when $p=0$ and zero otherwise, and also terms
\[E^{p,-d}_2 = \Ext^p_A(\Ext^d_A(N_{\bullet},\w_Y),\w_Y) = \Ext^p_A(\Ext^d_A(P_{\bullet},\w_Y),\w_Y) \]
which, by duality on $P_{\bullet}$ is $P_{\bullet}$ when $p=d$ and zero otherwise. There is thus an exact sequence of the form
\[  0 \lm P_{\bullet} \xrightarrow{\psi} N_{\bullet} \xrightarrow{\psi'} A \lm 0  .\]
Note that $\phi'\psi$ is injective otherwise $\ker \phi' \cap \im \psi$ is a non-zero submodule of $\ker \phi' = \im \phi$ which contradicts the fact that $\im \phi \simeq A$ is torsion-free. Now $P_{\bullet}$ is also locally finite so $\phi'\psi$ is an isomorphism which shows (*) splits in $\gr A$ and hence in $Y$. 

Finally we prove part ii). We may assume that $P_{\bullet}$ is $\mmm$-torsion-free and saturated in the sense that $\Ext^1_A(k,P_{\bullet}) = 0$. Suppose that $0=H^0(P) = P_0$. Since $A$ is generated in degree one and $P_{\bullet}$ is $\mmm$-torsion-free, we must have $P_{\leq 0} = 0$. By lemma~\ref{lproj}, $P_{\bullet}$ is noetherian as an $A$-module. The saturation and $\mmm$-torsion-free condition imply $\Gamma_{\mmm} P_{\bullet} = R^1\Gamma_{\mmm} P_{\bullet} = 0$. Also, since $c.\dim P_{\bullet} = 1$ we have $R^i\Gamma_{\mmm} P_{\bullet} = 0$ for $i >1$. Thus $R\Gamma_{\mmm}P_{\bullet} = 0$, a contradiction. 
\vspace{2mm}

We have thus proved
\begin{prop}   \label{pprojOK}  
A non-commutative smooth projective $d$-fold where $d\geq 2$ is a non-commutative smooth proper $d$-fold with a continuous finitely partitive dimension function. 
\end{prop}

Examples of non-commutative smooth projective surfaces are quantum planes $\proj A$ where $A$ is a 3-dimensional Sklyanin algebra and the smooth non-commutative quadrics of Smith-Van den Bergh [SV].

\vspace{2mm}

\part{Constructing Mori contractions}

In this part, we will show how a type of Fourier-Mukai transform can be used to construct morphisms in non-commutative algebraic geometry (section~\ref{sfibr}). We consider non-commutative analogues of $K$-negative curves with self-intersection zero in section~\ref{smoduli} as well as their flat deformations. The aim is to use the universal deformation as the kernel of the Fourier-Mukai transform to obtain a type of Mori contraction. Such a contraction theorem is given in section~\ref{smoricon}. We will need several results regarding generic behaviour of $Y$-modules in a flat family. These are given in section~\ref{sscont} and \ref{sbertini}. Section~\ref{sconverse} shows that most flat morphisms from a noetherian quasi-scheme to a commutative projective variety come from this Fourier-Mukai construction.

\vspace{2mm}
\section{Semicontinuity}  \label{sscont}  

In this section, we let $Y$ be any strongly noetherian, Ext-finite quasi-scheme. We prove some semicontinuity results. In the projective case, more general results can be found in [NV05, Appendix~A]. 
Given a flat family $\calm \in \mo Y_X$, we write loosely $M \in \calm$ to mean $M \in \mo Y$ is a closed fibre of $\calm$.

\begin{prop}  \label{pshiftflat}  
Let $X$ be a noetherian scheme and $\calm\in \mo Y_X$ a flat family of $Y$-modules. Suppose given an  exact functor $-\otimes I: \Mod Y \lm \Mod Y$ which commutes with arbitrary direct sums. Then $\calm \otimes_Y I$ is a flat family of $Y$-modules too.
\end{prop}
\textbf{Proof.} We may assume the flat family is defined over $\spec R$ where $R$ is a commutative noetherian ring. Now $\calm$ is an object of $\Mod Y$ with an action of $R$ so since $-\otimes I$ is a functor, $\calm \otimes I$ is also an object of $\Mod Y$ with an action of $R$. 

Recall that for an $R$-module $L$, the tensor product $\calm \otimes_R L$ is defined by considering a presentation
\[  \bigoplus_K R \lm \bigoplus_J R \lm L \lm 0 \]
and asserting exactness of
\[  \bigoplus_K \calm \lm \bigoplus_J \calm \lm \calm\otimes_R L \lm 0 .\]
Tensoring this sequence with $I$ shows that $\otimes I$ commutes with $\otimes_R L$. Hence $\calm \otimes I$ is also flat over $R$. 

\vspace{2mm}

We prove now a mild (semi-)continuity result. Recall for $N,N' \in \mo Y$ we defined  
\[ \xi(N,N') : = \sum_{i} (-1)^i\dim \Ext^i(N,N')    .\]

\begin{prop} \label{pcont}  
Let $R$ be a Dedekind domain of finite type over $k$ and $\caln \in \mo Y_R$ be a flat family of noetherian $Y$-modules over $R$. We fix $N' \in \mo Y$,  and let $N\in \caln$ vary over the closed fibres of $\caln$. 
\begin{enumerate}
\item The dimension of $\Ext^i(N',N)$ is upper semi-continuous as a function of $N$.
\item If $Y$ is smooth, then the number $\xi(N',N)$ is independent of $N$.
\end{enumerate}
Furthermore, if $Y$ is smooth, proper and Gorenstein, then  $\dim \Ext^i(N,N')$ is upper semi-continuous  and $\xi(N,N')$ continuous as functions of $N$. 
\end{prop}
\textbf{Proof.} The last assertion follows from i),ii) and proposition~\ref{pshiftflat}, for BK-Serre duality shows (semi)-continuity in the contravariant variable follows from (semi)-continuity in the covariant variable. Part ii) is an immediate consequence of the next result. In the course of its proof, we will also establish i). 

\begin{claim}
Let $\rk$ denote the rank of a noetherian $R$-module. Then 
\[ \xi(N',N) = \sum_i (-1)^i \rk \Ext^i(N',\caln)   .\]
\end{claim}
\textbf{Proof.} Given a flat $R$-algebra $S$, [AZ01, proposition~C3.1i)] implies that $\Ext^i(N',\caln \otimes_R S) \simeq \Ext^i(N',\caln)\otimes_R S$. We may thus assume that $R$ is a discrete valuation ring so that its uniformising parameter, say $u$, corresponds to the fibre $N \in \caln$. Flatness of $\caln$ ensures an exact sequence of the form
\[  0 \lm \caln \xrightarrow{\ u \ } \caln \lm N \lm 0  .\]
Applying $\Rhom(N',-)$ to this gives the long exact sequence 
\[ \ldots \lm \Ext^{i-1}(N',N) \lm \Ext^i(N',\caln) \xrightarrow{\ u \ } \Ext^i(N',\caln) \lm \Ext^i(N',N) \lm \Ext^{i+1}(N',\caln) \xrightarrow{\ u \ } \ldots \]
Ext-finiteness implies $E_i := \Ext^i(N',\caln)$ is a noetherian $R$-module. Note that $\Ext^i(N',N)$ is an extension of $\Hom_R(R/u,E_{i+1})$ by $E_i \otimes_R R/u$ so its dimension is upper semi-continuous as a function of $N$. This proves i).

The above sequence also shows 
\[ \xi(N',N) = \sum_i (-1)^i \dim \coker(E_i \xrightarrow{u} E_i) - \sum_i (-1)^i \dim \ker(E_i \xrightarrow{u} E_i)   .\]
The claim will follow when we show for arbitrary $E_i \in \mo R, i = 0,1, \ldots , d=\mbox{gl}.\dim Y$, that 
\[ (*) \hspace{1cm} \sum_i (-1)^i \dim \coker(E_i \xrightarrow{u} E_i) - \sum_i (-1)^i \dim \ker(E_i \xrightarrow{u} E_i) = \sum_i (-1)^i \rk E_i    .\]
As functions of the $E_i$, both sides of (*) take direct sums to sums so we need only check (*) for one non-zero $E_i$ which is either free or of the form $R/u^jR$. If $E_i$ is free, $\dim \coker(E_i \xrightarrow{u} E_i) = \rk E_i$ whilst $\dim \ker(E_i \xrightarrow{u} E_i) = 0$. Thus (*) holds in this case. If $E_i = R/u^jR$ then $\dim \coker(E_i \xrightarrow{u} E_i) =1 = \dim \ker(E_i \xrightarrow{u} E_i)$ so (*) holds again. This completes the proof of the claim and hence the proposition. 

\vspace{2mm}

\vspace{2mm}
\section{Bertini type theorems}   \label{sbertini}   

In this section, we continue our study of the generic behaviour of modules in a flat family. One version of the classic Bertini theorem states that given a morphism of projective varieties which is it own Stein factorisation, the fibres of the morphism are irreducible generically. We seek similar results in a non-commutative setting. Throughout, $Y$ will always be a strongly noetherian Hom-finite quasi-scheme.

We recall Artin-Zhang's Nakayama lemma [AZ01, theorem~C3.4]. Let $X$ be a noetherian scheme and $\calm\in \mo Y_X$. Then there is an open subset $U \subseteq X$ such that for any morphism $g:V \lm X$, we have $g^*\calm = 0$ if and only if $g$ factors through $U$.

\begin{defn} \label{dsupp}  
The closed set $X-U$ is called the {\em support of $\calm$ in $X$} and is denoted $\supp_X \calm$. 
\end{defn}

\begin{lemma}  \label{lmap0closed}  
Let $X$ be a scheme of finite type and $\phi: \caln \lm \calm$ be a non-zero morphism of objects in $\mo Y_X$ where $\calm$ is flat over $X$. The set $X_0\subseteq X$ of $x \in X$ where $\phi \otimes_X k(x): \caln \otimes_X k(x) \lm \calm \otimes_X k(x)$ is zero is a proper closed subset of $X$.
\end{lemma}
\textbf{Proof} by d\'{e}vissage. Thus we may assume $X$ is integral. Since $\otimes$ is right exact, we may replace $\caln$ with $\im \phi$ and so assume we have an exact sequence of the form 
$$0 \lm \caln \lm \calm \lm \mathcal{C} \lm 0 ,$$
where $\caln \neq 0$ by assumption. We claim $\supp_X \caln = X$. If this is not the case then there is some affine open $\spec R = U \subseteq X$ and a non-zero-divisor $r \in R$ with $\caln(U) \neq 0$ but $\caln(U) \otimes_U R[r^{-1}] = 0$ so $r^n \caln(U) = 0$ for $n \gg 0$. But $\calm/X$ is flat so this contradicts the fact that multiplication by $r^n$ induces an injection on $\calm(U)$. 

Now generic flatness of $\mathcal{C}$ [AZ01, theorem~C5.2] ensures by way of [AZ01, proposition~C1.4(i)] that there is a non-empty open set $V \subseteq X$ where $\caln \otimes_X k(x) \hookrightarrow \calm \otimes_X k(x)$ injects for all $x \in V$. Hence the locus where $\phi\otimes_X k(x)$ is zero must be contained in  the closed set $X - V$ and we are done by d\'{e}vissage.
\vspace{2mm}

Let $\oy$ be an object in $\mo Y$. If this is regarded as the structure sheaf then quotients like $\oy/I, \oy/J$ can be interpreted as subschemes of $Y$. Containment of subschemes can be expressed algebraically then via the condition $I \subseteq J$ or equivalently, the composite $I \lm \oy \lm \oy/J$ is zero. The next result shows that containing a given subscheme is a closed condition. It also contains a Bertini type result. We say that a set is {\em quasi-closed} if it is the countable union of closed sets and a {\em quasi-open} set is the complement of such a set.

\begin{cor} \label{cbertini} 
Let $\oy \in \mo Y$ and $\calm/X$ a flat family of quotients of $\oy$ parametrised by a scheme $X$ of finite type. 
\begin{enumerate}
\item For any $I < \oy$ the condition on $x \in X$ for $\oy \lm \oy/I$ to factor through $\oy \lm \calm \otimes_X k(x)$ is closed in $X$. 
\item Supppose that $Y$ has Halal Hilbert schemes, a continuous exact dimension function (that is hypotheses~3i),iii) of section~\ref{sass}) and no shrunken flat deformations. Fix some integer $n$. The condition for a fibre of $\calm$ to be $n$-critical is a quasi-open condition on $X$.
\end{enumerate}
\end{cor}
\textbf{Proof.} We shall only prove ii) since i) is easier and uses the same method. We may assume that $X$ is connected. 

Suppose some closed fibre $\calm \otimes_X k(x)$ is not critical but has some non-trivial $n$-dimensional quotient say $N$. Let $Z$ be the irreducible component of $\Hilb \oy$ such that the universal family $\caln/Z$ contains $N$ as a closed fibre. By assumption, the dimension of all the closed fibres of $\caln$ are also $n$. Also, the hypothesis of no shrunken flat deformations means that the irreducible components of $\Hilb \oy$ corresponding to $\calm/X$ are disjoint from $Z$. Thus a fibre $M=\oy/K \in \calm$ will certainly fail to be $n$-critical if there exists some $z \in Z$ such that the composite $K \lm \oy \lm \caln \otimes_Z k(z)$ is zero. We say in this case that $M$ is {\em $Z$-uncritical}.

We claim that the condition to be $Z$-uncritical is closed in $X$. Indeed, let $\mathcal{K} = \ker(\oy \otimes_k \ox \lm \calm)$. Applying the lemma to the composite $\mathcal{K}\otimes_k \calo_Z \lm \oy \otimes_k \calo_{X \times Z} \lm \caln \otimes_k \ox$ we see that the points $(x,z) \in X \times Z$ where
$\mathcal{K} \otimes_X k(x) \lm \caln \otimes_Z k(z)$ is zero is a closed set $C \subseteq X \times Z$. Now $Z$ is projective by our Halal Hilbert scheme assumption, so the $Z$-uncritical fibres of $\calm$ correspond to the closed image of $C$ under the projection map $X \times Z \lm X$. The corollary follows since failing to be critical corresponds to being $Z$-uncritical for one of the countably many possible components of Hilbert schemes $Z$ above.

\textbf{Remark:} We do not know if being critical is generic as we have no way of bounding the possibilities of the components $Z$ of the Hilbert scheme above. When $Y$ is projective, one might be able to approach the question by first showing that Hilbert schemes corresponding to a fixed Hilbert function are projective and then trying to bound the number of possible Hilbert functions of the quotients $N$. 

\vspace{2mm}

Recall that a module $M \in \mo Y$ is {\em $d$-pure} if for every proper submodule $N < M$ we have $\dim N = d$. We next give a result that guarantees that members of a flat family $\calm$ are generically 1-pure. 

\begin{prop}  \label{pgenpure}  
Let $Y$ be a smooth proper Gorenstein quasi-scheme with a compatible dimension function and classical cohomology. Let  $\calm$ be a flat family of $Y$-modules parametrised by a quasi-projective smooth curve $X$. If for some fibre $M \in \calm$ we have $H^0(M \otimes \w_Y) = 0$, then only a finite number of $M' \in \calm$ have non-trivial zero-dimensional submodules.
\end{prop}
\textbf{Proof.} If $P\neq 0$ is a zero-dimensional submodule of $M'$ then the assumption on classical cohomology ensures $H^0(P\otimes \omega_Y) \hookrightarrow H^0(M'\otimes \omega_Y)$ is non-zero. However, upper semicontinuity (proposition~\ref{pcont}) and $H^0(M \otimes \omega_Y) = 0$ means that this can only occur for finitely many $M'\in \calm$.

\vspace{2mm}

One expects the above purity result to hold more generally. In the projective case, we can prove the following.

\begin{prop}  \label{pprojgenpure}  
Let $Y=\proj A$ be a non-commutative smooth projective $d$-fold and $\calm$ be a flat family of $Y$-modules parametrised by a quasi-projective smooth curve $X$. If some fibre $M \in \calm$ has no zero-dimensional submodules, then only a finite number of $M' \in \calm$ have  zero-dimensional submodules.
\end{prop} 
\textbf{Proof.} We may repeat the proof above using the fact that the graded shift $\mo Y \lm \mo Y: N \mapsto N(n)$ preserves dimension and the following 

\begin{lemma}  \label{lGabbermax}  
Let $A$ be a noetherian graded algebra with an Auslander dualising complex and $M$ be a noetherian graded $A$-module which has no submodule of canonical dimension $\leq 1$. Then $H^0(M(-n))=0$ for $n\gg 0$. 
\end{lemma}
\textbf{Proof.} We can find a filtration of $M$ by pure modules so we may assume that $M$ is $m$-pure with respect to the canonical dimension, for some $m>1$. We write $\mmm$ for the augmentation ideal $A_{>0}$. Let $N$ be the maximal submodule of the injective hull $E(M)$ such that $N\geq M$ and $N/M$ is $\mmm$-torsion. Note that $N$ is also $m$-pure so by the Gabber maximality principle [YZ, theorem~2.19], $N$ is finitely generated and, in particular, left bounded. This completes the proof.

\vspace{2mm}
\vspace{2mm}
\section{Fibration over a Curve}  \label{sfibr}   

In this section, we introduce a method for constructing flat morphisms $f: Y \lm X$ from a strongly noetherian Hom-finite quasi-scheme $Y$ to a commutative projective curve. By a {\em flat morphism}, we mean a pair of adjoint functors $f^*:\Mod X \lm \Mod Y, f_*: \Mod Y \lm \Mod X$ such that $f^*$ is exact and preserves noetherian objects.


We use a type of Fourier-Mukai transform with kernel a flat family $\calm/X$ of objects in $\mo Y$ parametrised by a quasi-projective scheme $X$. Recall from section~\ref{sbbchange} that we have an exact functor $\calm \otimes_X -: \mo X \lm \mo Y_X$. We thus need to construct an auxiliary projection functor $\pi_*: \Mod Y_X \lm \Mod Y$ and its derived functors.  These are defined via ``relative Cech cohomology''. Let $\calu = \{U_i\}$ be a finite open affine cover of $X$ and $\caln \in \mo Y_X$. Given an open affine $U \subseteq X$ we let $\caln(U)$ denote the restriction of $\caln$ to $U$ as usual. We can define the usual Cech cohomology complex $C^\bullet (\calu,\caln)$ in $\Mod Y$ by 
$$ C^p(\calu,\caln) = \bigoplus_{i_0< \ldots < i_p} \caln(U_{i_0} \cap \ldots \cap U_{i_p})$$
for $p\geq 0$ and the usual differentials. We can thus define $R^p\pi_{\calu *}\caln := H^p(C^{\bullet}(\calu,\caln))\in \Mod Y$. If $X$ is affine then descent theory (see [AZ01, corollary~C8.2]) shows $R^p\pi_{\calu *}\caln =0$ for $p>0$ and $\pi_{\calu *} \caln:= R^0\pi_{\calu *}\caln = \caln$. As usual, this shows that $R^p\pi_{\calu *}$ is independent of $\calu$ so we may define $R^p\pi_*=R^p\pi_{\calu *}$. Indeed, let $\mathcal{V} = \{ V_j \}$ be another finite open affine cover of $X$. We can form the double Cech complex $C^\bullet(\calu,\mathcal{V},\caln)$ with 
$$ C^{p,q}(\calu,\mathcal{V},\caln) = \bigoplus \caln(U_{i_0} \cap \ldots \cap U_{i_p}\cap V_{j_0} \cap \ldots \cap V_{j_q})$$
where the direct sum is over indices satisfying $i_0< \ldots < i_p$ and $j_0< \ldots < j_q$. The columns are just direct sums of Cech complexes for $\caln$ restricted to $U_{i_0} \cap \ldots \cap U_{i_p}$ and the rows are direct sums of Cech complexes for $\caln$ restricted to $V_{j_0} \cap \ldots \cap V_{j_q}$. Since on $X$, the intersections of affines is affine, the rows and columns have cohomology only in degree 0 so the total complex abuts to both $C^\bullet(\calu,\caln)$ and $C^\bullet(\mathcal{V},\caln)$. This verifies that $R^p\pi_*$ is indeed well-defined. 

Finally, we can now define
\[  f^*: \mo X \lm \Mod Y: \mathcal{L} \mapsto \pi_*(\calm \otimes_X \mathcal{L})    .\]
We observe that $f^*$ is left exact and wish to show that, under suitable hypotheses, it is also right exact and preserves noetherian objects. This will follow from analogous statements about $\pi_*$. 

We start by recovering standard Serre theory of cohomology on projective schemes in this context.
Cohomology of some sheaves are easy to compute.

\begin{prop}  \label{pcohgen}  
Let $M \in \Mod Y, \mathcal{F} \in \Mod X$. Then 
$$R^i\pi_* (M \otimes_k \mathcal{F}) = M \otimes_k H^i(X,\mathcal{F}) .$$
\end{prop}
\textbf{Proof.} Let $\mathcal{U}$ be a finite affine open cover of $X$. We have an isomorphism of Cech complexes 
$$ C^{\bullet}(\mathcal{U},M \otimes_k \mathcal{F}) \simeq M \otimes_k C^{\bullet}(\mathcal{U}, \mathcal{F}) .$$
Since $M\otimes_k-$ is exact, the proposition follows.

\vspace{2mm}

Below is the usual Serre vanishing theorem in our context. We consider an ample line bundle $\ox(1)$ and, as usual write $\caln(n)$ for $\caln \otimes_X \ox(n)$.

\begin{thm}  \label{tserrevan}  
Let $X$ be a projective scheme and $\caln \in \mod Y_X$. Then i) $R^i\pi_*\caln$ is noetherian for all $i$ and, ii) $R^i\pi_*\caln(n) = 0$ for all $i>0, n \gg 0$.
\end{thm}
\textbf{Proof.} We check briefly that the standard proof works here. By proposition~\ref{pcohgen}, the theorem holds for any finite direct sum $\mathcal{F}$ of modules of the form $M\otimes_k \ox(j)$ where $M \in \mo Y$. Now proposition~\ref{pgenYX} ensures that $\caln$ is the quotient of such an $\mathcal{F}$. The result now follows from downward induction on $i$ and the long exact sequence in cohomology. 

\vspace{2mm} 
%
%

\begin{prop}  \label{pr1pi}  
Let $X$ be a projective curve, $\ox(1)$ an ample line bundle and let $H \subset X$ be the zeros of some non-zero section. Suppose $\caln \in \mo Y_X$ which is flat in a neighbourhood of $H$. Then $R^1\pi_* \caln$ has a finite filtration with factors isomorphic to quotients of $\caln \otimes_X \calo_H$.
\end{prop}
\textbf{Proof.} For $n \in \Z$ we consider the exact sequence
$$ 0 \lm \caln(n-1) \lm \caln(n) \lm \caln \otimes_X \calo_H \lm 0 $$ 
and the associated long exact sequence in cohomology
$$\pi_*(\caln \otimes_X \calo_H)  \lm R^1 \pi_*\caln(n-1) \lm R^1 \pi_*\caln(n)  \lm R^1\pi_*(\caln \otimes_X \calo_H) .$$
Choosing an affine open cover judiciously to compute Cech cohomology, we find that
$$ \pi_*(\caln \otimes_X \calo_H)  = \caln \otimes_X \calo_H\ ,\ \  R^1\pi_*(\caln \otimes_X \calo_H) =0.$$
The result follows by downward induction on $n$ and Serre vanishing theorem~\ref{tserrevan}.

\vspace{2mm} 

Let $N \in \mo Y$. Then we let $SQ(N)$ denote the set of simple quotients of $N$, that is, 
$$ SQ(N) = \{ S \in \mo Y| S\  \mbox{is simple and }\ \Hom(N,S) \neq 0\}.$$ 
This will be our replacement for the notion of the support of $N$. It is non-empty precisely when $N$ is non-zero since we are considering noetherian objects. 

\begin{cor}  \label{csq}  
\begin{enumerate}
\item Given an exact sequence in $\mo Y$ 
$$ 0 \lm N' \lm N \lm N'' \lm 0 $$
we have $SQ(N'') \subseteq SQ(N) \subseteq SQ(N') \cup SQ(N'')$.
\item Let $\caln \in \mo Y_X$ where $X$ is a projective curve. Then for every smooth closed point $p \in X$ such that $\caln$ is flat in a neighbourhood of $p$, we have $SQ(R^1\pi_* \caln) \subseteq SQ(\caln \otimes_X k(p))$. In particular, $R^1\pi_* \caln = 0$ if for every simple $S\in \mo Y$, there is a smooth point $p \in X$ with $\Hom(\caln \otimes_X k(p),S) = 0$.
\end{enumerate}
\end{cor}
\textbf{Proof.}
We omit the easy proof of part i) and proceed with ii). If $p$ is smooth, then it is the zero of some section of the ample line bundle $\ox(p)$. The corollary follows from proposition~\ref{pr1pi} and part i).

\vspace{2mm}

The obstruction to defining a morphism is given by the following non-commutative version of base points.
\begin{defn}  \label{dbasepoint}  
Let $X$ be a projective curve and $\calm \in \mo Y_X$, a flat family of $Y$-modules. A simple object $P \in \mo Y$ is a {\em base point} of $\calm$ if $\Hom(M',P) \neq 0$ for generic $M' \in \calm$. We say that $\calm$ is {\em base point free} if there are no base points.
\end{defn}

If $Y$ is smooth, proper and Gorenstein, then the semicontinuity results of section~\ref{sscont} hold, so to check that $\calm$ is base point free, it suffices to show that for any simple $S$, there is some smooth point $p \in X$ with $\Hom(\calm \otimes_X k(p),S) = 0$.

We can finally show right exactness of $f^*$. 
\begin{thm}   \label{tpiexact}   
Let $X$ be a projective curve which is generically smooth and $\calm \in \mo Y_X$ be a flat family of objects in $\mo Y$. Suppose that $\calm$ is base point free. Then the functor $f^*: \mo X \lm \mo Y$ is exact and so has a right adjoint $f_*$. Hence $f:Y \lm X$ is a flat morphism. 
\end{thm}
\textbf{Proof.} The long exact sequence in cohomology and adjoint functor theorem reduces the proof to showing 
\[ (*) \hspace{1cm} R^1\pi_*(\calm \otimes_X \call) = 0 \hspace{1cm} \] 
for any $\call \in \mo X$. If it is non-zero then there exists a simple object $S$ in $SQ(R^1\pi_*(\calm \otimes_X \call))$ so for every smooth point $p \in X$, corollary~\ref{csq}ii) ensures $\Hom(\calm\otimes_X \call \otimes _X k(p),S) \neq 0$. Thus $S$ is a basepoint which contradicts our assumption. The theorem is now proved. 

\begin{defn} \label{dFM}  
We call the morphism $f: Y \lm X$ in the theorem, the {\em Fourier-Mukai morphism associated to $\calm$}. 
\end{defn}

\begin{eg}  \label{eaffinemap} 
``Adjoint'' of closed embedding.
\end{eg}
Suppose $Y$ is a commutative smooth projective surface. Let $i:X \hookrightarrow Y$ be the closed embedding of a smooth curve $X$ and $\calm = \ox$, the structure sheaf of the graph of $i$. Then $\calm$ is base point free and the resulting morphism $f:Y \lm X$ is given by $f^* = i_*, f_* = i^!$. 
\vspace{2mm}

\noindent 
\textbf{Remark:} Given a base point free family $\calm/X$ of $Y$-objects and a non-constant morphism $X' \lm X$ of projective curves, the pull-back $\calm'/X'$ is also base point free. Thus there is also a map $f':Y \lm X'$. This should not be surprising since as in the previous example, given any affine morphism $X' \lm X$, there is an ``adjoint'' morphism of non-commutative schemes $X \lm X'$.

\vspace{2mm}
\section{Families from a flat morphism}  \label{sconverse}  

In this section, we give the converse construction to the one in the previous section. Consider a noetherian quasi-scheme $Y$, a commutative projective variety $X$ and a flat morphism $f:Y \lm X$. We show that it comes from a Fourier-Mukai transform as in section~\ref{sfibr}.

Given $\call \in \mo X$ and $U \subseteq X$ an affine open, we can consider $\call \otimes_X \calo_U$ as a quasi-coherent sheaf on $X$. We wish to construct an induced $X$-object $f^{\flat} \call$ as follows. Define 
$$(f^{\flat}\call)(U) = f^*(\call \otimes_X \calo_U) \in \Mod Y .$$
The $\calo_U$-action is given by 
$$\calo_U \lm \End_X (\call \otimes_X \calo_U) \xrightarrow{f^*} \End_Y (f^*(\call \otimes_X \calo_U)) .$$  
To verify compatiblity with restrictions, consider $V \subset U$ another affine open. Now $f^*$ is right exact and commutes with direct sums, so there is an isomorphism $f^*(\call \otimes_X \calo_V) \simeq  f^*(\call \otimes_X \calo_U) \otimes_U \calo_V$ which is canonical. We hence obtain a well-defined $X$-object $f^{\flat} \call$. Note that $f^{\flat} \call$ is naturally isomorphic to $(f^{\flat} \ox) \otimes_X \call$ since this holds for $X$ affine. Thus exactness of $f^{\flat}$ ensures that $f^{\flat}\ox$ is flat over $X$. 
\begin{prop}  \label{pinducednoeth} 
If $\call \in \Mod X$ is noetherian, then so is $f^{\flat} \call$.
\end{prop}
\textbf{Proof.} Let $H$ be an ample hyperplane section and $T = \oplus_{n} H^0(X, \ox(nH))$ be the corresponding homogeneous co-ordinate ring of $X$. Let $t \in H^0(X,\ox(H))$ correspond to $H$ so that $U:= X-H = \spec T[t^{-1}]_0$. It suffices to show that $f^{\flat} \call (U)$ is noetherian. We know $T[t^{-1}]_0$ is finitely generated, say by $T_mt^{-m}$. We will consider $\calo_U$ as the union of the sheaves $\ox(nH), n \in \mathbb{N}$. Then for $n \geq m$ we see that multiplication in $\calo_U = T[t^{-1}]_0$ induces surjections $\ox(nH) \otimes_k T_m t^{-m} \lm \ox((n+m)H)$ and hence surjections $f^*(\call \otimes_X \ox(nH)) \otimes_k T_m t^{-m} \lm f^*(\call \otimes_X \ox((n+m)H))$. This in turn gives a surjection $f^*(\call \otimes_X \ox(mH)) \otimes_k \calo_U \lm f^*(\call \otimes_X \calo_U)$. Now $f$ is flat, so by definition $f^*(\call \otimes_X \ox(mH))$ is noetherian and we are done by [AZ01, proposition~B5.1].
\vspace{2mm}

Let $\mathcal{U} = \{U_i\}$ be a finite affine open cover of $X$ with which we shall compute Cech cohomology. There is a Cech complex $C^{\bullet}(\mathcal{U},\call)$ of quasi-coherent sheaves on $X$ with 
$$ C^p(\mathcal{U},\call) = \bigoplus_{i_0< \ldots < i_p} \call \otimes_X \calo_{U_{i_0 \ldots i_p}} $$
where $U_{i_0 \ldots i_p} = U_{i_0} \cap \ldots \cap U_{i_p}$. The complex $C^{\bullet}(\mathcal{U},\call)$ only has cohomology in degree 0, and $H^0(C^{\bullet}(\mathcal{U},\call)) = \call$. Now $f^*$ is exact so we have 
$$H^0(C^{\bullet}(\mathcal{U},f^{\flat}\call)) \simeq f^* \call$$
from which follows
\begin{prop}  \label{pFMfrommap}  
Given a flat morphism $f:Y \lm X$ from a noetherian quasi-scheme to a commutative projective variety, we have a natural isomorphism of functors $f^* \simeq \pi_*(f^{\flat}\ox \otimes_X -)$.
\end{prop}

\section{$K$-non-effective rational curves of self-intersection 0} \label{smoduli}  

Recall that our aim in this part is to generalise theorem~\ref{tcommmoricon} about contracting a $K$-negative curve $C$ on a surface when $C^2 = 0$. In this section, we will define analogues of $\calo_C$ for non-commutative smooth proper surfaces. We wish to use Fourier-Mukai transforms as in section~\ref{sfibr} to obtain contractions. The requisite kernels will come from studying the Hilbert scheme. We start with a general result about Hilbert schemes, so $Y$ will denote a strongly noetherian quasi-scheme with Halal Hilbert schemes. We will specialise to non-commutative smooth proper surfaces later.

The next theorem nicely generalises the classical deformation theory of Hilbert schemes. 
\begin{thm}  \label{tobssmooth} 
Let $Y$ be a strongly noetherian Ext-finite quasi-scheme with Halal Hilbert schemes. Let 
\[0 \lm G \lm E \lm F \lm 0 \]
be an exact sequence in $\mo Y$ defining a closed point $p$ of the Hilbert scheme $\Hilb(E)$. Then the tangent space to $\Hilb$ at $p$ is $\Hom(G,F)$ and the Hilbert scheme is smooth at $p$ if $\Ext^1(G,F) = 0$.
\end{thm}
\textbf{Proof.} [AZ01, proposition~E2.4] gives the tangent space given above. It also gives an obstruction to smoothness, which we need to relate to our obstruction above. To recall their obstruction, let $R'$ be an artinian local $k$-algebra with residue field $k$ and maximal ideal $\mmm$. Let $I\triangleleft R'$ be an ideal with $\mmm I = 0$ and $R = R'/I$. Consider an $R$-flat deformation of $G \lm E \lm F$, say given by the exact sequence
$$ 0 \lm \Gt \lm E \otimes_k R \lm \Ft \lm 0 .$$
The obstruction in [AZ01, proposition~E2.4i)] to lifting $\Ft$ to $R'$ is given by an element $\eta \in \Ext^1_{Y_{R'}}(\Gt,F\otimes_k I)$. We are required to show that $\Ext^1(G,F) = 0$ ensures $\eta$ vanishes. We compute the obstruction space using the following result which, in the classical case can be obtained from a spectral sequence.

\begin{lemma}    \label{l3termexseq}  
Let $U \in \mo Y_{R'}$ and $V_0 \in \mo Y$. There is an exact sequence 
$$ 0 \lm \Ext^1_Y(U \otimes_{R'} k,V_0) \xrightarrow{\ \phi \ } \Ext^1_{Y_{R'}}(U,V_0) \xrightarrow{\ \psi\ } 
\Hom_Y(\Tor^{R'}_1(U,k),V_0)  $$
\end{lemma}
\textbf{Proof.} The map $\phi$ above can be defined by pull-back of extensions as in the commutative diagram with exact rows below
$$\diagram
\phi(E): \ 0 \rto & V_0 \rto \ddouble & L \dto \rto & U \rto\dto & 0 \\
E: \ 0 \rto & V_0 \rto & L_0 \rto & U \otimes_{R'} k \rto & 0 
\enddiagram$$
We may tensor this entire diagram with $k$ to obtain another commutative diagram with exact rows
$$\diagram
\phi(E)\otimes_{R'} k: &  & V_0 \rto^(.3){\iota} \ddouble & L\otimes_{R'} k \dto \rto & U\otimes_{R'} k \rto\dto & 0 \\
E: \ & 0 \rto & V_0 \rto & L_0 \rto & U \otimes_{R'} k \rto & 0 
\enddiagram.$$
Exactness of the bottom row ensures that $\iota$ is injective so the 5-lemma shows that $E = \phi(E) \otimes_{R'} k$.  We see then that if $\s:L \lm V_0$ splits $\phi(E)$ then $\s \otimes_{R'} k$ splits $E$. Thus $\phi$ is injective. To define $\psi$, consider the exact seqence in $\mo Y_{R'}$
$$ E': \ 0 \lm V_0 \lm L \lm U \lm 0 .$$
The long exact sequence in tor gives $\psi(E')$ as the connecting homomorphism below 
$$ \Tor^{R'}_1(U,k) \xrightarrow{\psi(E')} V_0 \lm L \otimes_{R'} k \lm U \otimes_{R'} k \lm 0  .$$
Note $\psi(E') = 0$ precisely when 
$$E' \otimes_{R'} k: 0 \lm V_0 \lm L \otimes_{R'} k \lm U\otimes_{R'} k \lm 0 $$
is exact. If $E'=\phi(E)$ then $E' \otimes_{R'} k = E$ which is exact so $\psi \circ \phi = 0$. Conversely if $E' \otimes_R k$ is exact then it is easy to see that $E' = \phi(E' \otimes_{R'} k)$. This completes the proof of the lemma.

\vspace{2mm}

To compute the tor term above we need
\begin{lemma}  \label{lcomputetor} 
Let $P=k[x_1,\ldots,x_r],R' = P/P_{>n}, R = P/P_{\geq n}$ so $I = P_{\geq n}/P_{>n}$. Suppose $U$ is an $R$-flat module in $\mo Y_R$. Then there are canonical isomorphisms
$$ \Tor^{R'}_1(U,k) \simeq U \otimes_R \Tor^{R'}_1(R,k) \simeq (U \otimes_R k) \otimes_k I  .$$
\end{lemma}
\textbf{Proof.} We only prove the case $r=2$, the general case is proved using the same methods with only more complicated bookkeeping required. We will only need the case $r=1$ in this paper. Consider the partial resolution of the $R'$-module $k$,  
$$ R'^{n+3} \xrightarrow{g_2} R'^2 \xrightarrow{g_1} R' \lm k \lm 0 $$
where 
$$g_2 = 
\begin{pmatrix}
 x_2 & x_1^n & x_1^{n-1}x_2 & \ldots & x_2^n & 0 \\
-x_1 & 0 & 0 & \ldots & 0 & x_2^n
\end{pmatrix}, \quad g_1 = (x_1 \ \ x_2).$$
Tensoring with $R$ over $R'$ we find the complex 
$$ C_{\bullet}: R^{n+3} \xrightarrow{g_2\otimes_{R'} R} R^2 \xrightarrow{g_1\otimes_{R'} R} R$$
has cohomology 
$$H:= 
k \begin{pmatrix}
   x_1^{n-1} \\ 0 
  \end{pmatrix} \oplus 
k \begin{pmatrix}
   x_1^{n-2}x_2 \\ 0 
  \end{pmatrix} \oplus \ldots \oplus
k \begin{pmatrix}
   x_2^{n-1} \\ 0 
  \end{pmatrix} \oplus 
k \begin{pmatrix}
   0 \\x_2^{n-1} 
  \end{pmatrix}$$
Note that $H$ is naturally isomorphic to $I$ via $g_1$. Also $U$ is flat over $R$ so tensoring with $C_{\bullet}$ shows 
$$ \Tor^{R'}_1(U,k) \simeq U \otimes_R I \simeq (U \otimes_R k) \otimes_k I  .$$
\vspace{2mm}

We return now to the proof of the theorem. Assume that $\Ext^1_Y(G,F) = 0$ so applying lemma~\ref{l3termexseq} to $U = \Gt,V_0 = F \otimes_k I$, we see it suffices to show that $\psi(\eta) = 0$. Let $P$ be the polynomial ring over $k$ in $\dim \Hom_Y(G,F)$ variables. To prove smoothness, it suffices to consider the case where $R' = P/P_{>n}, R = P/P_{\geq n}$.

We recall the definition of $\eta$ in [AZ01, proof of proposition~E2.4i)]. Let $\eta' \in \Ext^1_{Y_{R'}}(E \otimes_k R,E \otimes_k I)$ be given by the extension
$$ 0 \lm E \otimes_k I \lm E \otimes_k R' \lm E \otimes_k R \lm 0 .$$
Then $\eta$ is the image of $\eta'$ under the natural maps of Ext spaces
$$\Ext^1_{Y_{R'}}(E \otimes_k R,E \otimes_k I) \lm \Ext^1_{Y_{R'}}(\Gt,E \otimes_k I) \lm \Ext^1_{Y_{R'}}(\Gt,F \otimes_k I) .$$
Note that the maps of Ext spaces are induced by pull-back and push-forward of extensions. 

From the proofs of the lemmas above, we see that $\psi(\eta')\in \Hom_Y(\Tor^{R'}_1(E\otimes_k R,k),E \otimes_k I)$ is essentially the identity map 
$$ E \otimes_k I \xrightarrow{\ \ \sim\ \ } \Tor^{R'}_1(E\otimes_k R,k) \lm E \otimes_k I .$$
Now $\psi$ was defined as a connecting homomorphism, which is natural with respect to morphisms of extensions so $\psi(\eta)$ can be obtained from $\psi(\eta')$ as the composite
$$ \Tor^{R'}_1(\Gt,k) \simeq G \otimes_k I \hookrightarrow E \otimes_k I \xrightarrow{\psi(\eta')=\mbox{id}} E\otimes_k I \lm F \otimes_k I .$$
This is evidently zero so the theorem is proved. 

\vspace{2mm}

For the rest of this section, we restrict our attention to a non-commutative smooth proper surface $Y$. Recall that, as part of the data, we have a structure sheaf $\oy \in \mo Y$. Also, the (semi-)continuity results of section~\ref{sscont} are available to us. 

\begin{defn}  \label{dkneffint0}
Let $Y$ be a non-commutative smooth proper surface. A {\em curve} on $Y$, is a 1-pure module $M \in \mo Y$ which is a quotient of $\oy$. Furthermore, we say $M$ 
\begin{enumerate}
\item is {\em rational} if $h^0(M) = 1$ and  $h^1(M) = 0$. 
\item is {\em $K$-non-effective} if $h^0(M \otimes \w_Y) = 0$. 
\item {\em has self-intersection zero} if $M^2:= M.M = -\xi(M,M) = 0$. 
\end{enumerate}
\end{defn}
We will refrain from automatically assuming that the curve $M$ is 1-critical, which corresponds in the commutative case to an integral curve. Note that if $h^0(M) = 1$, the quotient map $\oy \lm M$ is unique up to scalar multiplication. Note also that we have not defined $K$-negative, but instead a related condition, that of $K$-non-effective, which does not appear in the commutative theory. Given a curve $C$ on a commutative smooth projective surface, $K$-negative curves are those with $\deg (\calo_C \otimes_Y \w_Y) < 0$ and this condition certainly implies the $K$-non-effective condition above. Unfortunately, we do not know if such an implication holds in the non-commutative case. The next result will show how it is the $K$-non-effective condition which will allow us to compute the flat deformations of $M$.

\begin{prop}  \label{pextmmp} 
Let $Y$ be a non-commutative smooth proper surface, $\calm/X$ be a flat family of quotients of  $\oy$ parametrised by a smooth quasi-projective curve $X$ and $M,M' \in \calm$. 
\begin{enumerate}
\item If $H^0(M\otimes \w_Y) = 0$ then $\Ext^2(M,M')=0$. 
\item If $H^0(M\otimes \w_Y) = 0$ and  $M^2=0$ then $\dim\Hom(M,M') =\dim \Ext^1(M,M')$.
\item  Suppose that $M$ is 1-critical and $M'$ 1-pure. Then $\Hom(M,M')= k$ if $M\simeq M'$ and is 0 otherwise. 
\item  Suppose that $h^0(M')=1$. Then $\Hom(M,M')= 0$ if $M,M'$ are non-isomorphic and is $k$ otherwise. 
\end{enumerate}
\end{prop}
\textbf{Proof.} 
First observe that $\Ext^2(M,M')^*=\Hom(M',M\otimes \w_Y)\subseteq H^0(M \otimes \w_Y) = 0$ so i) holds and ii) follows from $M.M'=0$ and proposition~\ref{pcont}. We now prove iii) and iv) and suppose that $\Hom(M,M') \neq 0$. If $M$ is 1-critical and $M'$ 1-pure, we have an embedding $M \hookrightarrow M'$. This must be an isomorphism by our assumption that there are no shrunken flat deformations. This also shows any non-zero endomorphism of $M$ must be an isomorphism. Thus $\Hom(M,M')$ must be a division ring that is finite dimensional over $k$. This can only be $k$ so iii) is proved. Suppose now instead that $h^0(M')=1$. Then any non-zero morphism $\phi:M \lm M'$ must be surjective as the composite $\oy \lm M \lm M'$ must be the unique, up to scalar, global section of $M'$. Our assumption that there are no shrunken flat deformations ensures that $\phi$ is an isomorphism. 
\vspace{2mm}

\begin{defn}  \label{dHilbsys} 
Let $Y$ be a non-commutative smooth proper surface and $M$ be a quotient of $\oy$. Let $X$ be the union of irreducible components of the Hilbert scheme $\Hilb \oy$ which contain the point corresponding to $M$. Then the {\em Hilbert system} of $M$ is the universal quotient $\calm/X$.
\end{defn}

The key result of this section is the following. 

\begin{cor}\label{cfamily} 
Let $Y$ be a non-commutative smooth proper surface and $\calm/X$, the Hilbert system of a $K$-non-effective rational curve $M$ with self-intersection zero. Then $X$ is a projective curve which is smooth at the point corresponding to $M$.
\end{cor}
\textbf{Proof.} Consider the exact sequence
\[ 0 \lm G \lm \oy \lm M \lm 0 .\]
Applying the functor $\Hom(-,M)$ yields the long exact sequence
\[ k= \Hom(M,M) \lm \Hom(\oy,M) \lm \Hom(G,M) \lm \Ext^1(M,M)  \]
\[  \lm H^1(M)  \lm \Ext^1(G,M) \lm \Ext^2(M,M) \lm H^2(M) = 0   \]
where the last equality follows from corollary~\ref{cvanish}. Rationality of $M$ and proposition~\ref{pextmmp} ensure that $\Hom(G,M) = \Ext^1(M,M) = k, \Ext^1(G,M) = 0$ so the result follows from theorem~\ref{tobssmooth}. 
\vspace{2mm}

\section{Mori contraction}  \label{smoricon} 

In this section, we wish to contract a 1-critical $K$-non-effective rational curve of self-intersection zero using a Fourier-Mukai morphism. From theorem~\ref{tpiexact}, we need its Hilbert system to be base point free. However, S. P. Smith has examples (see [SV]) which show that this will not always hold. However, for those examples, the base points are non-zero in the Picard group, that is, they are not in the radical of the intersection pairing. In fact, they have non-zero self-intersection. We will see that this is the only obstruction to base point freedom. In this section, $Y$ denotes a non-commutative smooth proper surface.

\begin{defn}  \label{dpt} 
A {\em point of $Y$} is a simple zero-dimensional $Y$-module $P$. We say {\em $P$ lies on $M\in \calm$} if $P \in SQ(M)$. 
\end{defn}

\begin{thm}  \label{tptmove} 
Let $\calm/X$ be a non-trivial flat family of quotients of $\oy$ parametrised by a projective curve $X$. Suppose one of the fibres $M \in \calm$ is a 1-critical $K$-non-effective rational curve with $M^2=0$ and that the corresponding point of $X$ is smooth. Suppose further that for every point $P$ on $M$ we have $\xi(M,P) \leq 0$. Then $\calm$ is base point free and so induces a Fourier-Mukai morphism $f: Y \lm X$. 
\end{thm}
\textbf{Proof.} 
Suppose to the contrary that $P\in \mo Y$ is a base point. If $\Hom(M,P)=0$ then semicontinuity ensures $\Hom(M',P)=0$ for generic $M' \in \calm$ so $P$ is not a base point. Also, $P=M$ cannot be a base point either as otherwise, for generic $M' \in M$ we obtain non-zero maps $M' \lm M$ which are isomorphisms since there are no shrunken flat deformations. Furthermore, $H^0(M)=k$ so the family must then be the trivial family, a contradiction. Thus $\Hom(M,P) \neq 0$ and $P$ is zero dimensional. Arguing by induction, we may assume that $\chi(P) = h^0(P)$ is minimal amongst all base points. 

Now $\Ext^2(M,P) = \Hom(P,M \otimes \w_Y)^*=0$, so $\xi(M,P) \leq 0$ ensures $\Ext^1(M,P)\neq 0$. Also, BK-Serre duality shows $\Ext^1(P,M \otimes \w_Y) \neq 0$ so there exists a non-split exact sequence
$$0 \lm M \otimes \w_Y \lm L \lm P \lm 0$$
where $L$ is 1-critical since $M \otimes \w_Y$ is. For generic $M' \in \calm$, applying $\Rhom(M',-)$ to this sequence and using proposition~\ref{pextmmp}, proposition~\ref{pgenpure} and BK-Serre duality shows that generically we have $\Hom(M',L) \neq 0$. Thus $\Hom(M,L) \neq 0$ by semicontinuity (proposition~\ref{pcont}) and we may construct an exact sequence 
$$ 0 \lm M \lm L \lm C  \lm 0  .$$ 
Note that $C$ is zero dimensional and 
$$h^0(C) = \chi(P) + \chi(M\otimes  \w_Y) - \chi (M) = h^0(P) -h^1(M \otimes \w_Y) - 1 < h^0(P) .$$
Minimality of $h^0(P)$ ensures that no simple subquotient of $C$ is a base point so generically we have $\Hom(M',C) = 0$. This means generically $\Hom(M',M) \neq 0$. As in the $P=M$ case, this contradicts the fact that the family $\calm$ is non-trivial and so proves the theorem. 

\vspace{2mm}

We have a converse result which relates base point freedom to the condition on intersection with points.

\begin{prop}  \label{pbptfree} 
Let $\calm/X$ be a flat family of noetherian $Y$-modules parametrised by a generically smooth projective curve. Suppose $\calm$ is base point free and that a generic fibre $M' \in \calm$ is 1-pure. Then for any point $P \in \mo Y$ and $M \in \calm$ we have $\xi(M,P)\leq 0$.
\end{prop}
\textbf{Proof.} Generically we have $\Hom(M',P) = 0$. Also $\Ext^2(M',P) = \Hom(P,M' \otimes \w_Y)^*$ which is generically 0 since $M' \otimes \w_Y$ is also 1-pure. Continuity of intersection numbers then shows $\xi(M,P) \leq 0$ for all $M \in \calm$.

\begin{defn}  \label{dmoricon}  
Let $\calm/X$ be a base point free Hilbert system of a $K$-non-effective rational curve $M$ with $M^2=0$. Then the associated morphism $f:Y \lm X$ is called the {\em non-commutative Mori contraction contracting $M$}.
\end{defn}

\part{Non-commutative $\PP^1$-bundles}  \label{Pncruled} 

The notion of non-commutative Mori contractions will hopefully be useful in understanding Artin's conjecture. In this part, we show how non-commutative ruled surfaces give non-trivial examples of the theory developed in part~II. 

Many of our results will hold more generally for non-commutative $\PP^1$-bundles so in this part, we will usually work in this setting. Non-commutative $\PP^1$-bundles were introduced in [VdB01p], and the version of his definition we shall follow, will be the one given in [Na, definition~2.4]. Here, we will only provide a brief description of them, mainly to fix notation for the rest of part~\ref{Pncruled}. 

Recall first Van den Bergh's notion of an $\ox$-bimodule. Let $R$ be a commutative $k$-algebra and $X$ an $R$-scheme. Then an {\em $R$-central coherent $\ox$-bimodule} is a coherent sheaf $\calb \in \mo X \times_R X$ such that the two projection maps $p_1,p_2:\supp \calb \lm X$ are finite. Its left and right $\ox$-module structures are $\ _{\ox}\!\calb = p_{1*} \calb, \calb_{\ox} = p_{2*}\calb$. We say $\calb$ is {\em locally free of rank $r$} if $\ _{\ox}\!\calb, \calb_{\ox}$ are. An {\em $R$-central quasi-coherent $\ox$-bimodule} is a direct limit of $R$-central coherent $\ox$-bimodule. These form a monoidal category, as does the full subcategory of coherent bimodules. Basic facts such as these about $\ox$-bimodules can be found in [VdB01p, section~3.1]. We will usually consider $k$-central $\ox$-bimodules so our default setting will be $R=k$.  

Let $X$ be a commutative smooth projective variety of dimension $d-1$. Recall that the point of departure for defining a non-commutative $\PP^1$-bundle is that a commutative $\PP^1$-bundle has the form $\proj_X \mbox{Sym}(\cale)$ where $\cale$ is a rank two vector bundle. Non-commutatively, we start with a locally free rank two $\ox$-bimodule $\cale$, interesting examples of which include line bundles supported on a $(2,2)$-divisor $D$ in $X \times X$. Then one can form a $\Z$-indexed algebra analogous to the tensor algebra (over $X$) $T(\cale)= \oplus_{i\leq j} T_{ij}$ whose degree 0 and 1 parts are given by 
$$ T_{ii} = \ox, \quad \mbox{and}\quad \ldots ,T_{01} = \cale, T_{12} = \cale^*, T_{23} = \cale^{**}, \ldots, 
T_{i,i+1} = \cale^{*i}, \ldots$$
where $\cale^{*i}$ is the $i$-th iterated right adjoint of $\cale$ (and left adjoints are used when $i<0$). 

The symmetric algebra $\cala=\mbox{Sym}(\cale)$ can be defined as the quotient of $T(\cale)$ by the ideal $(\mathcal{Q})$ generated in degree two, by the line bundles $\im(\ox \lm  \cale^{*i} \otimes_X \cale^{*(i+1)}) < T_{i,i+2}$. These line bundles determine the ``commutation relation''. 

Now $\cala$ is a $\Z$-indexed algebra so one can form the category of graded (right) modules $\Gr \cala$ as usual. Torsion modules, that is, direct limits of right bounded modules form a Serre subcategory tors and we can form the quotient category $Y=\proj \cala:=\Gr \cala/\mbox{tors}$ as usual. This is a noetherian quasi-scheme called a {\em non-commutative $\PP^1$-bundle} or, when $\dim X = 1$, a {\em non-commutative ruled surface}. There are the usual adjoint functors $\Psi:\Gr \cala \lm \proj \cala, \Omega:\proj \cala \lm \Gr \cala$. There exists a natural flat morphism $f:Y \lm X$ defined as follows. Firstly, for $m \in \Z$ we let $e_m$ denote $1 \in \cala_{mm} = \ox$ and $\oy :=\Psi(\ox \otimes_X e_0\cala)$. There are adjoint functors 
$$f_m^*:=\Psi(- \otimes_X e_m \cala): \Mod X \lm \Mod Y, f_{m*}:=\Omega(-)_m: \Mod Y \lm \Mod X.$$ 
The morphism $f$ is given by the adjoint functors $f_0^*,f_{0*}$. The bimodules $\cala_{mn}$ are all locally free so the functors $f_m^*$ are exact and, in particular, $f$ is flat. For the rest of part~\ref{Pncruled}, we will preserve the above notation whenever $Y$ is a non-commutative $\PP^1$-bundle.

The following was more or less proved in [VdB01p].
\begin{prop}  \label{pP1snoeth} 
A non-commutative $\PP^1$-bundle $Y = \proj \cala$ on a smooth projective variety $X$ is strongly noetherian.
\end{prop}
\textbf{Proof sketch.} Let $R$ be a commutative noetherian $k$-algebra. We will denote base extension from $k$ to $R$ by the subscript $R$. As in [AZ01, proposition~B8.1v)], we see that $\proj \cala_R$ is naturally equivalent to $(\proj \cala)_R$ so it suffices to show that the $\Z$-indexed $\calo_{X_R}$-bimodule algebra $\cala_R$ is noetherian. 

We use Van den Bergh's proof that $\cala$ is noetherian [VdB01p, section~6.3]. He first notes that the point scheme of $\cala$ is a projective variety say $E$, and the universal point module induces an $\calo_E$-bimodule algebra $\calb = \oplus \calb_{ij}$ which is {\em strongly $\Z$-indexed} in the sense that the multiplication maps $\calb_{ij} \otimes_E \calb_{jl} \lm \calb_{il}$ are surjective for all $i,j,k \in \Z$. In fact, they are all isomorphisms of locally free rank one bimodules. Now $\calb$ can be considered an $\ox$-bimodule algebra in the following way. There exist morphisms $\mu_i:E \lm X$ for each $i \in \Z$ such that $\mu_*\calb := \oplus (\mu_i,\mu_j)_* \calb_{ij}$ is an $\ox$-bimodule algebra. We next describe Van den Bergh's construction of a surjective morphism of  $\ox$-bimodule algebras $\cala \lm \mathcal{D}$ whose kernel $I$ is an invertible ideal. Let $t:E' \hookrightarrow E$ be the inclusion of the components of maximal dimension. We may pull back to obtain another $\calo_{E'}$-bimodule algebra $t^*\calb:=\oplus (t,t)^* \calb_{ij}$, which we can, as before, consider as the $\ox$-bimodule algebra $(\mu t)_* t^*\calb$. Let $\mathcal{D}$ be the positive $\Z$-indexed $\ox$-bimodule algebra which equals $(\mu t)_* t^*\calb$ in positive degrees but is $\ox$ in degree zero. Then Van den Bergh's morphism $\cala \lm \mathcal{D}$ is the one induced from the natural map $\cala \lm \mu_* \calb_{\geq 0}$. 

Now $I_R$ is an invertible ideal in $\cala_R$ so it suffices to show that $\mathcal{D}_R$ is noetherian. This will follow if $t^*\calb_{\geq 0}$ is noetherian. However, this is clear since $t^*\calb$ is a strongly indexed $\calo_{E'_R}$-bimodule algebra, so its graded category is naturally equivalent to the category of quasi-coherent sheaves on $E'_R$. This completes the proof of the proposition.

\section{Mori contractions for non-commutative ruled surfaces} 

Let $f: Y \lm X$ be the natural fibration of a non-commutative ruled surface. In this section, we show that the fibres of $f$ are $K$-non-effective rational curves with self-intersection zero, and that $f$ is the non-commutative Mori contraction associated to the Hilbert system of a fibre. 

We need some cohomology results of Izuru Mori (see [Mori, lemma~4.4]) which hold more generally for non-commutative $\PP^1$-bundles. Note that $\{f_m^*\ox(-l)\}_{m,l}$ is a set of generators for $\Mod Y$. 
\begin{lemma}  \label{limori} 
Let $f:Y \lm X$ be a non-commutative $\PP^1$-bundle and $\call \in \mo X$ be locally free. The following holds in $D^b(\Mod X)$. 
$$ Rf_{m*}(f_n^*\call) = \begin{cases}
\call \otimes_X \cala_{n,m}         & \text{if }  n \leq m .\\
0                                   & \text{if }  n = m+1  .\\
\call \otimes_X \cala_{m,n-2}^*[-1] & \text{if } n \geq m+2.
\end{cases}$$
\end{lemma}

\begin{prop}  \label{pcohncruled} 
Let $f:Y \lm X$ be a non-commutative $\PP^1$-bundle. Then i) $R^2f_{m*} = 0$ ii) $R^1f_{m*}f_n^* = \Psi(- \otimes_X \cala_{m,n-2}^*)$ and iii) $f_{m*}f_n^* = \Psi(-\otimes_X \cala_{n,m})$. Let $M \in \mo Y$. Then iv) $R^if_{m*} M \in \mo X$ for all $i,m$ and v) $R^1f_{m*} M = 0$ for all $m \gg 0$. 
\end{prop}
\textbf{Proof.}
Part i) follows from [Na, corollary~4.10, theorem~4.11]. We now show simultaneously that $R^1f_{m*}f_n^* \call = \Psi(\call \otimes_X \cala_{m,n-2}^*)$ and $f_{m*}f_n^* \call= \Psi(\call\otimes_X \cala_{n,m})$ for any $\call \in \mo X$. The case $\call$ locally free is just lemma~\ref{limori}. The proof for general $\call$ follows by induction on the length of a locally free resolution of $\call$ and the fact that $ - \otimes_X \cala_{m,n-2}^*$ and $-\otimes_X \cala_{n,m}$ are exact. This proves parts ii) and iii). Finally, part iv) is [Nb, corollary~3.3] and part v) is [Nb, lemma~3.4]. Both can also be proved by writing $M$ as a quotient of a direct sum of modules of the form $f_m^*\ox(-l)$ and appealing to i),ii) and iii). 

\vspace{2mm}

We have the following version of the Leray spectral sequence.

\begin{lemma}  \label{lleray}  
Let $f:Y \lm X$ be a flat morphism of noetherian quasi-schemes. Then for $N \in \Mod X, M \in \Mod Y$ we have the following convergent spectral sequence 
\[  \Ext_X^p(N,R^qf_*M) \Longrightarrow   \Ext_Y^{p+q}(f^*N,M)   .                 \]
In particular, if global cohomology on $X,Y$ are computed using $\ox$ and $\oy:=f^*\ox$, we have the Leray-Serre spectral sequence
\[  H^p(X,R^qf_*M) \Longrightarrow  H^{p+q}(Y,M)   .   \]
\end{lemma}
\textbf{Proof.} This follows from the Grothendieck spectral sequence since flatness and adjunction imply that $f_*$ takes injectives to injectives. 
\vspace{2mm}

\begin{cor}  \label{cpd2}  
Let $Y$ be a non-commutative $\PP^1$-bundle on a $(d-1)$-dimensional variety $X$. Then 
for $N \in \mo X$, the functor $\Ext^{d+1}_Y(f_m^* N ,-)= 0$.
\end{cor}

\begin{prop}  \label{pruledproper} 
A non-commutative $\PP^1$-bundle $Y$ is Ext-finite. 
\end{prop}
\textbf{Proof.} By [AZ, corollary~C3.14], it suffices to show that $Y$, i) is strongly noetherian, ii) has a set of flat generators of finite cohomological dimension and iii) $\Ext^i(M,N)$ is finite dimensional for every $M,N \in \mo Y$. Now i) is just proposition~\ref{pP1snoeth}. Also, $\{f_m^* \ox(-l)\}$ is a set of flat generators which, by corollary~\ref{cpd2} have finite cohomological dimension. Finally, iii) is [Nb, corollary~3.6]. 

We can now compute the cohomology of fibres. Some of part ii) below was proved in [Mori, theorem~5.2].
\begin{prop}  \label{prulings} 
Let $f: Y \lm X$ be the fibration of a non-commutive ruled surface $Y$ and $p \in X$ be a closed point. Then 
\begin{enumerate}
\item We have $\Ext^i_Y(f^*-,f^*-) = \Ext^i_X(-,-)$ so in particular, $f^*$ is fully faithful.
\item $ h^0(f^*\calo_p) = 1, h^1(f^*\calo_p) = 0, f^*\calo_p.f^*\calo_p = 0, \Ext^2(f^*\calo_p,f^*\calo_p) = 0, \Ext^2(f^*\calo_p,\oy) = 0.$
\item Distinct fibres of $f$ are non-isomorphic.
\end{enumerate}
\end{prop}
\textbf{Proof.} Putting together proposition~\ref{pcohncruled} with the Leray spectral sequence gives part i). Part ii) follows from part i) and the fact that $\oy = f^* \ox$. Finally, for closed points $p,q \in X$, we know $f_*f^*\calo_p = \calo_p$ so $f^*\calo_p \simeq f^*\calo_q$ if and only if $p=q$.
\vspace{2mm}

The partial BK-Serre duality of the next section (see theorem~\ref{thalfSerreduality}) means $\Ext^2(f^*\calo_p,\oy) = 0$ ensures $H^0(f^*\calo_p \otimes \w_Y) = 0$. We will also show in proposition~\ref{pdimOK} that $f^*\calo_p$ is 1-critical so it is indeed a $K$-non-effective rational curve with self-intersection zero. 

We can now check that the the fibration $f:Y \lm X$ of a non-commutative ruled surface is indeed a Fourier-Mukai morphism.
\begin{thm}  \label{trecover} 
Let $f:Y \lm X$ be the fibration of a non-commutative ruled surface and $p \in X$ a closed point. Suppose that $Y$ has Halal Hilbert schemes.
\begin{enumerate}
\item The Hilbert system of the natural quotient $\oy = f^* \ox \lm f^*\calo_p$ is $(f^{\flat} \ox)/X$. 
\item The functor $f^*$ is given by the Fourier-Mukai transform $\pi_*(- \otimes_X e_0\cala)$ (see section~\ref{sfibr} for the definition of $\pi_*$).
\end{enumerate}
\end{thm}
\textbf{Proof.} Part ii) follows from part i), proposition~\ref{pFMfrommap} and the fact that $f^{\flat}\ox = e_0 \cala$. We thus restrict our attention to proving i). We are assuming that Hilbert schemes exist and are well-behaved. Let $H$ be the connected component containing the point $q\in H$ corresponding to $f^*\calo_p$. The flat family $(f^{\flat}\ox)/X$ gives a morphism $\tau:X \lm H$ which is non-constant by proposition~\ref{prulings}iii). Also, deformation theory as in the proof of corollary~\ref{cfamily} and the cohomology computation of proposition~\ref{prulings}ii) ensure that $H$ is smooth curve in a neighbourhood of $q$. In fact, since $X$ must surject onto the irreducible component of $H$ containing $q$, we see $H$ itself must be a smooth curve and proposition~\ref{prulings} shows that $\tau$ is a bijection on closed points. 

To show $\tau$ is an isomorphism, it suffices to show it separates tangent vectors. The image under $\tau$ of a non-zero tangent vector to $X$ at $p$ can be represented by the natural quotient
$$f^*\ox \otimes_k \calo_{2p} \lm f^* \calo_{2p}.$$
Now the tangent space to $X$ is 1-dimensional so it suffices to show that this corresponds to a non-zero tangent vector in $H$. It suffices to show that the natural short exact sequence
$$0 \lm f^*\calo_p \lm f^*\calo_{2p} \lm f^*\calo_p \lm 0 $$
is not split. However, we know from proposition~\ref{prulings}i) that $f^*$ is fully faithful so $\Hom(f^*\calo_p,f^*\calo_{2p}) = k$ and the sequence cannot split. This completes the proof of the theorem. 

\vspace{2mm}

In section~\ref{shalal}, we will remove the Halal Hilbert scheme assumption in the theorem. In fact, the rest of this part will be devoted primarily to proving theorem~\ref{tqrsOK}. This means that our non-commutative ruled surface is indeed a non-commutative smooth proper surface and the fibration $f:Y \lm X$ is indeed the non-commutative Mori contraction contracting a fibre.

\section{Partial BK-Serre duality}

In this subsection we define the Serre functor for non-commutative $\PP^1$-bundles and show that a form of Serre duality holds. Full BK-Serre duality will be proved in section~\ref{sserre}. In the special case where $Y = \mathbb{P}_X(V)$ is a commutative ruled surface, we have the adjunction formula $\w_Y = f^*\w_X \otimes_Y \calo_{Y/X}(-2) \otimes_Y \det V$. Our Serre functor is reminiscent of this formula but it is interesting to note that the first two tensor factors do not exist individually in the non-commutative setting, only their combination. In this section, $Y=\proj \cala$ will always denote a non-commutative $\PP^1$-bundle over a variety $X$ of dimension $d-1$. 

Recall firstly that for any locally free $\ox$-bimodule $\mathcal{B}$ that we have a natural isomorphism  $\mathcal{B}^{**} \simeq \w_X^{-1} \otimes_X \mathcal{B} \otimes_X \w_X$ (see [VdB01p, lemma~3.1.8]). Thus we have 
\begin{prop} \label{pAshift2} 
For our non-commutative symmetric algebra $\cala$ over $X$, we have 
$$\w_X^{-1} \otimes_X \cala_{m-2,n-2} \otimes_X \w_X \simeq \cala_{m,n}.$$
\end{prop}
\textbf{Proof.}
Taking double right adjoints shifts all the generators $\cala_{m,m+1}$ of $\cala$ by two as well as the relations between them.
\vspace{2mm}

Hence, for any $\cala$-module $M$, there is a natural map 
$$ (M_{m-2} \otimes_X \w_X) \otimes_X \cala_{m,n} \lm M_{n-2} \otimes_X \w_X .$$
One stark difference between the $\Z$-indexed setting as opposed to the $\Z$-graded setting is that the shifts $M(n)$ are no longer $\cala$-modules. However, from what we have seen there is a natural $\cala$-module structure on $M(-2)\otimes_X \w_X$ and we define the Serre functor to be 
$$M \otimes_Y \w_Y := M(-2)\otimes_X \w_X.$$
It is clearly a functor which is an auto-equivalence on $\mo Y$.  
Comparing with the commutative adjunction formula, we see they are identical save for the missing $\det \cale$. One way to explain this is that in the $\Z$-indexed setting, one does not form the symmetric algebra of $\cale$ with tensor powers of $\cale$ with itself, but rather with its adjoints, and the determinant of $\cale \otimes_X \cale^*$ is trivial. 

\begin{prop} \label{pserreonfm} 
For $\call \in Mod X, m \in \Z$ we have $f_m^*\call \otimes \w_Y = f_{m+2}^*(\call \otimes_X \w_X)$.
\end{prop}
\textbf{Proof.} We just use proposition~\ref{pAshift2} to see 
$$ \call \otimes_X e_m \cala(-2) \otimes_X \w_X = \call \otimes_X \w_X \otimes_X e_{m+2} \cala $$
and then apply the quotient functor $\Psi: \Gr \cala \lm \proj \cala$ to obtain the proposition.
\vspace{2mm}

Below and for the rest of this section, we will omit the functor $\Psi$ when it is clear. Restricting the next result to $m=0, \call=\ox$ recovers [Mori, theorem~4.6] as a special case. 
\begin{thm} \label{thalfSerreduality} 
Let $\call \in \mo X$ be locally free and $M \in \mo \cala$. We have the following Serre duality  isomorphism 
$$ \Ext^i_Y(M,f^*_m\call \otimes \w_Y) \simeq \Ext^{d-i}_Y(f_m^*\call,M)^*$$
which is natural in $M$. In particular, the injective dimension of $f_m^*\call$ is $d$. The isomorphism is also natural with respect to morphisms of the form $f_m^* \call \lm f_{m'}^* \call'$. 
\end{thm} 
\textbf{Proof.} We mimic the standard proof of Serre duality in $\PP^n$ as, for example found in [Hart, chapter~III, theorem~7.1]. Incidentally, it is easy to prove the result for $M$ of the form $f^*_n \caln$ directly by using lemma~\ref{limori}  and the Leray spectral sequence. 

As usual, we start by establishing the perfect pairing in the case $i=0$. Below we construct a natural trace map $\tr: \Ext^d_Y(f^*_m\call,f^*_m\call \otimes \w_Y) \lm k$. Thus, to a morphism $\xi: M\lm f^*_m\call \otimes \w_Y$ we may associate the functional 
$$\Ext^d_Y(f_m^*\call,M) \xrightarrow{\Ext^d_Y(f_m^*\call,\xi)} \Ext^d_Y(f^*_m\call,f^*_m\call \otimes \w_Y) \xrightarrow{\tr} k$$
which gives a pairing in the $i=0$ case. 

We check the pairing is perfect first in the case where $M$ has the form $f_n^*\caln$ for $\caln \in \mo X$ locally free. Consider an element $\xi$ of 
$$ \Hom_Y(f^*_n \caln, f_m^*\call \otimes \w_Y) = \Hom_Y(f^*_n \caln, f_{m+2}^*(\call \otimes_X \w_X)) = \Hom_X(\caln, f_{n*}f_{m+2}^*(\call \otimes_X \w_X)) $$
which by propositions~\ref{pcohncruled} and \ref{pAshift2} is 
$$\Hom_X(\caln,\call \otimes_X \w_X \otimes_X \cala_{m+2,n}) = \Hom_X(\caln  \otimes_X \ ^*\!\cala_{m+2,n},\call \otimes_X \w_X) = \Hom_X(\caln  \otimes_X \cala^*_{m,n-2},\call \otimes_X \w_X).
$$
Proposition~\ref{pcohncruled} also gives 
$$
\begin{CD}
\Ext^d_Y(f_m^*\call,f^*_n \caln)  @=  \Ext^{d-1}_X(\call, \caln \otimes_X \cala^*_{m,n-2}) \\
@VVV  @VVV \\
\Ext^d_Y(f^*_m\call,f^*_{m+2}(\call \otimes_X \w_X))  @= \Ext^{d-1}_X(\call,\call\otimes_X \w_X)
\end{CD}
$$
But we just saw that $\xi$ is given naturally by a morphism $\caln  \otimes_X \cala^*_{m,n-2} \lm \call \otimes_X \w_X$ so we may use Serre duality on $X$ to get a trace map 
$$ \Ext^d_Y(f^*_m\call,f^*_{m+2}(\call \otimes_X \w_X))  \lm \Ext^{d-1}_X(\call,\call\otimes_X \w_X) \lm k$$ 
and the perfect pairing between $\Hom_X(\caln  \otimes_X \cala^*_{m,n-2},\call \otimes_X \w_X)$ and $\Ext^{d-1}_X(\call, \caln \otimes_X \cala^*_{m,n-2})$ induce the desired perfect pairing in the case $i=0$ and $M = f^*_n \caln$. 

First note that any $M \in \mo Y$ is a quotient of a finite direct sum of modules of the form $f^*_n \ox(-l)$ where we can assume $n,l$ are arbitrarily large. More precisely, given any $n_0 \in \Z$, one can find an $n> n_0$ such that there exists a surjection of the form $f^*_n \ox(-l)^r\lm M$ for any $l \gg 0$. Also the sequences of functors $\Ext^i_Y(-,f^*_{m+2}(\call \otimes_X \w_X))$ and $ \Ext^{d-i}_Y(f_m^*\call,-)^*$  are both contravariant $\delta$-functors and, when $i=0$, both are left exact by corollary~\ref{cpd2}. Thus Serre duality in the $i=0$ case is established. To show Serre duality for general $i$, it suffices to show both sides are co-effaceable. In fact, using lemma~\ref{limori} and the Leray spectral sequence we see for $i>0, n,l=l(n) \gg 0$ that 
$$\Ext^i_Y(f_n^*\ox(-l),f^*_{m+2}(\call \otimes_X \w_X))=0=\Ext^{d-i}_Y(f_m^*\call,f_n^*\ox(-l)).$$
This establishes the Serre duality isomorphisms in the theorem. Furthermore, corollary~\ref{cpd2} ensures $\id_Y f_m^*\call = d$ since $\Ext^i_Y(-,f_m^*\call) \simeq \Ext^{d-i}_Y(f_{m-2}^*(\call \otimes_X \w_X^{-1}),-)^*$. 

We now prove the final statement for which we need only establish naturality in the case $i=0$. Consider a morphism $\xi \in \Hom_Y(f_m^* \call, f_{m'}^*\call') = \Hom_X(\call, \call' \otimes_X \cala_{m'm})$. It suffices to construct a trace map $t:\Ext^d_Y(f_{m'}^*\call',f_m^*\call \otimes \w_Y) \lm k$ which makes the diagram below commute.
$$\begin{CD}
\Hom_Y(M,f^*_m\call \otimes \w_Y) \times \Ext^d_Y(f_m^*\call,M) @>>> 
\Ext^d_Y(f^*_m\call,f^*_m\call \otimes \w_Y) @>{\tr}>> k \\
  @AAA @AAA @AA{\id}A \\
\Hom_Y(M,f^*_m\call \otimes \w_Y) \times \Ext^d_Y(f_{m'}^*\call',M) @>>> 
\Ext^d_Y(f^*_{m'}\call',f^*_m\call \otimes \w_Y) @>{t}>> k \\
  @VVV @VVV @VV{\id}V \\
\Hom_Y(M,f^*_{m'}\call' \otimes \w_Y) \times \Ext^d_Y(f_{m'}^*\call',M) @>>> 
\Ext^d_Y(f^*_{m'}\call',f^*_{m'}\call' \otimes \w_Y) @>{\tr}>> k 
\end{CD}$$
We can use proposition~\ref{pcohncruled} to re-write the right hand squares as
$$\begin{CD}
\Ext^{d-1}_X(\call, \call \otimes_X \w_X) @>>> k \\
@AAA @AA{\id}A \\
\Ext^{d-1}_X(\call', \call \otimes_X \ ^*\cala_{m'm} \otimes_X \w_X) = 
\Ext^{d-1}_X(\call', \call \otimes_X \w_X \otimes_X \cala_{m'm}^*)@>t>> k \\
@VVV @VV{\id}V \\
\Ext^{d-1}_X(\call', \call' \otimes_X \w_X) @>>> k \\
\end{CD}$$
Serre duality on $X$ shows the middle term on the left to be 
$$\Ext^{d-1}_X(\call' \otimes_X \cala_{m'm}, \call \otimes_X \w_X) \simeq \Hom_X(\call, \call' \otimes_X \cala_{m'm})^*$$
and $t$ may be defined to be evaluation at $\xi \in \Hom_X(\call, \call' \otimes_X \cala_{m'm})$. Commutativity of the diagrams now follows from naturality of the Serre duality isomorphisms on $X$. This completes the proof of the theorem. 

\vspace{2mm}

\section{Internal Tor}

In this section, we consider the non-commutative symmetric algebra $\cala = \mbox{Sym}(\cale)$ as defined in the beginning of part~\ref{Pncruled}. 
Nyman defined internal Hom, $\otimes$ and Ext functors for $\cala$-modules. We show that there is similarly an internal Tor functor in special cases. Our setting is slightly different. We consider first a functor $-\intotimes - : \Mod \cala \times \Bimod (\cala,\ox) \lm \Mod X$ where $\Bimod (\cala,\ox)$ is the category of $(\cala,\ox)$-bimodules $B = \oplus_{m \in \Z} B_m$, that is, the $B_m$ are $\ox$-bimodules and there are multiplication maps $\cala_{m,n} \otimes_X \calb_n \lm \calb_m$ satisfying the usual associativity and unit axiom. For example $\cala e_m \in \Bimod (\cala,\ox)$ whilst $\ox \otimes_X e_m \cala \in \Mod \cala$. 

Let $M \in \Mod \cala, B \in \Bimod (\cala, \ox)$ and, abusing notation, $\mu$ denote scalar multiplication in either of these. Then one may define the internal tensor product as in [Na] by 
$$ M \intotimes B := \coker \bigl(\bigoplus_{l,m} M_l \otimes_X \cala_{l,m} \otimes_X B_m \xrightarrow{\mu \otimes 1 - 1 \otimes \mu} \bigoplus_n M_n \otimes_X B_n \bigr).$$
It is right exact being defined by cokernels. We may define $M$ or $B$ to be {\em internally flat} if $M \intotimes -$ or $- \intotimes B$ is exact respectively. Also, given $\call \in \Mod X, \calb \in \Bimod \ox$ we have as in [Na, proposition~3.5]
$$ (\call \otimes_X e_m \cala) \intotimes B = \call \otimes_X B_m, \ \ M \intotimes (\cala e_m \otimes_X \calb) = M_m \otimes_X \calb $$
so in particular, if $\call,\calb$ are locally free then $\call \otimes_X e_m \cala$ and $\cala e_m \otimes_X \calb$ are internally flat. 

The problem with defining internal Tor is that, although there are enough internally flat $\cala$-modules, we do not know if the same is true of internally flat $(\cala,\ox)$-bimodules. We thus restrict our attention to bimodules of {\em finite internal flat dimension}, that is, those which have a finite resolution by internally flat bimodules. As in [AZ01, section~C.2], our first step will be to define a set of $\intotimes$-acyclic $\cala$-modules which includes the set of internally flat modules. We will say $M \in \Mod \cala$ is $\intotimes$-acyclic if for any exact sequence of $(\cala,\ox)$-bimodules 
$$  0 \lm K \lm F \lm B \lm 0 $$
with $F$ internally  flat and $K$ of finite internal flat dimension, then the induced sequence 
$$ 0 \lm M \intotimes K \lm M \intotimes F \lm M \intotimes B \lm 0 $$
is also exact. 

\begin{prop} \label{ptensacy} 
Consider an exact sequence
$$ M^{\bullet}:  0 \lm M'' \lm M' \lm M \lm 0 $$
of $\cala$-modules where $M$ is $\intotimes$-acyclic.
\begin{enumerate} 
\item If $B$ is an $(\cala,\ox)$-bimodule of finite internal flat dimension then the sequence $M^{\bullet} \intotimes B$ is also exact.
\item The module $M'$ is $\intotimes$-acyclic iff $M''$ is.
\item Direct sums of $\intotimes$-acyclic modules are $\intotimes$-acyclic.
\item Any $\cala$-module has a resolution by $\intotimes$-acyclic modules.
\end{enumerate}
\end{prop}
\textbf{Proof.} Part i) follows on considering an exact sequence
$$ 0 \lm K \lm F \lm B \lm 0 $$
with $F$ internally flat and chasing the diagram one obtains on internal tensoring this with $M^{\bullet}$. Part ii) then follows from part i) and a similar diagram chase. Part iii) follows from definition whilst part iv) holds since for any $\call \in \Mod X$ locally free, we know $\call \otimes_X e_m \cala$ is internally flat and hence $\intotimes$-acyclic. 
\vspace{2mm}

The proposition shows that for any bimodule $B$ of finite internal flat dimension, we can define the derived functors $\inttor_i(-,B)$ by using $\intotimes$-acyclic or internally flat resolutions. Looking at double complexes, one sees immediately that these internal tor functors can also be computed using internally flat resolutions of $B$. The bimodule of interest for us is $\cala_0 := \oplus_m \cala_{mm}$. It is an $\cala$-bimodule and so, naturally, also an $(\cala,\ox)$-bimodule.

\begin{prop}  \label{pA0finflat} 
There is an exact sequence of $(\cala,\ox)$-bimodules
$$0 \lm \cala e_{n-2} \lm \cala e_{n-1} \otimes_X \cale^{*(n-1)} \lm \cala e_n \lm \ox \lm 0 .$$
Hence, the $\cala$-bimodule $\cala_0$ has finite internal flat dimension and $\rk \cala_{mn} = n-m+1$.
\end{prop}
\textbf{Proof.} This is essentially [VdB01p, theorem~6.1.2]. 

\section{Smoothness}  \label{sruledsmooth} 

In this section, we finally prove that a non-commutative $\PP^1$-bundle $Y \lm X$ is smooth of dimension $\dim X + 1$. The method is to construct finite resolutions by modules of the form $\oplus_m f_m^*\call_m$ which, by  theorem~\ref{thalfSerreduality} have injective dimension $\dim X + 1$. One could also use $\Z$-indexed versions of [MS06, lemma~3.4]. As usual, in this section, $\cala$ will be the non-commutative symmetric algebra used to construct $Y$. 

\begin{defn}  \label{dXtorsion} 
Let $M = \oplus_m M_m$ be an $\cala$-module. 
\begin{enumerate}
\item We say $M$ is {\em $X$-torsion} or {\em $X$-torsion-free} if all the sheaves $M_m$ are torsion or torsion-free respectively. 
\item We say $M$ is {\em $X$-induced} if it is the direct sum of modules of the form $\call \otimes_X e_m \cala$ for some $m \in \Z, \call \in Mod X$.
\end{enumerate}
\end{defn}
Similarly, we say that a $Y$-module is {\em $X$-torsion(-free)} or {\em $X$-induced} if it can be represented by such an $\cala$-module. 

We mimic the standard proof of graded Hilbert syzygies theorem to show that noetherian $\cala$-modules have finite resolutions by noetherian $X$-induced $\cala$-modules. The main difference with the classical theory is that $\cala_0$ is not semisimple so one cannot uniformly bound the length of the resolution. First note the following standard version of the Nakayama lemma.

\begin{lemma}  \label{lNAK}  
Let $M,F \in \mo \cala$. Then $M \intotimes \cala_0 = 0 \implies M=0$. Hence, if $F \lm M$ is a morphism with $F \intotimes \cala_0 \lm M \intotimes \cala_0$ surjective, then $F \lm M$ is also surjective.
\end{lemma}

Given an $\cala$-module $M$, we note that $M \intotimes \cala_0$ is coherent. We define the {\em non-zero degrees of $M\intotimes \cala_0$} to be the finite set of integers $\{m| (M \intotimes \cala_0)_m \neq 0\}$.

\begin{thm}  \label{tsyzygy} 
Let $M= M_d \oplus M_{d+1} \ldots \in \Gr \cala$ be a noetherian module with $M_d \neq 0$. 
\begin{enumerate}
\item Then there exists an $X$-induced $\cala$-module $F = \oplus_m \calf_m \otimes_X e_m \cala$ and a surjective map $F \lm M$ such that a) $\calf_m \leq M_m$, b) the non-zero degrees of $F \intotimes \cala_0$ and $M \intotimes \cala_0$ are the same and c) $F \intotimes \cala_0 \lm M \intotimes \cala_0$ is an isomorphism in degree $d$. We call any such $F$ a {\em tight cover}.
\item $M$ has a finite resolution by noetherian $X$-induced $\cala$-modules. If $M$ is $X$-torsion then we can even assume the resolution is by $X$-torsion $X$-induced modules.
\end{enumerate}
\end{thm}
\textbf{Proof.} We prove part i) first. For any $m \in \Z$ we consider the surjective map $M_m \lm (M \intotimes \cala_0)_m$. We may thus find a subsheaf $\calf_m$ of $M_m$ such that the induced map $\calf_m \lm (M \intotimes \cala_0)_m$ is surjective. Note that $(M \intotimes \cala_0)_d = M_d =\calf_d$ and that $F \intotimes \cala_0 \simeq \oplus \calf_m$. There is a natural map $F \lm M$ and, by construction $F \intotimes \cala_0 \lm M \intotimes \cala_0$ is surjective so $F \lm M$ is surjective too by the Nakayama lemma~\ref{lNAK}. This proves part i)

We now prove the first statement of part ii) and start by showing that $X$-induced $\cala$-modules are $\intotimes \cala_0$-acyclic. Let $\mathcal{F} \in \Mod X$ and $\call^{\bullet} \lm \mathcal{F}$ be a finite locally free resolution. We obtain an internally flat resolution 
$$\call^{\bullet} \otimes_X e_m \cala \lm \mathcal{F} \otimes_X e_m \cala $$ 
from which we can compute $\inttor_i(\mathcal{F} \otimes_X e_m \cala, \cala_0)=0$  for $i>0$.
Thus any $X$-induced $\cala$-module is $\intotimes \cala_0$-acyclic. Also, $\cala_0$ has an internally  flat resolution of length 2 by proposition~\ref{pA0finflat}, so replacing $M$ with an appropriate syzygy with respect to a resolution by $X$-induced modules, we may assume $M$ is $\intotimes \cala_0$-acyclic. We now argue by induction on the number of non-zero degrees of $M \intotimes \cala_0$. Consider an exact sequence
$$ 0 \lm N \lm F \lm M \lm 0 $$
where $F$ is a tight cover as in part i). Note that $N$ is also $\intotimes \cala_0$-acyclic and that we have an exact sequence
$$ 0 \lm N \intotimes \cala_0 \lm F \intotimes \cala_0 \lm M \intotimes \cala_0 \lm 0 .$$
Part i) ensures that $N \intotimes \cala_0$ has strictly fewer non-zero degrees than that of $M \intotimes \cala_0$ so the resolution exists by induction. 

We now prove the final statement of part ii). Going through the proof in the previous paragraph, we see that a resolution can be constructed by repeatedly taking tight covers. Since the tight cover of an $X$-torsion $\cala$-module is $X$-torsion, we are done.
\vspace{2mm}

This theorem together with theorem~\ref{thalfSerreduality} immediately give

\begin{thm}  \label{tncruledsmooth} 
A non-commutative $\PP^1$-bundle $Y \lm X$  is smooth of dimension $\dim X + 1$.
\end{thm}

\section{BK-Serre duality}  \label{sserre} 

Now that we know non-commutative $\PP^1$-bundles are smooth, we can now prove BK-Serre duality in full. The proof is virtually the same as that of theorem~\ref{thalfSerreduality} so we will concentrate on the required modifications to that proof. Throughout this section, $Y$ denotes a non-commutative $\PP^1$-bundle on a $(d-1)$-dimensional smooth projective variety $X$. 

For technical purposes only, we introduce the following 
\begin{defn}  \label{dsuffneg} 
A  noetherian $X$-induced $Y$-module is said to be {\em sufficiently negative} if it is the direct sum of modules of the form $f^*_m\ox(-l)$ for $m,l=l(m) \gg 0$.
\end{defn}

\begin{lemma} \label{lcohvan} 
For $M \in \mo Y$ we have $\Ext^i(F,M) = 0=\Ext^{d-i}(M,F)$ for $i>0$ and $F$ a  sufficiently negative $X$-induced $Y$-module.
\end{lemma}
\textbf{Proof.} When $M$ has the form $f_n^* \call$ with $\call \in \mo X$ locally free then, as we have already seen, the lemma follows from lemma~\ref{limori}. The general case follows by taking a resolution of $M$ by sufficiently negative $X$-induced $Y$-modules and the fact that $Y$ is smooth.
\vspace{2mm}

\begin{cor}  \label{ccoefface} 
Let $M \in \mo Y$. The two sequences of contravariant $\delta$-functors 
$$\Ext^i(-,M), \Ext^{d-i}(M,- \otimes \w_Y)^*: \mo Y \lm \mo k$$
are co-effaceable.
\end{cor}

\begin{thm}   \label{tSerreduality} 
For a non-commutative $\PP^1$-bundle $Y$, the functor $- \otimes \w_Y$ is a Serre functor. Hence $Y$ is Gorenstein. 
\end{thm}
\textbf{Proof.} For $M,N \in \mo Y$, we need to show there is a ``Serre duality isomorphism'' $\Ext^i(N,M) \simeq \Ext^{d-i}(M,N\otimes \w_Y)^*$ which is natural in $M,N$. Theorem~\ref{thalfSerreduality} gives the case at least when $N$ is a sufficiently negative $X$-induced module. We start by fixing $M$ and exact sequences 
$$ 0 \lm M' \lm F_0 \lm M \lm 0, \ \ F_1 \lm M' \lm 0 $$
where $F_1,F_0$ are sufficiently negative $X$-induced $Y$-modules. By the corollary~\ref{ccoefface} and [BK, proposition~3.4], it suffices to show that there is a natural isomorphism $\Hom(-,M) \simeq \Ext^d(M,- \otimes \w_Y)^*$. Consider the following diagram with exact rows
$$\diagram
0 \rto & \Hom(M,M) \rto \dto^{\phi} & \Hom(F_0,M) \rto \dto & \Hom(F_1,M)\dto  \\
0 \rto & \Ext^d(M,M \otimes \w_Y)^* \rto & \Ext^d(M,F_0 \otimes \w_Y)^* \rto & \Ext^d(M,F_1 \otimes \w_Y)^*
\enddiagram$$
where the right hand square commutes by the natural Serre duality isomorphisms in theorem~\ref{thalfSerreduality} and $\phi$ is chosen to be the isomorphism which makes the left hand square commute. As in the proof of theorem~\ref{thalfSerreduality}, we need a trace map which will be given by $\phi(\id_M)$. In other words, the pairing between $\Hom(N,M)$ and $\Ext^d(M, N \otimes \w_Y)$ is given by associating to $\xi: N \lm M$, the image of $\id_M$ under
$$ \Hom(M,M) \lm \Ext^d(M,M \otimes \w_Y)^* \xrightarrow{\Ext^d(M, \xi \otimes\w_Y)^*} \Ext^d(M,N \otimes \w_Y)^*.$$
We show the pairing is perfect first for $N=F$ a sufficiently negative $X$-induced $Y$-module. It suffices to show that in this case, the pairing recovers the Serre duality isomorphism of theorem~\ref{thalfSerreduality}. By lemma~\ref{lcohvan}, we may assume $F$ is sufficiently negative that $\Ext^1(F,M') = 0$ so any morphism $\xi:F \lm M$ lifts to $\xi':F \lm F_0$. We now have the following commutative diagram 
$$\diagram
\Hom(M,M) \rto \dto^{\phi} & \Hom(F_0,M) \rto \dto & \Hom(F,M)\dto  \\
\Ext^d(M,M \otimes \w_Y)^* \rto & \Ext^d(M,F_0 \otimes \w_Y)^* \rto & \Ext^d(M,F \otimes \w_Y)^*
\enddiagram$$
which shows that our new pairing associates to $\xi$, its image under the old Serre duality isomorphism $\Hom(F,M) \lm \Ext^d(M,F \otimes \w_Y)^*$. To show that the pairing is perfect for general $N$ now follows by taking a resolution by sufficiently negative $X$-induced $Y$-modules as in theorem~\ref{thalfSerreduality}. This completes the proof of the theorem.

\section{Dimension}

In this section we define a dimension function for a non-commutative ruled surface $Y = \proj \cala$ where $\cala$ is the non-commutative symmetric algebra $\mbox{Sym} \cale$. The definition is a little roundabout.

Our starting point is to generalise the theory of Hilbert polynomials. Below, if $\call \in \mo X$ then we let $[\call]$ denote its image in the Grothendieck group $K_0(X)$. We start with the additive function 
$$\chi_n: K_0(Y) \lm K_0(X): M \mapsto [f_{n*} M] - [R^1f_{n*} M].$$
Suppose now that we have another additive function $\chi':K_0(X) \lm \Z$ such that for any $M \in \mo Y	$, the function $p_M:n \mapsto \chi'\chi_n(M)$ is a polynomial. Then we obtain an exact dimension function on $\mo Y$ by defining $\dim M = \deg p_M$ as usual. Alternatively, since $R^1f_{n*} M = 0$ for $n \gg 0$, we can also compute the dimension using the growth rate of $\chi'(f_{n*}M)$. In fact the following lemma shows we can even use the growth rate of $\chi'(M_n)$ where $\oplus M_n$ is a noetherian $\cala$-module representing $M$. 

\begin{lemma}  \label{lamp} 
Let $M_{\bullet}= \oplus M_n$ be a noetherian $\cala$-module representing $\Psi(M_{\bullet}) = M$. Then $f_{n*}M = M_n$ for $n \gg 0$.
\end{lemma}
\textbf{Proof.} Consider the $\cala_{>0}$-torsion functor $\tau := \varinjlim_i 
\shom_{\cala}(\cala/\cala_{\geq i}, -)$. It suffices to show that $R^i\tau M_{\bullet}, i\geq 0$ is right bounded, that is, zero in large degrees. [Nb, theorem~2.6] gives the lemma for $M_{\bullet}$ an $X$-induced $\cala$-module and the general case follows by writing $M_{\bullet}$ as a quotient of such a module. 

\vspace{2mm}

There are two natural condidates for $\chi'$, the rank and degree. Using the rank recovers the usual (Hilbert) dimension of the $\cala\otimes_X k(X)$-module $M \otimes_X k(X)$ so we shall denote it $\dim_{k(X)} M \otimes_X k(X)$. Recall the degree of a coherent sheaf $\call$ on $X$ is $\deg \call = \chi(\call) - (\rk{\call}) \chi(\ox)$. It is tempting to use $\chi' = \deg$ to define dimension for non-commutative ruled surfaces, but unfortunately, if $\cale$ is a line bundle on a $(2,2)$-divisor, then the degrees of $\ox \otimes_X \cala_{mn}$ grow cubically with $n-m$ and not quadratically, as occurs in the commutative case. However, this function is well-behaved for $X$-torsion modules.
\begin{prop}  \label{ptXdim} 
Let $M \in \mo Y$ be an $X$-torsion module and consider the function $p(n) = \deg \chi_n (M)$. 
\begin{enumerate}
\item $p(n)$ is a polynomial function whose degree, denoted $\dim_{\tau} M$, is $-\infty,0$ or 1.  
\item If $M = f_m^*\calo_p$ then $p(n)=n-m+1$. 
\end{enumerate}
\end{prop}
\textbf{Proof.} Part ii) follows from proposition~\ref{pcohncruled} and proposition~\ref{pA0finflat}. This in turn gives i) since we may resolve $M$ by $X$-torsion $X$-induced modules using theorem~\ref{tsyzygy}ii).
\vspace{2mm}

Let $\tau_X: \Mod X \lm \Mod X$ denote the $X$-torsion functor and note that the $X$-torsion submodule of an $\cala$-module is an $\cala$-submodule. We may now give the 
\begin{defn}  \label{ddim} 
For $M \in \mo Y$ we define its {\em dimension} to be 
$$\dim M = \max\{\dim_{\tau} \tau_X M, 1 + \dim_{k(X)} M \otimes_X k(X)\}.$$
\end{defn}
In other words, $\dim M$ is $\dim_{\tau} M$ if $M$ is $X$-torsion, and is $1 + \dim_{k(X)} M \otimes_X k(X)$ otherwise. 

\begin{prop}  \label{pdimexact} 
The dimension function $\dim$ is a compatible dimension function on $Y$.
\end{prop}
\textbf{Proof.} We first prove $\dim$ is exact. Let 
$$0 \lm M' \lm M \lm M'' \lm 0$$
be an exact sequence in $\mo Y$. Suppose first that $M \otimes_X k(X)\neq 0$ so $\dim M \geq 1$. Then at least one of $M' \otimes_X k(X), M'' \otimes_X k(X)$ is also non-zero. The fact that $\dim M = \max\{\dim M',\dim M''\}$ now follows from exactness of $\dim_{k(X)}$ and the fact that the dimension of $X$-torsion modules is at most 1. If $M \otimes_X k(X) = 0$ then both $M',M''$ are also $X$-torsion so we are done by exactness of $\dim_{\tau}$.

 We now show the Serre functor preserves dimension. Recall the Serre functor is given by the formula $M \otimes \w_Y = M(-2) \otimes_X \w_X$. If $M$ is not $X$-torsion, then $M \otimes_X k(X)$ is unaffected by $- \otimes_X \w_X$ and the shift by 2 does not change the degree of the Hilbert polynomial so $\dim$ is preserved in this case. Suppose now that $M$ is $X$-torsion and $M_{\bullet}$ is an $X$-torsion $\cala$-module representing $M$. Then $\dim M$ can be computed using the growth rate of the function $n \mapsto \mbox{length}\ M_n$. Again in this case, $M_n$ is unaffected by $- \otimes_X \w_X$ so $\dim$ is preserved as before.

\textbf{Remark:} The proof for exactness works with $\dim_{k(X)}$ replaced with any exact dimension function on $\mo Y\otimes_X k(X)$ and $\dim_{\tau}$ any exact dimension function on $X$-torsion modules.
\vspace{2mm}

\begin{prop}  \label{pdimOK} 
\begin{enumerate}
\item For any zero-dimensional module $P\neq 0$ we have a) $f_{n*}P \neq 0$ so in particular $h^0(P) \neq 0$ and b) $\Ext^1(P,\oy)=0$.
\item The fibre $f^*\calo_p$ is 1-critical.
\item The module $\oy$ is 2-critical. 
\end{enumerate}
In particular, a non-commutative ruled surface has classical cohomology.
\end{prop}
\textbf{Proof.} We first prove i)a).  For $n \gg 0, R^1f_{n*} P = 0$ whilst  $f_{n*} P \neq 0$. Thus the Hilbert polynomial $p(n) = h^0(f_{n*}P) - h^0(R^1f_{n*}P)$ must be a positive constant. It follows that for any $n$, $h^0(f_{n*}P) \geq p(n) > 0$ so $f_{n*}P \neq 0$. To prove i)b), it suffices by BK-Serre duality and proposition~\ref{pdimexact} to show that $\Ext^1(\oy,P) = 0$. Suppose this is not the case so there is an exact sequence of noetherian $\cala$-modules
$$ (*) \hspace{1cm} 0 \lm P_{\bullet} \lm M_{\bullet} \lm \ox \otimes_X e_0\cala \lm 0 $$
where $\Psi P_{\bullet}=P$. Looking in degree zero we obtain an exact sequence 
$$ 0 \lm  P_0 \lm M_0 \lm \ox \lm 0 $$
which splits and thus induces a morphism $\ox \otimes_X \cala \lm M$ which splits (*). This completes the proof of part i).

We now prove ii). Proposition~\ref{ptXdim}ii) shows that $\dim f^*\calo_p = 1$. We first show that $f^*\calo_p$ is 1-pure for suppose $P$ is a 0-dimensional submodule. We consider the injection $f_{n*}P \hookrightarrow f_{n*}f^*\calo_p$. Now proposition~\ref{pcohncruled} shows that $f_{n*}f^*\calo_p=0$ for $n< 0$ so we see $P=0$ by part i). Let $N$ be a non-zero submodule of $f^*\calo_p$. Its Hilbert polynomial must linear with leading co-efficient at least one, so the quotient $f^*\calo_p/N$ must have constant Hilbert polynomial. This proves part ii)

To prove part iii), note first that $\oy$ is $X$-torsion-free so any non-zero submodule $N$ must have dimension $\dim N = \dim_{k(X)} N \otimes_X k(X) \geq 1$. Furthermore, the argument in part ii) shows 
$\oy \otimes_X k(X)$ is 1-critical with respect to $\dim_{k(X)}$ so $\oy$ must be 2-critical.

\textbf{Remark:} We do not know if the dimension function is finitely partitive. The problem is that we do not know if there could be 1-dimensional $X$-torsion-free modules with arbitrarily long filtrations with 1-dimensional $X$-torsion composition factors. This cannot happen if the Hilbert polynomials with respect to $\chi' = \deg$ are quadratically bounded.

\section{Hilbert schemes for $\Gr \calb$}

In this section we work more generally with an arbitrary $\Z$-indexed $\ox$-bimodule algebra $\calb$. We wish to show that the Hilbert functors in $\Gr \calb$ are representable by projective schemes. The proof is completely analogous to the one given in [AZ01, sections~E.4,E.5] so we shall only briefly sketch their proof, giving details only when there is a need to show how the proof must be modified. They deal with the case where $\calb$ is a graded $k$-algebra. As we shall soon see, the extension from the graded to the indexed setting comes for free, but the extension to the $\ox$-bimodule algebra setting requires some thought.

As usual we shall assume that $X$ is a projective scheme, say with chosen very ample line bundle $\ox(1)$. This allows us to define, for any coherent sheaf $\calf \in \mo X$, its Hilbert polynomial $h(\calf;n)$. Now polynomials $h(n)$ can be ordered by their behaviour as $n \lm \infty$, or equivalently, by the lexicographic ordering on the co-efficients. Thus we can talk about semi-continuity of functions whose values are rational polynomials.

We need the following elementary results from commutative algebraic geometry.
\begin{prop}  \label{philbpoly} 
Let $S$ be a noetherian scheme. Then 
\begin{enumerate}
\item Given a coherent sheaf $\calf$ on $X \times S$, the set function $h:S \lm \mathbb{Q}[n]:s \mapsto h(\calf \otimes_S k(s);n)$ is upper semi-continuous.
\item Let $\phi:\calf \lm \calf'$ be a map of coherent sheaves on $X \times S$, both of which are flat over $S$. Then the locus in $S$ where $\phi$ is zero is scheme-theoretically closed.
\end{enumerate}
\end{prop}
\textbf{Proof.} i) Firstly, by generic flatness, there is an open set $U \subseteq S$ such that $\calf|_U$ is flat over $U$. The Hilbert polynomials of the closed fibres in $U$ are all the same, and we denote it $h_U$. We need to show for any $s \in S$, the corresponding Hilbert polynomial $h_s:= h(\calf \otimes_S k(s);n)$ satisfies $h_s \geq h_U$. Let $\pi: X \times S \lm S$ denote projection. By the stable version of ``cohomology commutes with base-change'', there is an $m_0\gg 0$ such that for any $m \geq m_0$ and $x \in U \cup \{s\}$ we have 
$$ \pi_*\calf(m) \otimes_S k(x) \simeq H^0(X,\calf(m) \otimes_S k(x)).$$
The result now follows by the classical upper semi-continuity result.

ii) By flatness, we may pick $m_0$ large enough so that for any $m \geq m_0$, the cohomology of $\calf(m),\calf'(m)$ commute with arbitrary base change $\rho:S' \lm S$, that is 
$\rho^*R^q\pi_* \calf (m) \simeq R^q\pi_* (\id_X \times \rho)^*\calf(m)$ and similarly for $\calf'(m)$. Let $\Gamma_*$ denote the functor $\oplus_{m \in \Z} \pi_*(- \otimes_X \ox(m))$. We may assume $m_0$ large enough that $\pi_* \calf(m_0), \pi_*\calf'(m_0)$ generate, modulo torsion, the $\Gamma_*\calo_{X \times S}$-modules $\Gamma_*\calf, \Gamma_* \calf'$. Finally, pick $m_0$ large enough so $\pi_*\calf(m_0),\pi_*\calf'(m_0)$ are locally free. Then the zero locus of $\phi$ is the same as the zero locus of $\pi_*(\phi \otimes_X \ox(m_0)): \pi_*\calf(m_0) \lm \pi_*\calf'(m_0)$ which we know is a closed subscheme of $S$.

\textbf{Remark:} Note that these are relative versions of [AZ01, lemma~E5.1i) and v)]. The relative versions of the other parts of their lemma are well-documented in the literature.
\vspace{2mm}

We assume henceforth in this section that $\calb$ satisfy the following conditions.
\begin{enumerate}
\item $\calb$ is {\em connected}, in the sense that all the $\calb_{ii}= \ox$ and $\calb_{ij} = 0$ if $i>j$.
\item $\calb$ is {\em locally finite}, in the sense that the $\calb_{ij}$ are all locally free bimodules of finite rank.
\item $\calb$ is {\em strongly noetherian} in the sense that $\Gr \calb$ is a strongly locally noetherian category.
\end{enumerate}
Given a noetherian graded $\calb$-module $M$, we can associate to it its ``double'' Hilbert function $h(M;j,n)$, which assigns to any $j \in \mathbb{N}$ the Hilbert polynomial $h(M_j;n)$ of the coherent sheaf $M_j$. Also, if $R$ is a commutative noetherian ring and $M$ a noetherian graded $\calb \otimes_k R$-module which is flat over $R$, we can talk about the double Hilbert function of $M$ since for every $j$, $M_j$ is flat over $R$ too so the Hilbert polynomials of its closed fibres are all the same.

To define the Hilbert functor, we first fix a double Hilbert function $h$ and $F \in \Gr \calb$. We consider the Hilbert functor $\Hilb(F,h)$ which assigns to any commutative noetherian ring $R$, the set of $R$-flat objects in $\Gr \calb$ whose double Hilbert function is $h$. As usual, we will drop the notation $F,h$ from $\Hilb(F,h)$ when they are understood. Just as in [AZ01, section~E.4], the functor can be extended to non-noetherian rings by declaring that it is limit preserving.

Let $I \subset \Z$ be a variable finite subset, and for any $M \in \Mod \calb$ we define $M_I:= \oplus_{i \in I} M_i$. Continuing as in [AZ01], we will define a subfunctor $\Hilb^I$ and quotient functor $\Hilb_I$ of $\Hilb$ as follows. First consider the following three conditions on a quotient map $q:F \otimes_k R \lm Q$ in $\Gr \calb_R$.
\begin{itemize}
\item[a)] The kernel $K:=\ker q$ is generated in degree $I$, that is, by $K_I$.
\item[b)] For each $i \in I$, $Q_i$ is $R$-flat with Hilbert polynomial $h(i,-)$.
\item[b')] $Q$ is $R$-flat with double Hilbert function $h$.
\end{itemize}
We define $\Hilb_I(R)$ to be the set of isomorphism classes of quotients $q:F \otimes_k R \lm Q$ satisfying a) and b), while $\Hilb^I(R)$ is the set of those satisfying a) and b'). 		

The indexed $\ox$-bimodule version of [AZ01, theorem~E4.3] is

\begin{thm}  \label{thilbgr} 
Let $\calb$ be a connected, locally finite, strongly noetherian $\ox$-bimodule algebra where $X$ is a projective scheme. Fix $F \in \Gr \calb$ and a double Hilbert function $h$. Then the Hilbert functor $\Hilb(F,h)$ is representable by a projective scheme, and in fact can be identified with the quotient functor $\Hilb(F,h)_I$ for some finite set $I$. 
\end{thm}
\textbf{Proof Sketch.} Most of the proof in [AZ01, section~E.4] carries over with the obvious word substitutions. Consequently, we will only outline their proof and expand when modifications are required. The reader with a copy of [AZ01] in hand should then have no problems determining the proof of the theorem. 

The finite subsets $I\subset \Z$ form a direct system by inclusion and induce an inverse system on the $\Hilb_I$. As usual [AZ01, lemma~E4.6ii)] we have $\Hilb = \varprojlim \Hilb_I$. We wish to show first that each $\Hilb_I$ is representable by a projective scheme (this is the analogue of [AZ01, lemma~E4.6i)]). Let $q:F \otimes_k R \lm Q$ represent an $R$-point of $\Hilb_I$. For each $i \in I$, the quotient map $F_i \otimes_k R \lm Q_i$ defines an $R$-point of the Hilbert scheme of quotients of $F_i$ with Hilbert polynomial $h(i,-)$. Hence $q$ determines, and is uniquely determined by, an $R$-point of some finite product of classical projective Hilbert schemes. We wish to show that the condition to be in $\Hilb_I$ is scheme-theoretically closed. If $K := \ker q$, then this condition is precisely that for every $i,j \in I$ we have that the map $K_i \otimes_X \calb_{ij} \lm Q_j$ is zero. This is scheme-theoretically closed by proposition~\ref{philbpoly}ii).

The next step is to show that the inverse system of $\Hilb_I$ stabilises. To do so, we need to show $\Hilb^I$ is represented by a locally closed subscheme of $\Hilb_I$ (this is [AZ01, proposition~E4.10]). As usual, we prove this by proving the analogue of [AZ01, proposition~E4.8], which in turn depends on our upper semi-continuity result proposition~\ref{philbpoly}i). The rest of the proof is purely formal and involves analysing the inverse system of constructible sets $\Hilb_I - \Hilb^I$. We leave it to the reader to verify that the rest of the proof in [AZ01, section~E.4] carries over.

\section{Halal Hilbert schemes}  \label{shalal} 
In this section, $X$ is a projective scheme and $\calb$ is a $\Z$-indexed $\ox$-bimodule algebra. 
As usual, we have the quotient category $Y = \proj \calb:= (\Gr \calb)/\mbox{tors}$ and there are adjoint functors $\Psi:\Gr \calb \lm \proj \calb, \Omega: \proj \calb \lm \Gr \calb$ and $f_m^* = \Psi(- \otimes_X e_m \calb), f_{m*} = \Omega(-)_m$. We show, \`a la [AZ01, section~E5], that under natural hypotheses, the Hilbert schemes of $\proj \calb$ are countable unions of projective schemes. This is then used to show that non-commutative $\PP^1$-bundles have Halal Hilbert schemes and finally, prove theorem~\ref{tqrsOK}. 

\begin{thm}  \label{thilbproj} 
Let $\calb$ be a connected, locally finite, strongly noetherian $\Z$-indexed $\ox$-bimodule algebra such that $\proj \calb$ is Ext-finite and adically complete. Then the Hilbert functor $\Hilb$ in $\proj \calb$ is representable by a separated scheme, locally of finite type which is a countable union of projective schemes. 
\end{thm}
\textbf{Comments on Proof.} The proof in [AZ01, theorem~E5.1] carries over painlessly. For future reference, we recall some parts of the proof. Firstly, representability of the Hilbert functor by a separated algebraic space, locally of finite type follows directly from [AZ01, theorem~E3.1] so the remainder of the proof involves showing that this algebraic space is a countable union of projective schemes. Let $F$ be a noetherian $\calb$-module and $h$ a double Hilbert function. Let $F_{>j},h_{>j}$ be the truncations of $F,h$, that is, in degrees $\leq j$, they are zero, and in other degrees they are the same as $F$ or $h$ respectively. Since Hilbert schemes in $\Gr\calb$ are projective by theorem~\ref{thilbgr}, the key step is to show $\Hilb(\Psi F)$ is a countable union of the functors $\Hilb(F_{>j},h_{>j})$. This follows from [AZ01, lemma~E5.3] which in our setting is

\begin{lemma}  \label{lflatinproj} 
Let $\calb$ be a connected, locally finite, strongly noetherian $\Z$-indexed $\ox$-bimodule algebra. Suppose $R$ is a commutative noetherian ring and $M$ is a noetherian graded $\calb_R$-module such that $\Psi M$ is flat over $R$. Then for $j\gg 0$, the truncation $M_{>j}$ is flat over $R$ too. 
\end{lemma}

We have the following immediate 
\begin{cor}  \label{ccompdim} 
The dimension function on a non-commutative ruled surface is continuous.
\end{cor}
\textbf{Proof.} This follows from the lemma and continuity of double Hilbert functions. 
\vspace{2mm} 

Finally, there is also a version of Grothendieck's existence theorem.

\begin{thm}  \label{tgrothexist} 
Let $\calb$ be a connected, locally finite, strongly noetherian $\Z$-indexed $\ox$-bimodule algebra. Suppose that $\proj \calb$ is Ext-finite, $\Omega: \proj \calb \lm \Gr \calb$ has finite cohomological dimension and that for any noetherian $M \in \proj \calb$ we have i) $R^if_{m*} M = 0$ for $m \gg 0,i>0$ and ii)  $ R^if_{m*}M \in \mo X$ for all $i,m$. Then $\proj \calb$ is adically complete. 
\end{thm}
\textbf{Proof.} We first need
\begin{lemma}
For any commutative noetherian $k$-algebra $R$ and $\calm \in \proj \calb_R$ we have i) $R^if_{m*} \calm = 0$ for $m \gg 0,i>0$ and ii)  $ R^if_{m*}\calm \in \mo X_R$ for all $i,m$.
\end{lemma}
\textbf{Proof lemma.} Since $\Omega$ has finite cohomological dimension, we may proceed by downward induction on $i$. Assumptions~i),ii) in the theorem, give i) and ii) for $R$-objects of the form $\calm = N \otimes_k R$. We may pick $N \in \proj \calb$ such that there is an exact sequence
$$ 0 \lm \mathcal{K} \lm N \otimes_k R \lm \calm  \lm 0 .$$
The lemma now follows from the long exact sequence in cohomology. 

\vspace{2mm}\noindent
\textbf{Proof} of theorem continued. 
The proof in [AZ01, theorem~D6.1] reduces the theorem to establishing the following statement: for any finitely generated graded commutative $k$-algebra $R$ and $\calm \in \proj \calb_R$, we have $\Ext^1_Y(f_m^*\ox(-l),\calm) = 0$ for $m \gg 0$ and $l = l(m) \gg 0$. First we use the lemma to pick $m$ large enough so $R^1f_{m*} \calm = 0$. For such an $m$ we have $\Ext^1_X(\ox(-l),f_{m*} \calm) = 0$ for $l \gg 0$. The desired statement and hence theorem, now follows from the Leray spectral sequence.
\vspace{2mm}

\begin{prop}  \label{phalal} 
A non-commutative $\PP^1$-bundle $Y = \proj \cala$ has Halal Hilbert schemes.
\end{prop}
\textbf{Proof.} We use the criterion of theorem~\ref{thilbproj}. All the hypotheses have been checked except adic completion which follows from theorem~\ref{tgrothexist} and proposition~\ref{pcohncruled}.

\vspace{2mm}
We can finally prove that a non-commutative ruled surface is indeed a non-commutative smooth proper surface.

\noindent
\textbf{Proof theorem~\ref{tqrsOK}}. Let $Y  = \proj \cala$ be a non-commutative ruled surface. We know that $Y$ is strongly noetherian (proposition~\ref{pP1snoeth}), Ext-finite (proposition~\ref{pruledproper}), has compatible dimension function (proposition~\ref{pdimexact}), classical cohomology (proposition~\ref{pdimOK}) and Halal Hilbert schemes (proposition~\ref{phalal}). We need only show there are no shrunken flat deformations, so suppose that $M,M'\in \proj \cala$ are members of a flat family parametrised by a connected scheme of finite type. Suppose furthermore that there is a morphism $\phi: M \lm M'$. We need to show that it is an isomorphism if it is either injective or surjective. Suppose that it is injective and $\phi$ is represented by an injective morphism of noetherian $\cala$-modules $M_{\bullet} \lm M_{\bullet}$. For large $i$, the Hilbert polynomials of $M_i,M'_i$ are the same so $M_i \lm M'_i$ must be an isomorphism. Hence $\phi$ is an isomorphism and the same argument yields the case where $\phi$ is assumed to be surjective.

\part{Properties of Mori contractions}

Morphisms of non-commutative schemes, unlike their commutative counterparts, have very little structure, so it is hard to extract information from them. For example, suppose that $Y,X$ are quasi-schemes equipped with distinguished objects $\oy,\ox$ respectively. Then given a morphism $f:Y \lm X$, one does not expect any relationship between $\oy$ and $f^* \ox$ in general, though in the cases of interest, one would hope they were isomorphic. For example, when $\oy \simeq f^*\ox$ and $f$ is flat, then the Leray spectral sequence can be used to link the cohomology on $Y$ with that on $X$. 

Let $Y$ be a strongly noetherian Hom-finite quasi-scheme and $\calm\in \mo Y_X$ be a base point free Hilbert system parametrised by a generically smooth curve $X$. Rather than looking at arbitrary morphisms, we will restrict our attention to Fourier-Mukai morphisms $f : Y \lm X$ of such a Hilbert system. Recall this is defined via $f^*\call = \pi_*(\calm \otimes_X \call)$ in the notation of section~\ref{sfibr}. We may apply $\pi_*$ to the surjection $ \oy \otimes_k \ox \lm \calm$ to obtain a natural map $\nu: \oy \lm f^* \ox$. We see immediately that the composite $\oy \lm f^*\ox \lm f^* \calo_p$ is the usual quotient map so is in particular surjective. As seen in example~\ref{eaffinemap}, the map $\nu$ need not be an isomorphism. The driving question in this part will be to find conditions for $\nu$ to be an isomorphism in the case where $f$ is a non-commutative Mori contraction.

Section~\ref{sdisjoint} concerns sufficient criteria for the map $\nu$ to be injective. To get anywhere, we will need to assume that the dimension function on $Y$ is finitely partitive (see section~\ref{sass}). We have seen that the fibres of $\calm$ are given in terms of $f$ by $f^*\calo_p$. In section~\ref{scohmori}, we will see how information about these fibres can be used to tell us about $R^if_*$. This will be used in section~\ref{snuonto} to give sufficient criteria for $\nu$ to be an isomorphism. 

\section{Disjoint Fibres}  \label{sdisjoint} 

In this section, we let $Y$ be a non-commutative smooth proper surface. Consider a base point free Hilbert system $\calm/X$ parametrised by a generically smooth projective curve $X$ and its associated Fourier-Mukai morphism $f:Y \lm X$. We wish to relate the concept of fibres of $f$ being generically disjoint with the condition that $\nu: \oy \lm f^*\ox$ is injective.

The next result explains why we wish to look at Hilbert systems as opposed to more general flat families of quotients of $\oy$. 
\begin{prop}  \label{p1dimhilbsys} 
Let $X$ be a projective curve and $f:Y \lm X$ a Fourier-Mukai morphism associated to some base point free Hilbert system. Suppose that $p \in X$ is a smooth point such that the corresponding fibre satisfies $h^0(f^*\calo_p)=1$. Then the natural map $T_p X \hookrightarrow \Ext^1(f^*\calo_p,f^*\calo_p)$ is injective and the image is spanned by the extension
$$0 \lm f^* \calo_p \lm f^* \calo_{2p} \lm f^* \calo_p \lm 0 .$$
\end{prop}
\textbf{Proof.} This follows from the usual deformation theory as in the proof of corollary~\ref{cfamily}. \vspace{2mm}

We note some consequences of the condition that $\nu: \oy \lm f^* \ox$ is an isomorphism.

\begin{prop}  \label{pifnuiso} 
Let $f:Y \lm X$ be the Fourier-Mukai morphism associated to a base point free Hilbert system $\calm/X$ parametrised by a projective curve $X$. If further $\nu: \oy \lm f^*\ox$ is an isomorphism, then the following hold.
\begin{enumerate}
\item For any closed subscheme $D \subseteq X$, the natural map $\oy \lm f^* \calo_D$ is surjective.
\item For any smooth closed point $p \in X$ with $h^0(\calm \otimes_X \calo_p) = 1$, the natural map 
$$ \Ext^1(\calm \otimes_X \calo_p, \oy) \lm \Ext^1(\calm \otimes_X \calo_p,\calm \otimes_X \calo_p)$$
is non-zero. In particular, $H^1(\calm \otimes_X \calo_p \otimes \w_Y) \neq 0$. 
\item For distinct $M_1,\ldots,M_r \in \calm$ and a point $P \in\mo Y$, the vector spaces $\Hom(M_i,P)$ considered as subspaces of $\Hom(\oy,P) = H^0(P)$ are linearly independent. In particular, if $h^0(P)=1$ then $P$ lies on at most one fibre. 
\end{enumerate}
\end{prop}
\textbf{Proof.} Part i) follows from the fact that $f^*\ox \lm f^* \calo_D$ is surjective. For part ii), observe that we have the following morphism of extensions in $\mo X$
$$\diagram
0 \rto & \ox \rto \dto & \ox(p) \rto \dto & \calo_p \rto \dto^{\id} &  0 \\
0 \rto & \calo_p \rto & \calo_{2p} \rto & \calo_p \rto & 0
\enddiagram$$
Applying $f^*$ to the whole diagram yields a similar morphism of extensions where the bottom one corresponds to a non-trivial element of $\Ext^1(\calm \otimes_X \calo_p,\calm \otimes_X \calo_p)$ by proposition~\ref{p1dimhilbsys}. This and BK-Serre duality proves part ii). 

For part iii), choose $p\in X$ so that $M_i=f^*\calo_p$. Then by adjunction, $\Hom(M_i,P) = \Hom(\calo_p,f_*P)$ so corresponds to the sections of $\Hom(\ox,f_*P) = H^0(P)$ supported (scheme-theoretically) at $p$. 

\vspace{2mm}

The key to verifying that $\nu: \oy \lm f^* \ox $ is an isomorphism is to show that the conclusions of proposition~\ref{pifnuiso}i) hold. If $P$ is a point with $h^0(P)=1$ then part iii) shows that $P$ can lie on at most one fibre if $\nu$ is an isomorphism. This corresponds to the fact that given a morphism of commutative schemes, the fibres are disjoint. However, when $h^0(P) >1$, part iii) shows the correct way to generalise this ``disjointness'' of fibres. 

We look at the question of showing $\nu$ is injective. The following result gives a sufficient condition.

\begin{prop}  \label{poyinject} 
Let $f:Y \lm X$ be the Fourier-Mukai morphism associated to a base point free Hilbert system $\calm/X$  parametrised by a projective curve $X$. 
\begin{enumerate}
\item Let $M_1, \ldots, M_n \in \calm$ be such that for any simple $Y$-module $P$ we have that the subspaces $\Hom(M_i,P)':= \im(\Hom(M_i,P) \hookrightarrow H^0(P))$ are linearly independent in $H^0(P)$. Then $\oy \lm \oplus_{i=1}^n M_i$ is surjective. 
\item Suppose the dimension function on $Y$ is finitely partitive. Let $M_1,M_2,\ldots \in \calm$ be 1-dimensional fibres such that $\oy \lm \oplus_{i=1}^n M_i$ is surjective for all $n$. Then $\oy \lm f^*\ox$ is injective.
\end{enumerate}
\end{prop}
\textbf{Proof.}
For i), we show by induction on $n$ that $\oy \lm M_1 \oplus \cdots \oplus M_n$ is surjective. Consider inductively, the exact sequence
$$ 0 \lm K \lm \oy \lm \oplus_{i=1}^{n-1} M_i\lm 0 .$$
We are done if the natural map $K \lm M_n$ is surjective so suppose its cokernel $C$ has a simple quotient $P$. We have thus a commutative diagram
$$\begin{CD}
\oy @>>> M_n \\
@VVV @VVV  \\
\oplus_{i=1}^{n-1} M_i @>>> P
\end{CD}$$
This gives a non-zero element of $H^0(P)$ which is both in $\Hom(M_n,P)'$ and $\sum_{i=1}^{n-1} \Hom(M_i,P)'$.  

To prove ii), consider the exact sequence 
$$0 \lm K_n \lm \oy \lm \oplus_{i=1}^n M_i \lm 0 $$
and let $K = \cap_n K_n$. It suffices to show that $K=0$ since $\ker (\oy \lm f^*\ox)$ is contained in every $K_n$. Suppose this is not the case so $\oy/K$ must be 1-dimensional. Then $\{K_n/K\}_{n\in \mathbb{N}}$ is a strictly decreasing sequence with 1-dimensional factors. This contradicts the fact that $\dim$ is finitely partitive so we are done.

\vspace{2mm}

Unfortunately, to be able to use this result, we need to impose more hypotheses and take our cue from proposition~\ref{pifnuiso}. We start with the case which is very close to the commutative one where any point $P$ on a fibre $M$ has $h^0(P)=1$ and $H^1(M \otimes \w_Y)\neq 0$. 

\begin{prop}  \label{ph0P1} 
Let $M$ be a 1-critical $K$-non-effective rational curve with $M^2=0, H^1(M \otimes \w_Y) \neq 0$ and a base point free Hilbert system  $\calm/X$. Let $P$ be a point of $M$ with $h^0(P) =1$. Then $P$ lies on no other 1-critical fibre of $\calm$.
\end{prop}
\textbf{Proof.} Suppose $P$ also lies on some distinct fibre $M'\in \calm$ which is 1-critical. Consider an exact sequence of the form
$$ 0 \lm N' \lm M' \lm P \lm 0 .$$
Applying $\Rhom(M,-)$ to this sequence, we see from proposition~\ref{pextmmp} and BK-Serre duality that 
$$\Ext^1(M,P) \simeq \Ext^2(M,N') \simeq \Hom(N',M \otimes \w_Y)^*$$
Since $\calm$ is base point free, this is non-zero by proposition~\ref{pbptfree}. Now $N'$ is also 1-critical so there is an embedding $N' \hookrightarrow M \otimes \w_Y$. Then 
$$0 = \chi(M') - h^0(P) = \chi(N') \leq \chi(M \otimes \w_Y) = -h^1(M \otimes \w_Y)$$
which contradicts our assumption that $H^1(M \otimes \w_Y) \neq 0$. 
\vspace{2mm}

If we do not assume that points $P$ satisfy $h^0(P) = 1$, then we have to impose the condition of surjectivity of $\Ext^1(M,\oy) \lm \Ext^1(M,M)$ as found in the conclusion of proposition~\ref{pifnuiso}ii). We study this condition more, and in particular show it is generic.  

\begin{lemma}  \label{lontoextmm} 
Let $\calm/X$ be the Hilbert system of a $K$-non-effective rational curve $M$ with $M^2=0$. Suppose that the natural morphism $\e:\Ext^1(M,\oy) \lm \Ext^1(M,M)$ is surjective. Then $\e':\Ext^1(M',\oy) \lm \Ext^1(M',M')$ is surjective for generic $M' \in \calm$.
\end{lemma}
\textbf{Proof.} Shrinking $X$ to an appropriate affine neighbourhood $\spec R$ of the point $p\in X$ corresponding to $M$, we may suppose that $X$ is smooth and $M'$ has the following cohomological properties of $M$. We can assume $H^0(M') = k$ by semicontinuity (and the fact that we must have $h^0(M') >0$) so also $\Hom(M',M') = k$. Also, we can assume $H^0(M' \otimes \w_Y) = 0$ so $\Ext^2(M',M')=0$. In other words, every fibre of $\calm$ is $K$-non-effective rational with self-intersection zero. Since $M'^2=0$ we also have $\Ext^1(M',M') = k$ so surjectivity of $\e'$ just means it is non-zero.

Consider the exact sequence of $R$-modules
$$\Ext^1_{Y_R}(\calm,\oy \otimes_k R) \lm \Ext^1_{Y_R}(\calm,\calm) \lm C \lm 0 .$$
We need 
\begin{claim}
Let $\caln = \oy \otimes_k R$ or $\calm$. Then for any closed point $q \in X$, there are natural isomorphisms 
$$\Ext^1_{Y_R}(\calm,\caln) \otimes_R \calo_q  \simeq \Ext^1(\calm\otimes_R \calo_q,\caln\otimes_R \calo_q)$$  
\end{claim}
Note that the lemma follows from the claim since surjectivity of $\e$ ensures that $C$ is torsion and thus, that $\e'$ is surjective for generic $M'$.

\noindent
\textbf{Proof claim.} We use the Tor-Ext spectral sequence of [AZ01, corollary~C3.9] 
$$ \Tor^R_{-i}(\Ext^j_{Y_R}(\calm,\caln),\calo_q) \implies \Ext^{i+j}_{Y_R}(\calm, \caln \otimes_R \calo_q)$$
This yields an embedding
$$ \Ext^2_{Y_R}(\calm,\caln) \otimes_R \calo_q \hookrightarrow \Ext^2_{Y_R}(\calm,\caln \otimes_R \calo_q) \simeq \Ext^2(\calm \otimes_R \calo_q,\caln \otimes_R \calo_q)$$
where the isomorphism comes from [AZ01, proposition~C3.4i), proposition~C2.6iii)]. Now the right hand term vanishes, since $0=\Ext^2(M',\oy) = \Ext^2(M',M')$ for all $M' \in \calm$. Since this is true for all $q$, we see $\Ext^2_{Y_R}(\calm,\caln)=0$. Thus the Tor-Ext spectral sequence also shows that 
$$\Ext^1_{Y_R}(\calm,\caln) \otimes_R \calo_q  \simeq \Ext^1_{Y_R}(\calm,\caln\otimes_R \calo_q).$$
Now [AZ01, proposition~C3.1v)] shows also 
$$ \Ext^1_{Y_R}(\calm,\caln \otimes_R \calo_q) \simeq \Ext^1(\calm \otimes_R \calo_q,\caln \otimes_R \calo_q) $$
so the claim, and hence lemma are proved. 
\vspace{2mm}

We will say a condition holds {\em quasi-generically} if it holds on some non-empty quasi-open subset. The condition we will be interested in will be for fibres to be critical. Recall from corollary~\ref{cbertini} that this holds quasi-generically if it holds somewhere and the dimension function is continuous. 

\begin{thm}  \label{tptfibre} 
Let $Y$ be a non-commutative smooth proper surface with a finitely partitive dimension function. 
Let $\calm/X$ be a base point free Hilbert system of a $K$-non-effective rational curve with self-intersection zero and let $f:Y \lm X$ be its non-commutative Mori contraction. Assume either that a) the generic fibre of $\calm$ is critical or b) the quasi-generic fibre is critical and the ground field $k$ is uncountable. Suppose further that one of the following two conditions hold.
\begin{enumerate}
 \item There exists a $K$-non-effective rational $M \in \calm$  with $M^2=0$ such that the natural map $\Ext^1(M,\oy) \lm \Ext^1(M,M)$ is non-zero.
 \item The generic fibre $M \in \calm$ is ``fat-free'' in the sense that $H^1(M \otimes \w_Y) \neq 0$ and any point $P$ on a fibre must have $h^0(P) = 1$. 
\end{enumerate}
Then the natural map $\nu:\oy \lm f^* \ox$ is injective.
\end{thm}
\textbf{Proof.} We use the criterion of proposition~\ref{poyinject}. Arguing by semicontinuity as in lemma~\ref{lontoextmm}, there exists a sequence $M_1,M_2,\ldots \in \calm$ of distinct fibres which are all 1-critical $K$-non-effective rational curves with self-intersection zero. We need to check linear independence of the subspaces $\Hom(M_i,P)':= \im(\Hom(M_i,P) \hookrightarrow H^0(P))$. This is clear if $P = M_j$ for some $j$, since then proposition~\ref{pextmmp}iv) ensures that $\Hom(M_i,P) = 0$ unless $i=j$. We may thus assume that $P$ is a point. If the hypotheses of ii) hold, then the ``disjointness'' of fibres result of proposition~\ref{ph0P1} guarantees linear independence in this case too and we are done. 

We prove the theorem now under the assumptions of i). By lemma~\ref{lontoextmm} we can assume, on deleting a finite number of $M_i$, that $\Ext^1(M_n,\oy) \lm \Ext^1(M_n,M_n)$ is surjective for all $n$. 
Suppose that $\Hom(M_1,P)',\ldots, \Hom(M_n,P)'$ are linearly dependent and that $n$ is minimal with respect to this condition. Proposition~\ref{poyinject}i) ensures that $\oy \lm \oplus_{i=1}^{n-1} M_i$ is surjective. Now $\Hom(M_n,P) \neq 0$ while base point freedom ensures (propositions~\ref{pgenpure},\ref{pbptfree}) that $\xi(M_n,P)\leq 0$. Thus $\Ext^1(P,M_n\otimes \w_Y) \simeq \Ext^1(M_n,P)^* \neq 0$ and we may thus choose a non-split extension
$$  E: 0 \lm M_n \otimes \w_Y \lm L \lm P \lm 0 .$$
where $L$ is 1-critical. 

Suppose $\phi_i:M_i \lm P$ are morphisms, not all zero, such that their images $\phi'_i:\oy \lm M_i \lm P$ in $H^0(P)$ satisfy $\phi_n'= \sum_{i=1}^{n-1} \phi_i' \neq 0$. The ext computation in proposition~\ref{pextmmp} and BK-Serre duality show that for $i < n$, $\phi_i$ lifts to $\psi_i: M_i \lm L$. Thus we see that $\phi_n'= \sum_{i=1}^{n-1} \phi_i' \in H^0(P)$ lifts to a unique $\psi_n' \in H^0(L)$ since $H^0(M_n \otimes \w_Y) = 0$. We wish also to lift $\phi_n$ to $\psi_n:M_n \lm L$. Consider the commutative diagram of natural maps of ext spaces below.
$$\diagram
\Ext^1(P,M_n\otimes \w_Y) \rto^{\phi_n^*} \drto &   \Ext^1(M_n,M_n\otimes \w_Y) \dto \\
 & \Ext^1(\oy,M_n\otimes \w_Y) 
\enddiagram$$
If $E \in \Ext^1(P,M_n\otimes \w_Y)$ is the extension above, then we need to show $\phi_n^*(E) = 0$ for which it suffices, by our assumption i),  to show its image is zero in $\Ext^1(\oy,M_n\otimes \w_Y)$. This holds since  $\phi_n'$ lifts to $\psi_n'$. Thus $\phi_n$ lifts to $\psi_n:M_n\lm L$ which is injective since $M_n$ and $L$ are 1-critical. Now $\oy \lm \oplus_{i=1}^{n-1} M_i$ is surjective so $\im \psi_i \leq \im \psi_n$ which contradicts the hom computations of proposition~\ref{pextmmp}. The theorem is finally proved. 

\vspace{2mm}

The proof above immediately gives
\begin{schol}  \label{schptfibre}
Let $f: Y \lm X$ be a non-commutative Mori contraction and suppose the dimension function on $Y$ is finitely partitive. Let $f^*\calo_{p_1}, \ldots , f^*\calo_{p_n}$ be distinct 1-critical $K$-non-effective rational curves such that the natural maps $\Ext^1(f^*\calo_{p_i},\oy) \lm \Ext^1(f^*\calo_{p_i},f^*\calo_{p_i})$ are all surjective. Then the natural map $\oy \lm \oplus_{i=1}^n f^*\calo_{p_i}$ is surjective.
\end{schol}
\vspace{2mm}

The next result allows us in one special case, to relax the surjectivity of $\Ext^1(M,\oy) \lm \Ext^1(M,M)$ condition in i) of the theorem to $\Ext^1(M,\oy) \simeq H^1(M \otimes \w_Y)^* \neq 0$ as found in ii). 
\begin{prop}  \label{ph1oyis0} 
Let $Y$ be a non-commutative smooth projective surface with $H^1(\oy) = 0$ and $M$ a $K$-non-effective curve with $M^2=0, H^0(M)=k, H^1(M\otimes \w_Y)\neq 0$. Then $\Ext^1(M,\oy) \lm \Ext^1(M,M)$ is an isomorphism so in particular, $H^1(M\otimes \w_Y)=k$.
\end{prop}
\textbf{Proof.}  Consider the exact sequence
$$0 \lm I \lm \oy \lm M \lm 0 .$$
Our cohomology assumptions ensure $H^0(I) = H^1(I)=0$ and $\Ext^1(M,M) = \Hom(M,M) = k$. We have a long exact sequence in cohomology
$$ 0 \lm \Hom(M,M) \lm \Ext^1(M,I) \lm  \Ext^1(M,\oy) \lm  \Ext^1(M,M) $$
so it suffices to show that $\Ext^1(M,I)=k$. We consider another long exact sequence in cohomology
$$ 0=H^0(I) \lm \Hom(I,I) \lm \Ext^1(M,I) \lm H^1(I) = 0 .$$
Suppose that $Y = \proj A$ and let $I_{\bullet}$ be the ideal of $A$ corresponding to $I$ so that $\Hom(I,I) = \Hom_{\Gr A}(I_{\bullet},I_{\bullet}) \subset Q(A)$ where $Q(A)$ is the field of fractions of $A$. Properness ensures $\Hom(I,I)$ is a finite extension of $k$ which, being a domain must be $k$ itself. This completes the proof.

\vspace{2mm}

We note that the condition $H^1(\oy) = 0$ corresponds in the commutative case to surfaces ruled over $\PP^1$. The commutative proof simplifies in this case since the fibration is constructed from the linear system $|M|$ rather than some multiple of $M$.

\vspace{2mm}
\section{Cohomology of a Mori contraction}  \label{scohmori} 

 In this section, we wish to compute the higher direct images of various sheaves with respect to a non-commutative Mori contraction $f:Y \lm X$. This will be useful in showing that $\nu: \oy \lm f^* \ox$ is an isomorphism in certain cases. We consider first the question of properness of $f$ which, as one expects should follow from properness of $Y$. Hence, in this section, we will only assume $Y$ is a noetherian Ext-finite quasi-scheme.

Recall that a morphism $f:Y\lm X$ of noetherian quasi-scemes is {\em proper} if $R^if_*$ preserves noetherian modules for all $i$. 
\begin{lemma}  \label{lproper}  
Let $f:Y \lm X$ be a flat morphism to a smooth projective curve $X$. Then $R^if_*M$ is noetherian for any $M \in \mo Y$. 
\end{lemma}
\textbf{Proof.} Let $N \in \mo X$. The Ext spectral sequence for $f$ in lemma~\ref{lleray} collapses to give exact sequences 
\[ 0 \lm \Ext^1_X(N,R^if_*M) \lm \Ext^{i+1}_Y(f^*N,M) \lm \Hom_X(N,R^{i+1}f_*M) \lm 0  .\]
The middle term is finite dimensional by Ext-finiteness so the outer terms are too. It suffices to show that for any $F \in \Mod X$ such that $\Ext^1_X(N,F), \Hom_X(N,F)$ are finite dimensional for all $N \in \mo X$, we must have $F$ coherent. Now $\Hom_X(\ox,F)$ is finite dimensional so the torsion part of $F$ is coherent and we may thus henceforth assume that $F$ is torsion-free. 

Suppose $F$ is not coherent and pick any coherent subsheaf $F' < F$. Now $H^0(F/F')$ is also finite dimensional so its torsion subsheaf is coherent. We can thus find a coherent sheaf $F_0 < F$ containing $F'$ such that $F/F_0$ is torsion free. Repeating this procedure allows us to build a strictly increasing chain $F_0 < F_1 < F_2 < \ldots \ $ of coherent subsheaves of $F$ with $F_i$ a sub-bundle of $F_j$ whenever $i<j$. 

Pick a closed point $p \in X$. Now $\calo_p$ is noetherian so using Serre duality we find
\[ \Ext_X^1(\calo_p,F) = \varinjlim \Ext_X^1(\calo_p,F_m) = \varinjlim \Hom_X(F_m,\calo_p)^* = 
\varinjlim F_m \otimes_X \calo_p   .\]
The limit on the right is a union of vector spaces of strictly increasing dimension. This contradiction completes the proof of the lemma. 
%
\vspace{2mm}

 For $Y$ a smooth proper non-commutative surface, the Ext spectral sequence furnishes us with the following vanishing cohomology results.

\begin{prop}  \label{pkillcoh}  
Let $f:Y \lm X$ be a flat morphism to a smooth projective curve $X$. If $Y$ is also smooth of dimension two in the sense that $\Ext^3_Y=0$, then 
\begin{enumerate}
\item for $i \geq 2$ we have $R^if_*=0.$
\item for $M \in \mo Y$ and  $p \in X$ a closed point, we have $R^1f_*M \otimes_X \calo_p \simeq \Ext^2(f^*\calo_p,M)$ so $p$ lies in the support of $R^1f_*M$ iff $\Ext^2(f^*\calo_p,M) \neq 0$.  
\end{enumerate}
\end{prop}
\textbf{Proof.} By lemma~\ref{lproper}, we have $R^if_*M$ is coherent for any $M \in \mo Y$. Serre duality on $X$ shows that for $p \in X$ closed we have 
\[\Ext^1_X(\calo_p,R^if_*M) = \Hom_X(R^if_*M,\calo_p)^* \simeq R^if_*M \otimes_X \calo_p .\]
Now the Ext spectral sequence of lemma~\ref{lleray} for $f$ shows that when $i\geq 2$, the left hand term vanishes so i) follows. Hence it also shows that 
\[ \Ext^2_Y(f^*\calo_p,M)  \simeq \Ext^1_X(\calo_p,R^1f_*M)  \]
which gives us ii).
\vspace{2mm}

We now need to specialise to the case where $Y$ is a non-commutative smooth proper surface and $f:Y \lm X$ is the Fourier-Mukai morphism associated a base point free flat family $\calm/X$ parametrised by a smooth projective curve. 

\begin{prop}  \label{pfOY} 
Let $Y$ be a non-commutative smooth proper surface and $\calm/X$ be a base point free flat family of quotients of $\oy$ parametrised by a smooth projective curve $X$. Let $f:Y \lm X$ be the associated Fourier-Mukai morphism and $r_0,r_1$ be the minimal values of $h^0(M \otimes \w_Y), h^1(M \otimes \w_Y)$ as $M$ varies over the fibres of $\calm$. Then the open set $U$ where the minimal values of $r_0$ and $r_1$ are obtained are the same, and on $U$ we have $f_*\oy,R^1f_*\oy$ are locally free of rank $r_1,r_0$. 
\end{prop}
\textbf{Proof.}
Note first that the fibres of $\calm$ have to have dimension less than 2 for otherwise, since $\oy$ is 2-critical, all the fibres are $\oy$ which contradicts base point freedom. Thus for all $M \in \calm$, we have $H^2(M \otimes\w_Y) \simeq \Hom(M,\oy)^* = 0$ so continuity of Euler characteristic now ensures that $h^0(M \otimes \w_Y),h^1(M \otimes \w_Y)$ minimise on the same open set $U$. 

For a closed point $p \in U$, we compute using proposition~\ref{pkillcoh}
$$R^1f_*\oy \otimes_X \calo_p \simeq \Ext^2_Y(f^*\calo_p,\oy) \simeq H^0(f^*\calo_p \otimes \w_Y)^* \simeq k^{r_0}.$$
This shows $R^1f_*\oy$ is locally free of rank $r_0$ on $U$. 

Similarly, we compute
$$f_* \oy \otimes_X \calo_p \simeq \Ext^1_X(\calo_p,f_* \oy).$$
The Leray spectral sequence gives the exact sequence
$$ 0 \lm \Ext^1_X(\calo_p,f_* \oy) \lm \Ext^1_Y(f^*\calo_p,\oy) \lm \Hom_X(\calo_p,R^1f_*\oy) \lm 0.$$
The right hand term is zero since $R^1f_* \oy$ is locally free at $p$. Hence  BK-Serre duality gives
$$f_* \oy \otimes_X \calo_p \simeq \Ext^1_Y(f^*\calo_p,\oy) \simeq H^1(f^*\calo_p \otimes\w_Y)^*.$$
The last term is $r_1$-dimensional so we are done.
\vspace{2mm}

We next compute $R^if_*f^* \calo_p,R^if_*(f^* \calo_p\otimes \w_Y)$. In the following proposition, $D,D^{\vee}$ will denote ``bad'' sets where we cannot recover the expected commutative behaviour. 

\begin{prop} \label{pRffop} 
Let $Y$ be a non-commutative smooth proper surface, $X$ a smooth projective curve and $f:Y \lm X$ the Fourier-Mukai morphism associated to a base point free Hilbert system. 
\begin{enumerate}
\item Let $D \subset X$ be the closed set of points $p \in X$ where $h^0(f^*\calo_p)> 1$ and $D^{\vee}\subseteq X$ be the closed set where $h^0(f^*\calo_p \otimes \w_Y) > 0$.
\begin{enumerate}
\item For $p \in X$ we have $\supp R^1f_*f^*\calo_p \subseteq D^{\vee}$.
\item For $p \in X - D$, we have $R^1f_*(f^*\calo_p \otimes \w_Y) = \calo_p \oplus \mathcal{F}$ where $\supp \mathcal{F} \subseteq D$. 
\end{enumerate}
\item Suppose now that for some (and hence every) $p \in X$ we have $f^*\calo_p.f^*\calo_p = 0$. Then 
\begin{enumerate}
\item For $p \in X - D$, we have $f_*f^* \calo_p = \calo_p \oplus \mathcal{F}$  where $\supp \mathcal{F} \subseteq D^{\vee}$.
\item For $p \in X - D - D^{\vee}$, we have $\supp f_*(f^*\calo_p \otimes \w_Y) \subseteq D$.
\end{enumerate}
\end{enumerate} 
\end{prop}
\textbf{Proof.} For any $p\in X, q \in X-D^{\vee}$, we have $\Ext^2_Y(f^*\calo_q,f^*\calo_p) = 0$ so part i)a) follows from proposition~\ref{pkillcoh}. 

We now prove i)b). Let $q \in X-D$ be a closed point. From proposition~\ref{pkillcoh}ii), we have 
$$R^1f_*(f^*\calo_p \otimes \w_Y) \otimes\calo_q \simeq \Ext^2_Y(f^*\calo_q,f^*\calo_p \otimes\w_Y) \simeq \Hom_Y(f^*\calo_p,f^*\calo_q)^*.$$
Since $\calm/X$ is a Hilbert system, proposition~\ref{pextmmp}iii) tells us this is $k$ if $p=q$ and 0 otherwise. Hence $R^1f_*(f^*\calo_p \otimes \w_Y)= \calo_{mp} \oplus \mathcal{F}$ for some positive integer $m$ and sheaf $\mathcal{F}$ supported in $D$. We need to show that $m=1$ which will follow from showing 
$$ \Hom_X(R^1f_*(f^*\calo_p \otimes \w_Y),\calo_{2p}) \simeq k .$$
We use the Leray spectral sequence as in the proof of proposition~\ref{pkillcoh}ii). This time it gives 
$$\Hom_X(R^1f_*(f^*\calo_p \otimes \w_Y),\calo_{2p}) \simeq \Ext^2_Y(f^*\calo_{2p},f^*\calo_p \otimes\w_Y)^* \simeq \Hom_Y(f^*\calo_p,f^*\calo_{2p}).$$
To compute this last term, recall that the theory of Hilbert schemes (proposition~\ref{p1dimhilbsys}) tells us there is a non-split extension
$$ (*) \hspace{2cm} 0 \lm f^*\calo_p \lm f^*\calo_{2p}  \lm f^*\calo_p  \lm 0.$$
In the long exact sequence in cohomology obtained by applying $\Hom_Y(f^*\calo_p,-)$, the connecting homomorphism $\Hom_Y(f^*\calo_p,f^*\calo_p) \lm \Ext^1_Y(f^*\calo_p,f^*\calo_p)$ is non-zero. This ensures that $\Hom_Y(f^*\calo_p,f^*\calo_{2p})=k$, thus proving part i)b).

For part ii), we will use the ext computations of proposition~\ref{pextmmp}. To prove part ii)a), note that $\Ext^1_Y(f^*\calo_q,f^*\calo_p) = 0$ for $p \neq q \in X -D^{\vee}$ so the Leray spectral sequence then gives $\Ext^1_X(\calo_q,f_*f^*\calo_p) = 0$. Thus $f_*f^*\calo_p = \calo_{np} \oplus \mathcal{F}$ for some integer $n >0$ and $\mathcal{F}$ supported in $D^{\vee}$. To show $n=1$, we need to show $\Hom_Y(f^*\calo_{2p},f^*\calo_p) = k$. As in the proof of part i)b), this follows on applying $\Rhom_X(-,f^*\calo_p)$ to the non-split exact sequence (*) above. This proves part ii)a). 

For part ii)b), we compute $f_*(f^*\calo_p\otimes\w_Y) \otimes_X\calo_q \simeq \Ext^1_X(\calo_q,f_*(f^*\calo_p\otimes\w_Y))$ for $q \in X-D$. The Leray sequence gives us this time the exact sequence 
$$0 \lm \Ext^1_X(\calo_q,f_*(f^*\calo_p\otimes\w_Y)) \lm \Ext^1_Y(f^*\calo_q,f^*\calo_p\otimes\w_Y) \lm 
\Hom_X(\calo_q,R^1f_*(f^*\calo_p\otimes\w_Y))\lm 0.$$
Part ii)b) will be proved as soon as we can show that the middle and last terms of the sequence have the same dimension. Now part i)b) shows that the last term is $\Hom_X(\calo_q,\calo_p)$ which is $k$ if $p=q$ and 0 otherwise. The same is true of $\Ext^1_Y(f^*\calo_q,f^*\calo_p\otimes\w_Y) \simeq \Ext^1_Y(f^*\calo_p,f^*\calo_q)^*$ by proposition~\ref{pextmmp}. This completes the proof of part ii)b) and the proposition.



\vspace{2mm}

Unfortunately, we can only hypothesise away the bad sets in the previous proposition. 
\begin{defn}  \label{duniform} 
We say  that a flat family of $Y$-modules $\calm/X$ is {\em uniform} if for every $M \in \calm$, i) $h^0(M)=1$  and ii) $h^0(M \otimes \w_Y) = 0$. Suppose that $\calm/X$ is base point free so induces a Fourier-Mukai morphism $f:Y \lm X$. We say that $f$ is {\em uniform} if $\calm/X$ is. 
\end{defn}
If $\calm$ is a flat family of quotients of $\oy$, then semi-continuity implies that we always have $h^0(M)\geq 1$  and  $h^0(M \otimes \w_Y) \geq 0$ so i) and ii) assert equality here. Hypothesis ii) in the commutative case follows from the fact that one should assume $K$-negativity i.e. $M.K <0$ and intersection products are continuous. However, as we have already remarked, we do not know if $K$-negative implies $K$-non-effective. 

\begin{prop}  \label{pXsmooth} 
Let $\calm/X$ be the Hilbert system of a $K$-non-effective rational curve with self-intersection zero. Suppose that  $\calm/X$ is uniform. Then i) every $M \in \calm$ is a $K$-non-effective rational curve with self-intersection zero and ii) the curve $X$ is smooth.
\end{prop}
\textbf{Proof.} 
Part i) follows from continuity of Euler characteristics and the fact that $H^2$ vanishes in this case. Part ii) follows from corollary~\ref{cfamily}.

\begin{lemma}  \label{lffox} 
Let $f:Y \lm X$ be a uniform non-commutative Mori contraction. Then $R^1f_*f^*\ox = 0$ and $f_*f^*\ox = \ox$.
\end{lemma}
\textbf{Proof.}
For any $p \in X$ closed we have by propositions~\ref{pkillcoh},\ref{pRffop},
$$R^1f_*f^*\ox \otimes_X \calo_p \simeq \Ext^2_Y(f^*\calo_p,f^*\ox) \simeq \Hom_Y(f^*\ox,f^*\calo_p \otimes \w_Y)^* \simeq \Hom_X(\ox, f_*(f^*\calo_p\otimes\w_Y))^* = 0 .$$
Thus $R^1f_*f^*\ox = 0$.

We next show $\dim_k (f_*f^*\ox \otimes_X \calo_p) = 1$ for all $p \in X$ so that $f_*f^*\ox$ must be locally free of rank 1. Using $R^1f_*f^*\ox = 0$, the Leray spectral sequence this time gives an isomorphism 
$$\Ext^1_X(\calo_p,f_*f^*\ox) \simeq \Ext^1_Y(f^*\calo_p,f^*\ox) \simeq \Ext^1_Y(f^*\ox,f^*\calo_p\otimes\w_Y)^*$$
so it suffices to show the last term is 1-dimensional. This follows from the exact sequence
$$ 0 \lm \Ext^1_X(\ox,f_*(f^*\calo_p \otimes \w_Y)) \lm \Ext^1_Y(f^*\ox,f^*\calo_p \otimes \w_Y) \lm \Hom_X(\ox,R^1f_*(f^*\calo_p \otimes \w_Y)) \lm 0 $$
and proposition~\ref{pRffop}. Thus $f_*f^*\ox$ is locally free of rank one. 

Adjunction gives a non-zero map $\ox \lm f_*f^*\ox$ and hence the commutative diagram below on the left where $\xi:\ox(D) \lm f_*f^*\ox$ is an isomorphism for some effective divisor $D$. 
$$\diagram
\ox \dto_{\eta}  \rto & f_*f^*\ox \\
\ox(D) \urto_{\xi} & 
\enddiagram\ \ \ 
\diagram
f^*\ox \dto_{f^*\eta}  \rto^{\id} & f^*\ox \\
f^*\ox(D) \urto & 
\enddiagram$$
Adjunction also gives us the commutative diagram on the right which shows $f^*\eta$ is split injective. It suffices to show that $D=0$ so suppose to the contrary that $p \in D$. Pushing forward and pulling back the extension
$$ 0 \lm f^*\ox \lm f^*\ox(D) \lm f^*\calo_D \lm 0 $$
by $f^*\ox \lm f^*\calo_p$ and $f^*\calo_p \hookrightarrow f^*\calo_D$ we see that the extension
$$ 0 \lm  f^* \calo_p \lm f^* \calo_{2p} \lm f^* \calo_p  \lm 0 $$
must also split. This contradicts proposition~\ref{p1dimhilbsys} and completes the proof of the lemma.

\vspace{2mm}

If $f^*\ox = \oy$ then the lemma tells us in particular that $f_*\oy = \ox$. In the commutative case, this would be immediate from the fact that $f$ is its own Stein factorisation since, for example, $X$ is normal.

\begin{cor} \label{ch1Mw} 
Let $f:Y \lm X$ be a uniform non-commutative Mori contraction. If the natural map $\nu: \oy \lm f^*\ox$ is injective, then for generic $p \in X$, we must have $h^1(f^*\calo_p \otimes \w_Y) \leq 1$. 
\end{cor}
\textbf{Proof.}
Note that $f_*$ is left exact so $f_*\oy$ embeds in $f_*f^*\ox$. But $f_*f^*\ox$ is locally free of rank 1 so proposition~\ref{pfOY} finishes the proof.


\vspace{2mm}
\section{Sufficient criteria for $\oy \simeq f^* \ox$}    \label{snuonto}  

In this section, $Y$ will be a non-commutative smooth proper surface with finitely partitive dimension function. We consider a non-commutative Mori contraction $f:Y \lm X$ and the associated natural map $\nu:\oy \lm f^* \ox$. Our goal will be to extend the injectivity of $\nu$ results obtained in section~\ref{sdisjoint}, to sufficient criteria for $\nu$ to be an isomorphism. Ideally, one would like a criterion which depends only on the contracted $K$-non-effective rational curve $M$ with self-intersection zero. For example, theorem~\ref{tptfibre} provides such a criterion for injectivity of $\nu$, at least when the ground field is uncountable. Unfortunately, in this part we will need to assume  at the very least that $f$ is uniform, a condition which presumably can only be checked if you actually have a complete 1-parameter family of deformations of $M$. 

The simplest criterion to guarantee $\oy \simeq f^*\ox$ is the following result.

\begin{prop}  \label{poysurjects} 
Let $f:Y \lm X$ be a uniform non-commutative Mori contraction of a curve $M$ with $H^1(M \otimes \w_Y) \neq 0$. 
\begin{enumerate}
\item Suppose the natural map $\nu:\oy \lm f^*\ox$ is injective.
\begin{enumerate}
\item There exists some effective divisor $D \subset X$ such that $f_*\oy = \ox(-D)$.
\item The natural map $f^*\ox(-D) \hookrightarrow f^*\ox$  factors through $\nu$, so in particular, $C:=\coker \nu$ is a quotient of $f^*\calo_D$. 
\item $\nu$ is an isomorphism if and only if $D = 0$. 
\end{enumerate}
\item Suppose that for every closed subscheme $D \subsetneq X$ we have that the natural map $\oy \lm f^*\calo_D$ is surjective. Then $\nu$ is an isomorphism. 
\end{enumerate}
\end{prop}
\textbf{Proof.} 
Assuming i), we have an exact sequence
$$ 0 \lm \oy \lm f^* \ox \lm C  \lm 0 .$$
Now proposition~\ref{pfOY} and our uniform assumption ensures $f_*\oy \neq 0, R^1f_* \oy = 0$ so we have an exact sequence
$$0 \lm f_*\oy \lm \ox \xrightarrow{\gamma} f_*C \lm 0 $$
by lemma~\ref{lffox}. This identifies $f_*\oy$ with $\ox(-D)$ for some effective divisor $D \subset X$ and $f_*C$ with $\calo_D$. Part c) holds since $\nu$ is an isomorphism if and only if $\gamma$, or equivalently $D$ is zero. Part b) follows from the commutative diagram below.
$$\diagram
\oy \rto^{\nu} & f^*\ox \dto \rto^{\phi} & C \\
& f^*\calo_D \urto &  
\enddiagram$$

We now prove part ii). We already know from proposition~\ref{poyinject} that $\nu$ is injective so we need to show $C=0$. Now $\oy \xrightarrow{a} f^*\calo_D \xrightarrow{b} C$ is a factorisation of the zero map. By assumption, $a$ is surjective while $b$ is surjective since $\phi:f^*\ox \lm C$ is. This completes the proof of the proposition.



\vspace{2mm}

We need some more hypotheses. 
\begin{defn}  \label{dextrem} 
Let $M$ be a $K$-non-effective rational curve $M$ with $M^2=0$ and $\calm/X$ be its Hilbert system. We say that $M$ is {\em extremal} if $h^1(M \otimes \w_Y) = 1$ and for every $M' \in \calm$ we have $M'$ is 1-critical.
\end{defn}
To explain these, recall that in the commutative setting, extremal curves cannot be algebraically deformed into two curves (that is, ``bent and broken''). This corresponds roughly to the hypothesis that all fibres are 1-critical. The condition $h^1(M \otimes \w_Y) = 1$ does not have anything to do with the extremal condition but is automatic in the commutative case by the genus formula. 

We say that a non-commutative Mori contraction $f:Y \lm X$  is {\em extremal}, if it contracts an extremal $K$-non-effective rational curve with self-intersection zero.

\begin{thm}  \label{tontofstartod} 
Let $f:Y \lm X$ be an extremal uniform non-commutative Mori contraction of the curve $M$. Suppose furthermore that any point $P$ on a fibre of $f$ satisfies $h^0(P)=1$ (so all fibres are ``fat-free''). Then $\nu:\oy \lm f^*\ox$ is an isomorphism. 
\end{thm}
\textbf{Proof.}
We verify the criterion of proposition~\ref{poysurjects}. We start with
\begin{lemma}
Let $p\in X$ be a closed point and $m$ be an integer. Then any quotient map $\phi:f^*\calo_{mp} \lm P$ where $P$ is a point of $Y$, factors through the natural surjection $f^*\calo_{mp} \lm f^*\calo_p$. 
\end{lemma}
\textbf{Proof.}
We argue by induction on $m$, the case $m=1$ being clear. Consider the exact sequence
$$ 0 \lm M \lm f^*\calo_{mp} \lm f^*\calo_{(m-1)p} \lm 0 .$$
By induction, it suffices to show that $\phi|_M = 0$. Restricting $\phi$ to 
$$\ker(f^*\calo_{mp} \lm f^*\calo_{(m-2)p}) \simeq f^*\calo_{2p}$$
we see that we need only prove the lemma for the case $m=2$. 

Note that $M \simeq f^*\calo_p$ so letting $M_2 = f^*\calo_{2p}$, we may re-write the sequence above as
$$ E:0 \lm M \lm M_2 \lm M  \lm 0$$
which we note is non-split. Also, $P$ must be a point of $M$ so has $h^0(P)=1$. Now $\phi|_M$ maps to zero under the connecting homomorphism $\Hom_Y(M,P) \lm \Ext^1_Y(M,P)$. Thus the extension $E$ lies in
$$\ker ( \Ext^1_Y(M,M) \xrightarrow{\Ext^1_Y(M,\phi|_M)} \Ext^1_Y(M,P))$$
In particular, if $J:= \ker \phi|_M$, then we see that $\Ext^1(M,J) \neq 0$. 

Base point freedom ensures that $\xi(M,P) \leq 0$ so 
$$\xi(M,J) = \xi(M,M) - \xi(M,P) \geq 0.$$
Now $\Hom_Y(M,J) = 0$ so we must have 
$$0 \neq \Ext^2_Y(M,J) \simeq \Hom_Y(J,M \otimes\w_Y)^*.$$
This means there is a non-zero map $J \lm M \otimes\w_Y$ which, since $J$ and $M\otimes\w_Y$ are 1-critical, must be an injection with 0-dimensional cokernel. Thus
$$\chi(M)-\chi(P) = \chi(J) \leq \chi(M \otimes\w_Y)$$
which forces $h^0(P) > 1$, a contradiction. This proves the lemma.

\vspace{2mm}

We return now to the proof of the theorem. Suppose that $D = m_1p_1 + \ldots +m_jp_j$ where $p_1,\ldots,p_j$ are distinct points of $X$. We show by induction on $j$ that $\oy\lm f^*\calo_D$ is surjective. The case $j=1$ follows from the lemma and the fact that $\oy \lm f^*\calo_p$ is surjective. We may assume that we have an exact sequence
$$0 \lm K \lm \oy \lm \oplus_{i=1}^{j-1} f^*\calo_{m_ip_i} \lm 0 .$$
If the composite $\psi:K\lm \oy \lm f^*\calo_{m_jp_j}$ is surjective then we are done. Suppose this is not the case and that $\psi$ has a non-zero cokernel $C$. Pick a simple quotient $C \lm P$ which must correspond to a point on $f^*\calo_{p_j}$ by the lemma. The composite map $\oy \lm C \lm P$ must also factor through $\oplus_{i=1}^{j-1} f^*\calo_{m_ip_i}$. This is impossible as distinct fibres of $f$ have distinct simple quotients by proposition~\ref{ph0P1}. This proves the theorem.
\vspace{2mm}

The next result states that if $\nu$ is not an isomorphism, then one of two perverse phenomena occur.

\begin{thm}  \label{tnunotiso} 
Let $f:Y \lm X$ be an extremal uniform non-commutative Mori contraction of the curve $M$. Suppose that $\nu:\oy \lm f^*\ox$ is not an isomorphism. Then there exists $p \in X$ such that one of the following must occur.
\begin{enumerate}
\item $\Ext^1(f^*\calo_p,\oy) \lm \Ext^1(f^*\calo_p,f^*\calo_p)$ is zero.
\item $f^*\calo_{2p}$ is not uniform in the sense that it contains the direct sum of two non-zero submodules. 
\end{enumerate}
\end{thm}
\textbf{Proof.}
Suppose condition i) does not hold so we need to show condition ii) holds. By scholium~\ref{schptfibre}, we know for any reduced divisor $D$, the map $\oy \lm f^*\calo_D$ is surjective. We also know $\nu$ is injective so we may apply the results of proposition~\ref{poysurjects}. Consider the exact sequence
$$ 0 \lm \oy \lm f^* \ox \lm C  \lm 0 $$
as usual. The sequence is non-split since it is non-split when you apply $f_*$ to it. Our classical cohomology assumption ensures that $\dim C = 1$. Let $f_*\oy = \ox(-D)$ as in proposition~\ref{poysurjects}. We consider $C_B:= \coker \oy \lm f^*\calo_B$ for $D_{red} \leq B \leq D$. When $B = D_{red}$ we see $C_B = 0$ whilst when $B = D$ we get $C_B = C$. Hence we can find a divisor $B$ and point $p\in X$ with $D_{red}\leq B  \leq D-p$ and $\dim C_B \leq 0$ but $\dim C_{B + p} = 1$. Consider the following commutative diagram with exact rows.
$$\diagram
0 \rto & \rto M & f^*\calo_{B+p} \rto^{\pi} \dto & f^*\calo_B \rto \dto & 0 \\
       &        & C_{B+p}        \rto      & C_B        \rto      & 0
\enddiagram$$
where $M \simeq f^*\calo_p$. Let $K = \ker (f^*\calo_{B+p} \lm C_{B+p})$. Now $\ker(C_{B+p} \lm C_B) \simeq M/M\cap K$ must be 1-dimensional so, since $M$ is 1-critical, $K \cap M = 0$. Also, $B \geq D_{red}$ so there exists a submodule $N \leq f^*\calo_B$ which is isomorphic to $f^*\calo_p$ and such that $\pi^{-1}(N) \simeq f^*\calo_{2p}$. Now $\dim C_B \leq 0$ means that we must have $\pi(K) \cap N \neq 0$. Hence $\pi^{-1}(N) \cap K$ and $M$ are two non-trivial submodules of $\pi^{-1}(N)$ which intersect trivially. This shows $f^*\calo_{2p}$ is not uniform.

\vspace{2mm} Note that by proposition~\ref{ph1oyis0}, condition i) is excluded if $Y$ is projective and $H^1(\oy) = 0$.

\vspace{5mm}

\textbf{\large References}

\begin{itemize}
\item [{[A]}] M. Artin, ``Some problems on three-dimensional graded domains '', Representation theory and algebraic geometry, London Math. Soc. Lecture Notes Series \textbf{238}, Cambridge Univ. Press (1995), p.1-19
\item [{[ATV]}] M. Artin, J. Tate, M. Van den Bergh, ``Some algebras associated to automorphisms of elliptic curves'', The Grothendieck Festschrift vol. 1 p.333-85, Prog. Math. \textbf{86}, Birkhauser Boston, (1990)
  \item  [{[AV]}] M. Artin, M. Van den Bergh,  ``Twisted Homogeneous Coordinate Rings'', Journal of Algebra \textbf{133} (1990) p.249-71
       {\em J. of Algebra}, \textbf{133}, (1990) , p.249-271
\item [{[AZ94]}] M. Artin, J. Zhang, ``Noncommutative projective schemes'' Adv. Math. \textbf{109} (1994), p.228-87
  \item  [{[AZ01]}] M. Artin, J. Zhang,  ``Abstract Hilbert Schemes'', Alg. and Repr. Theory \textbf{4} (2001), p.305-94
\item [{[BK]}] A. Bondal, Kapranov, ``Representable functors, Serre functors and mutations'', Math USSR-Izv.  \textbf{35} (1990) p.519-41
\item [{[CI]}] D. Chan, C. Ingalls, ``Minimal model program for orders over surfaces'' Invent. Math. \textbf{161} (2005) p.427-52
\item [{[Hart]}] R. Hartshorne, ``Algebraic geometry'', Graduate Texts in Mathematics \textbf{52}, Springer Verlag, New York, (1977)
\item [{[KL]}]  G. Krause, T. Lenagan, ``Growth of algebras and Gelfand-Kirillov dimension'' Graduate Studies in Math. \textbf{22}, American Mathematical Society, Providence (2000)
\item [{[KM]}] J. Kollar, S. Mori, ``Birational geometry of algebraic varieties'', Cambridge Tracts in Math. \textbf{134} Cambridge Univ. Press (1998)
\item [{[McR]}] J. McConnell, J. Robson, ``Noncommutative noetherian rings'' Graduate Studies in Math. \textbf{30} American Mathematical Society, Providence (2001)
  \item  [{[Mori]}]  I. Mori, ``Intersection theory over quantum ruled surfaces'', J. Pure \& Applied Algebra, \textbf{211}, (2007), p.25-41
\item [{[MS01]}]  I. Mori, P. Smith, ``B\'{e}zout's theorem for non-commutative projective spaces'' J. Pure \& Applied Algebra \textbf{157} (2001) p.279-99
\item [{[MS06]}]  I. Mori, P. Smith, ``The Grothendieck group of a projective space bundle'' $K$-theory \textbf{37} (2006) p.263-89
\item  [{[Na]}] A. Nyman, ``Serre duality for non-commutative $\PP^1$-bundles'', Trans. AMS, 
            \textbf{357}, (2004), p.1349-416
\item [{[Nb]}] A. Nyman, ``Serre finiteness and Serre vanishing doe non-commutative $\PP^1$-bundles'', J. of Algebra \textbf{278} (2004), p.32-42
\item [{[NV04]}] K. de Naeghel, M. Van den Bergh, ``Ideal classes of three dimensional Sklyanin algebras'' J. Algebra \textbf{276} (2004), p.515-51
\item [{[NV05]}] K. de Naeghel, M. Van den Bergh, ``Ideal classes of three dimensional Artin-Schelter regular algebras'' J. Algebra \textbf{283} (2005), p.399-429
\item [{[Pat]}] D. Patrick, ``Noncommutative ruled surfaces'' PhD thesis MIT, June 2007
\item [{[SV]}] S. P. Smith, M. Van den Bergh, ``Non-commutative quadrics'' preprint
\item  [{[VdB96]}]   M. Van den Bergh, ``A Translation Principle for the
          Four-dimensional Sklyanin Algebras'', J. of Algebra,
          \textbf{184}, (1996), p.435-490
\item [{[VdB97]}] M. Van den Bergh, ``Existence theorems for dualizing complexes over non-commutative graded and filtered rings'' J. Algebra \textbf{195} (1997), p.662-79
\item [{[VdB01]}] M. Van den Bergh, `` Blowing up of non-commutative surfaces'' Memoirs AMS \textbf{159} (2001)
  \item  [{[VdB01p]}]   M. Van den Bergh, ``Non-commutative
          $\mathbb{P}^1$-bundles over Commutative Schemes'',
          math.RA/0102005 1 Feb. 2001
\item [{[Vistoli]}] A. Vistoli, ``Intersection theory on algebraic stacks and their moduli spaces'', Invent. Math. \textbf{97} (1989), p.613-70
  \item  [{[YZ]}] A. Yekutieli, J. J. Zhang, ``Rings with Auslander dualizing complexes'' J. Algebra \textbf{213} (1999), p.1-51
\end{itemize}

\end{document}